%% file: GEOIIpp2.tex
\input zineb.tex

 \auteurcourant={ M. Kaddar }
 \pageno=1
\baselineskip=5mm

 \hoffset=-4mm
  \voffset=8mm
 \hsize=1cm
\vsize=180mm \hfuzz=-2cm
\input amssym.def
\input amssym.tex%
\magnification = 1200 \hoffset=-4mm
  \voffset=8mm
 \hsize=14cm
\vsize=190mm
\input xy
\xyoption{all}
\par
\topskip = 10 truemm \overfullrule=0pt \rightskip=10mm

\def \cqfd{\unskip\kern 6pt\penalty 500
 \raise -2pt\hbox{\vrule\vbox to10pt{\hrule width 4pt
\vfill\hrule}\vrule}\par}

\def\adots{\mathinner{\mkern2mu\raise1pt\hbox{.}
\mkern3mu\raise4pt\hbox{.}\mkern1mu\raise7pt\hbox{.}}}
\def\pmb#1{\setbox0=\hbox{#1}%
\kern-.025em\copy0\kern-\wd0 \kern.05em\copy0\kern-\wd0
\kern-.025em\raise.0433em\box0}

 \font \tite=cmbxti10 scaled 1400


\def\build#1_#2^#3{\mathrel{\mathop{\kern 0pt#1}\limits_{#2}^{#3}}}

\def\build#1_#2^#3{\mathrel{\mathop{\kern 0pt#1}\limits_{#2}^{#3}}}
\def\hfl#1#2{\smash{\mathop{\hbox to 18mm{\rightarrowfill}}\limits^{\scriptstyle#1}_{
\scriptstyle#2}}}

\def\boxit#1#2{\setbox1=\hbox{\kern#1{#2}\kern#1}%
\dimen1=\ht1 \advance\dimen1 by #1 \dimen2=\dp1 \advance\dimen2 by
#1
\setbox1=\hbox{\vrule height\dimen1 depth\dimen2\box1\vrule}%
\setbox1=\vbox{\hrule\box1\hrule}%
\advance\dimen1 by .4pt \ht1=\dimen1 \advance\dimen2 by .4pt
\dp1=\dimen2 \box1\relax} \titrecourant={Morphismes
g\'eom\'etriquement plats  et classes fondamentales relatives en
g\'eom\'etrie analytique complexe.}

\input xy

\headline={\ifnum\pageno=1 {\hfill} \else{\hss \tenrm -- \folio\ --
\hss} \fi} \footline={\hfil}



\input xy
\xyoption{all}
\par
\par\centerline{{\tite{PLATITUDE GEOMETRIQUE ET }}}\bigskip
\centerline{\tite{ CLASSES FONDAMENTALES
RELATIVES}}\bigskip\noindent\centerline{\tite{ PONDEREES
II.}}\bigskip\bigskip\bigskip
\centerline{MOHAMED  KADDAR}\bigskip\bigskip\bigskip\bigskip
\noindent $${\vbox{\hsize=12cm{{\bf Abstract.~}\sevenrm{ This second
part is devoted to the proof of all main results that we have
mentionned in [KI].}}}}$$
\bigskip\bigskip\noindent
AMS Classification (2000): 32C15, 32C30, 32C35, 32C37,
14C05.\bigskip\noindent Key words: Analytic spaces, Integration,
cohomology, dualizing sheaves.\bigskip\bigskip\noindent
\centerline{{\tite Table des mati\`eres:}}\bigskip\noindent \S 0.
{\bf Introduction et \'enonc\'es des
r\'esultats}\page{1}\smallskip\noindent \S 1. {\bf La preuve du
th\'eor\`eme 1}\par
 \leftskip=1em
{\bf 1.0. Quelques annulations de cohomologie.}\page{6}\par {\bf
1.1. Preuve du th\'eor\`eme 1.}\page{10}\par {\bf 1.2.
 Preuve du corollaire 1.1.}\page{15}
\smallskip\noindent \S 2. {\bf La preuve du th\'eor\`eme
2} \par
 \leftskip=1em
{\bf 2.0.Platitude g\'eom\'etrique analytique et int\'egration
globale.}\page{16}\par {\bf 2.1.Platitude g\'eom\'etrique analytique
et morphisme classe fonda\par mentale relative
pond\'er\'e.}\page{38}\smallskip\noindent \S 3. {\bf La preuve des
th\'eor\`emes 3 et 4.}\par
 \leftskip=1em
{\bf 3.1. Preuve du th\'eor\`eme 3.}\page{55}\par {\bf 3.2. Preuve
du corollaire 3.1.}\page{57}\par {\bf 3.3. Platitude g\'eom\'etrique
analytique, faisceaux dualisants et formes\par m\'eromorphes
r\'eguli\`eres: le th\'eor\`eme 4.}\page{58}
\smallskip\noindent
\S 4.  {\bf Bibliographie}\page{68}
\bigskip\bigskip\noindent
 \$ 0. Dans cette partie, qui
fait suite \`a [KI], nous allons donner la preuve
 des r\'esultats qui y \'etaient annonc\'es. Avant cela,
 rappelons quelques notations introduites dans [KI]. Soient  $n$  un entier naturel
  donn\'e et $S$  un espace analytique complexe r\'eduit de dimension
  (localement) finie. Un morphisme $\pi:X\rightarrow S$ d'espaces complexes est dit
  {\it universellement}-$n$-{\it \'equidimensionnel} s'il est ouvert et \`a fibres de
   dimension pure constante $n$. Nous avons d\'esign\'e par:\par\noindent
   $\bullet$  ${\cal E}(S,n)$ l'ensemble de tous
    les morphismes d'espaces analytiques complexes $\pi:X\rightarrow
    S$ universellement $n$-\'equidimensionnels,\par\noindent
    $\bullet$ ${\cal E}_{pond}(S,n)$ l'ensemble des morphismes
    universellement $n$-\'equidimensionnels munis d'une certaine
    pond\'eration $\goth{X}$,\par\noindent
    $\bullet$ ${\cal G}_{c}(S,n)$ (resp. ${\cal G}_{a}(S,n)$) le
    sous ensemble de ${\cal E}_{pond}(S,n)$ constitu\'e
    d'\'el\'ements dont la pond\'eration est contin\^ument
    g\'eom\'etriquement plate (resp analytiquement
    g\'eom\'etriquement plate). \par\noindent
    On a (cf [KI]) les inclusions g\'en\'eralement strictes\par
    \centerline{${\cal G}_{a}(S,n)\subset {\cal G}_{c}(S,n) \subset
    {\cal E}_{pond}(S,n)\subset {\cal E}(S,n)$}\par\noindent
On supposera $X$ d\'enombrable \`a l'infini (ou paracompact et
compl\`etement paracompact quand il sera n\'ecessaire de le
pr\'eciser). Par ailleurs, il est plus commode de travailler avec
$S$ de dimension pure (ce qui entrainera automatiquement que $X$
l'est aussi puisque le morphisme est ouvert et \`a fibres de
dimension pure constante) car le  cas g\'en\'eral n'apporte que
tr\`es peu de substance au texte si ce n'est la lourdeur des
notations et les discussions techniques sans grand
inter\^et.\smallskip\noindent
On divisera le texte en trois parties I, II et III consacr\'ees respectivement aux preuves des th\'eor\`emes 1, 2 et 3  dont nous allons rappeler les \'enonc\'es pour plus de commodit\'e.\smallskip\noindent
 \Th{1}{} Soit $\pi$ un \'el\'ement de
  ${\cal E}(S,n)$.  Alors, il existe un unique faisceau de ${\cal O}_{X}$-modules
${\omega}^{n}_{\pi}$, v\'erifiant :\par\noindent {\bf (i)} Si
$\xymatrix{X\ar@/_/[rr]_{\pi}\ar[r]^{\sigma}&Z\ar[r]^{q}&S}$ est une
factorisation locale de $\pi$
 dans laquelle  $\sigma$ est un plongement local dans $Z$ lisse sur $S$ et
 de dimension relative $n+p$, alors\par
\centerline{$\omega^{n}_{\pi}\simeq \sigma^{*}{\cal E}xt^{p}({\cal
O}_{X}, \Omega^{n+p}_{Z/S})$}\par\noindent {\bf (ii)} il est ${\cal O}_{X}$
coh\'erent, de profondeur au moins deux fibres par fibres sur
$S$\par\noindent {\bf(iii)} $\omega^{n}_{\pi}$ et $\Omega^{n}_{X/S}$
coincident canoniquement sur  la partie r\'eguli\`ere du morphisme
$\pi$.\smallskip\noindent {\bf (iv) } sa construction est compatible
aux inclusions ouvertes sur $X$ dans le sens de Verdier :\par\noindent
 si   $\pi_{i}:X_{i}\rightarrow S,\,\,i=1,2$ sont deux morphismes analytiques  universellement $n$-\'equidimensionnels  et $U$ un ouvert de $X_{1}$ muni
  de deux inclusions ouvertes $j_{i}:U\rightarrow X_{i}$
  tels  que le diagramme
  $$
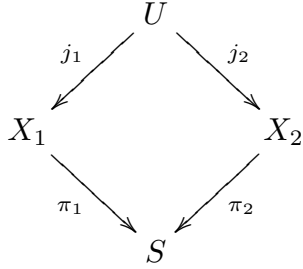
\xymatrix{&U\ar[ld]_{j_{1}}\ar[rd]^{j_{2}}&\\
  X_{1}\ar[rd]_{\pi_{1}}&&X_{2}\ar[ld]^{\pi_{2}}\\
  &S&}$$
  soit commutatif, on a $\ds{j_{1}^{*}(\omega^{n}_{\pi_{1}}) =
  j_{2}^{*}(\omega^{n}_{\pi_{2}})}$
et donc, en particulier,  pour tout ouvert  $U$ de $X$ muni de
l'injection naturelle $j:U\rightarrow X$ et de la restriction de
$\pi$ \`a $U$ que l'on note $\pi_{U}$, on a $\omega^{n}_{\pi}|_{U}=\omega^{n}_{\pi_{U}}$. De plus, il est  stable par tout changement de base plat entre espaces complexes.\par\noindent
{\bf (v)}  il munit $\pi$ d'un morphisme canonique $\ds{{\rm I}\!{\rm
R}^{n}\pi_{!}\omega^{n}_{\pi}\rightarrow {\cal O}_{S}}$ de formation
compatible aux restrictions ouvertes sur $X$ (resp. $S$) et aux changements de bases plats sur $S$.\par\noindent
{\bf (vi) } si $\pi$ est propre, $\omega^{n}_{\pi}= {\cal H}^{-n}(\pi^{!}({\cal O}_{S}))$\par\noindent
{\bf (vii)}
Tout  diagramme commutatif d'espaces analytiques complexes
$$
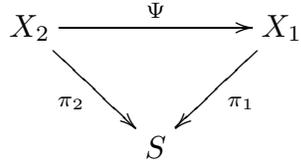
\xymatrix{X_{2}\ar[rr]^{\Psi}\ar[rd]_{\pi_{2}}&&X_{1}\ar[ld]^{\pi_{1}}\\
&S&}$$
 avec $\pi_{1}$ (resp. $\pi_{2}$)  universellement \'equidimensionnels de dimension relative $n_{1}$ (resp. $n_{2}$)  et $\Psi$ propre de dimension relative born\'ee par l'entier   $d:=n_{2}-n_{1}$,
 donne un morphisme canonique ${\rm I}\!{\rm
 R}^{d}{\Psi}_{*}\omega^{n_{2}}_{\pi_{2}}\rightarrow \omega^{n_{1}}_{\pi_{1}}$ induisant le diagramme commutatif de faisceaux analytiques
$$
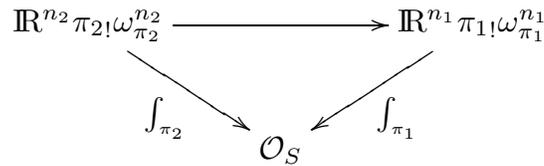
\xymatrix{{\rm I}\!{\rm
 R}^{n_{2}}{\pi_{2}}_{!}\omega^{n_{2}}_{\pi_{2}}\ar[rr]\ar[rd]_{\int_{\pi_{2}}}&&
 {\rm I}\!{\rm
 R}^{n_{1}}{\pi_{1}}_{!}\omega^{n_{1}}_{\pi_{1}}\ar[ld]^{\int_{\pi_{1}}}\\
&{\cal O}_{S}&}$$
{\bf (viii) }  si $S$ est un point, il coincide avec  le faisceau
dualisant de Grothendieck, Andr\'eotti- Kas ou Golovin :\par
\centerline{$\omega^{n}_{\pi}\simeq {\cal H}_{n}({\cal O}_{X})\simeq {\cal D}^{n}({\cal O}_{X})=
\omega^{n}_{X}$}\rm\smallskip\noindent
 \Cor{1.1}{}Soient $\pi\in{\cal E}(S,n)$ et $k$ un entier naturel. Alors, il existe un unique faisceau ${\cal O}_{X}$-coh\'erent  $\omega^{k}_{\pi}$  v\'erifiant les propri\'et\'es suivantes :\smallskip\noindent
 {\bf(i)} Pour toute installation locale  de $\pi$ (cf Thm 1), on a\par
\centerline{$\omega^{k}_{\pi}:=\sigma^{*}({\cal E}xt^{p}(\sigma_{*}\Omega^{n-k}_{X/S},
  \Omega^{n+p}_{Z/S}))$}
\par\noindent {\bf(ii)} on a des isomorphismes canoniques\par
  \centerline{$\omega^{k}_{\pi}\simeq {\cal
H}om_{{\cal O}_{X}}(\Omega^{n-k}_{X/S},
\omega^{n}_{\pi})$}\par\noindent {\bf(iii)} $\omega^{k}_{\pi}$ est
un ${\cal O}_{X}$- module coh\'erent de profondeur au moins $2$
fibre par fibre sur $X$.\par\noindent
\par\noindent
{\bf(iv)} Si $S$ est un point, ces faisceaux coincident avec les
dualis\'es de Andr\'eotti-Kas-Golovin  des faisceaux
$\Omega^{n-k}_{X}$. \rm\smallskip\smallskip\noindent  \Th{2}{} Soit  $\pi$ un \'el\'ement de
${\cal E}(S,n)$ muni d'une pond\'eration ${\goth X}$.  Alors, on a
les \'equivalences\smallskip\noindent {\bf(i)} $\pi$ est
analytiquement g\'eom\'etriquement plat\smallskip\noindent {\bf(ii)}
il existe un unique morphisme ${\cal O}_{S}$-lin\'eaire continu
 $\ds{{\int}_{\pi,{\goth X}}:{\rm I}\!{\rm
R}^{n}\pi_{!}\Omega^{n}_{X/S}\rightarrow {\cal O}_{S}}$  de
formation compatible \`a l'additivit\'e des pond\'erations sur
$\pi$, stable par changement de base entre espace complexes
r\'eduits, donnant, en particulier,  l'int\'egration usuelle dans le
cas o\`u $S$ est un point.\par\noindent
{\bf(iii)}  il existe  un  unique morphisme canonique  de faisceaux coh\'erents
 ${\cal C}_{\pi,\goth{X}}:\Omega^{n}_{X/S}\rightarrow \omega^{n}_{\pi}$,
prolongeant  le morphisme
naturel $j_{*}j^{*}\Omega^{n}_{X/S}\rightarrow
j_{*}j^{*}\omega^{n}_{\pi}$ et
 v\'erifiant les propri\'et\'es suivantes:\par
{\bf(1)} il est de formation compatible \`a l'additivit\'e des pond\'erations,   aux
 changement de base r\'eduit et satisfait  {\it{la propri\'et\'e de la
 trace relative}},\par
{\bf(2)} il induit  un morphisme de complexes  diff\'erentiels gradu\'es\par
 \centerline{${\cal C}^{\bullet}_{\pi, \goth{X}}:\Omega^{\bullet}_{X/S}\rightarrow \omega^{\bullet}_{\pi}$}
{\bf(3)} $\bullet$ si $\goth{X}$ est la pond\'eration standard,  ${\cal C}_{\pi,\goth{X}}$  donne  la classe fondamentale relative \par de
 [B4],\par
 $\bullet$ si  $\goth{X}$ est la pond\'eration alg\'ebrique (i.e si $\pi$ est plat ou plus g\'en\'eralement de\par {\it Tor-dimension} finie),  ${\cal C}_{\pi,\goth{X}}$ est le morphisme classe fondamentale relative de [A.E] \par \'etendu \`a ce cadre.
\rm\smallskip\medskip\noindent
\Cor{2.1}{}Soit
$\pi_{2}:X_{2}\rightarrow X_{1}$ (resp. $\pi_{1}:X_{1}\rightarrow
S$) un morphisme analytiquement g\'eom\'etriquement plat de
dimension relative $n_{2}$ (resp. $n_{1}$). Alors, le morphisme
compos\'e $\pi:X_{2}\rightarrow S$ est analytiquement
g\'eom\'etriquement plat si et seulement si le morphisme naturel
${\rm I}\!{\rm R}^{n_{1}+n_{2}}{\pi}_{!}(
\Omega^{n_{2}}_{X_{2}/X_{1}} \otimes
\pi_{2}^{*}\Omega^{n_{1}}_{X_{1}/S})\rightarrow{\cal O}_{S}$ se
prolonge en unique morphisme ${\rm I}\!{\rm
R}^{n_{1}+n_{2}}{\pi}_{!}\Omega^{n_{1}+n_{2}}_{X_{2}/S}\rightarrow{\cal
O}_{S}$ ayant les propri\'et\'es du {\it th\'eor\`eme 2} et rendant
commutatif le diagramme
$$\xymatrix{ {\rm I}\!{\rm R}^{n_{1}+n_{2}}{\pi}_{!}(\Omega^{n_{2}}_{X_{2}/X_{1}} \otimes
\pi_{2}^{*}\Omega^{n_{1}}_{X_{1}/S} )\ar[rd]\ar[rr]&&
 {\rm I}\!{\rm R}^{n_{1}+n_{2}}{\pi}_{!}\Omega^{n_{1}+n_{2}}_{X_{2}/S}\ar[ld]^{\Phi}\\
&{\cal O}_{S}&}$$
\Th{3}{} Soit  $\pi$ un \'el\'ement de ${\cal E}(S,n)$ muni d'une
pond\'eration ${\goth X}$. Alors, les assertions suivantes sont
\'equivalentes
\smallskip\noindent
 {\bf(i)} $\pi$ est analytiquement g\'eom\'etriquement plat\smallskip\noindent
{\bf (ii)} le faisceau  $\omega^{n}_{\pi}$ est caract\'eris\'e par
la propri\'et\'e de la trace (relative) et stable par changement de
base entre espaces complexes r\'eduits.\rm\smallskip\noindent On a,
alors, le
\Cor{3.1}{} Soient\par\noindent $\bullet$ $\pi:X\rightarrow
S$ un morphisme $n$-analytiquement g\'eom\'etriquement plat
d'espaces analytiques complexes r\'eduits de dimension localement
fini,\par\noindent $\bullet$ $X_{0}$ (resp. $S_{0}$) l'ouvert dense
de $X$ (resp. $S$) sur lequel $\pi$ est plat et
 $\tilde{\omega}^{n}_{X_{0}/S_{0}}$ le faisceau des
 formes m\'eromorphes r\'eguli\`eres de Kunz-Waldi-Kersken.\smallskip\noindent
 Alors, $\omega^{n}_{\pi}$ est l'unique prolongement coh\'erent du
 faisceau  $\tilde{\omega}^{n}_{X_{0}/S_{0}}$. De plus, la famille de
faisceaux ${\cal O}_{X}$-coh\'erents $\omega^{\bullet}_{\pi}:={\cal
H}om(\Omega^{n-\bullet}_{X/S}, \omega^{n}_{\pi})$ est munie d'une
diff\'erentielle non triviale ${\rm D}$, faisant de
$(\omega^{\bullet}_{X/S}, {\rm D})$  un complexe diff\'erentiel de
$(\Omega^{\bullet}_{X/S}, d_{X/S})$-modules.
\rm\bigskip\noindent
On peut aussi \'enoncer de fa\c con ``presque'ind\'ependante'' du {\it th\'eor\`eme.1}, le
\Th{4}{} Soit $\pi:X\rightarrow S$ un morphisme analytiquement
g\'eom\'etriquement plat. Alors, il existe un unique faisceau ${\cal
O}_{X}$ coh\'erent, $\Lambda^{n}_{X/S}$ v\'erifiant:\par\noindent
{\bf(i)} il est  de profondeur au moins deux fibre par fibre sur
$S$, caract\'eris\'e par la propri\'et\'e de la trace relative et
stable par changement de base entre espaces complexes r\'eduits de
dimension finie,\par\noindent {\bf(ii)} il munit $\pi$ d'un
morphisme d'int\'egration $\int_{\pi}:{\rm I}\!{\rm
R}\pi_{!}\Lambda^{n}_{X/S}\rightarrow {\cal O}_{S}$ compatible avec
l'additivit\'e des cycles, stable par changement de
base entre espaces complexes r\'eduits. De plus, tout  diagramme commutatif d'espaces analytiques complexes
$$\xymatrix{X_{2}\ar[rr]^{\Psi}\ar[rd]_{\pi_{2}}&&X_{1}\ar[ld]^{\pi_{1}}\\
&S&}$$
 avec $\pi_{1}$ (resp. $\pi_{2}$) analytiquement g\'eom\'etriquement plats  de dimension relative $n_{1}$ (resp. $n_{2}$)  et $\Psi$ universellement \'equidimensionnel et propre de dimension relative born\'ee par l'entier   $d:=n_{2}-n_{1}$,  donne un morphisme canonique ${\rm I}\!{\rm
 R}^{d}{\Psi}_{*}\Lambda^{d+n_{1}}_{X_{2}/S}\rightarrow \Lambda^{n_{1}}_{X_{1}/S}$ induisant le diagramme commutatif de faisceaux analytiques
$$\xymatrix{{\rm I}\!{\rm
 R}^{n_{2}}{{\pi_{2}}_{!}}\Lambda^{n_{2}}_{X_{2}/S}\ar[rr]\ar[rd]_{\int_{\pi_{2}}}&&
 {\rm I}\!{\rm
 R}^{n_{1}}{{\pi_{1}}_{!}}\Lambda^{n_{1}}_{X_{1}/S}\ar[ld]^{\int_{\pi_{1}}}\\
&{\cal O}_{S}&}$$
\rm\bigskip\par\noindent
Rappelons que le terme  `` continue'' sera justifi\'e par le fait
que  les groupes  des  sections globales des faisceaux qui
interviennent dans le th\'eor\`eme 2  peuvent \^etre munis de
structures d'espaces vectoriels topologiques {\bf{F.S}} ou
{\bf{Q.F.S}}. Rappelons aussi que le morphisme gradu\'e dont il est
question dans ce m\^eme th\'eor\`eme n'est pas un morphisme
d'alg\`ebres diff\'erentielles car  $\omega^{\bullet}_{\pi}$ n'en
est pas une si $X$ n'est pas normal. On s'en convainc facilement,
dans la situation absolue, avec l'exemple simple donn\'e par
$X=\{(x,y)\in {\Bbb C}^{2}: x^{2}-y^{3}=0\}$. En effet, la fonction
m\'eromorphe ${y\over{x}}$ (resp. la forme ${dx\over{y}}$) d\'efinit
une section de $\omega^{0}_{X}$ (resp. $\omega^{1}_{X}$) mais leur
produit ${dx\over{x}}$ ne d\'efinit pas une section de
$\omega^{1}_{X}$ !
\bigskip\bigskip\noindent
\centerline{{\tite{I. La preuve  du th\'eor\`eme
1.}}}\smallskip\bigskip\noindent {\bf 1.0. Quelques annulations de
cohomologie.}\par\noindent Comme nous le verrons dans un article
ult\'erieur dont le sujet est la notion de paire dualisante en
g\'eom\'etrie analytique complexe, on peut adopter une approche par
la dualit\'e analytique relative pour montrer ce th\'eor\`eme et
constater que le groupe des sections globales de ce faisceaux
apparait comme le ${\cal O}_{S}$- dual topologique du faisceau
analytique ${\rm I}\!{\rm R}^{n}\pi_{!}{\cal O}_{X}$.\par\noindent
 \Lemme{1}{} Soit  $q:Z\rightarrow S$ un morphisme lisse de dimension relative $n+p$  d'espaces analytiques complexes r\'eduits et de dimension
 pure. Soient  $X$ un sous ensemble analytique de $Z$,
 $n$-universellement \'equidimensionnel sur $S$ et ${\cal F}$  un faisceau coh\'erent, localement libre sur  $Z$. Alors
$${\cal H}^{j}_{X}({\cal F})= 0,\,\,\forall\,j<p$$
\dem Si $S$ est r\'eduit \`a un point, c'est un r\'esultat de
Siu-Trautmann [S.T]. Le passage au cas relatif a \'et\'e trait\'e
dans [K1] p.292 o\`u l'on a proc\'ed\'e  par  r\'ecurrence
ascendante sur la dimension de $S$ et descendante sur la codimension
du support $X$. Rappelons bri\`evement la strat\'egie.\par\noindent
Le probl\`eme \'etant de nature locale, on peut supposer donn\'ee
une installation locale
$\xymatrix{X\ar@/_/[rr]_{\pi}\ar[r]^{\sigma}&Z\ar[r]^{q}&S}$ dans
laquelle $\sigma$ est  un plongement de $X$ dans $Z=S\times U\times
B$ avec  $U$ (resp. $B$) un polydisque ouvert relativement compact
de ${\Bbb C}^{n}$ (resp. ${\Bbb C}^{p}$) et $S$ un ouvert de Stein.
On peut, de plus,  supposer ${\cal F}={\cal O}_{Z}$. Il existe un
ouvert dense de $X$ sur lequel le morphisme $\pi$ est une fibration
localement triviale et m\^eme submersif (simple application du
th\'eor\`eme des fonctions implicites).
 Mais si $X$ est de la forme $S\times X_{0}$ (avec $S$ et $X_{0}$ lisses), une formule de type
 K\"unneth donne imm\'ediatement
$${\cal  H}^{j}_{X}({\cal O}_{Z})\simeq {\cal O}_{S}{{\widehat\otimes}_{\atop{{\Bbb C}}}}{\cal H}^{j}_{X_{0}}({\cal O}_{U\times B})$$
qui est donc nul pour $j<p$.\par\noindent Soit $Y$, le
compl\'ementaire de cet ouvert, qui est un sous ensemble analytique
ferm\'e d'int\'erieur vide de $X$. Soit $S'$ le ferm\'e
d'int\'erieur vide de $S$ constitu\'e des points $s$ de $S$ tel que
$Y\cap X_{s}$ ne soit pas d'int\'erieur vide dans $X_{s}$. Comme au
dessus de $S-S'$, $Y$ est de codimension au moins $p+1$ et que la
propri\'et\'e de {\it{locale libert\'e}} est stable par image
r\'eciproque et changement de base, l'hypoth\`ese de r\'ecurrence
sur la codimension du support et la suite classique de
Mayer-Vietoris, nous ram\`ene \`a l'\'etude de ${\cal
H}^{j}_{X'}({\cal F}')$  o\`u $X':=\ds{X{\times_{\atop S}} S'}$ et
$Z':= \ds{Z{\times_{\atop S}} S'\simeq S'\times U\times B}$ et
${\cal F}'$ est l'image r\'eciproque de ${\cal F}$ dans ce
changement de base.\par\noindent Alors, $S'$ \'etant de dimension
strictement plus petite que celle de $S$,  l'hypoth\`ese de
r\'ecurrence faite sur $S$ permet de conclure
$\blacksquare$\rm\par\noindent
 \Lemme{2}{}
 Soit   $q:Z\rightarrow S$ un morphisme lisse de dimension relative $n+p$  d'espaces analytiques complexes r\'eduits et de dimension pure.
 Soit ${\cal F}$ un faisceau coh\'erent sur $Z$ dont le support $X$ est $n$-universellement \'equidimensionnel sur $S$. Alors
$${\cal E}xt^{k}_{{\cal O}_{Z}}({\cal F}, \Omega^{n+p}_{Z/S}) = 0,\,\,\,\forall\,k<p $$
\dem  Dans le cas absolu, l'annulation des faisceaux coh\'erents
${\cal E}xt^{k}_{{\cal O}_{Z}}({\cal F}, \Omega^{n+p}_{Z})$,  pour
tout entier $k<p $, d\'ecoule du r\'esultat classique ( [Gr1]
p.30-33,th 2.2 et prop.2.9) disant que si {\it{$A$ est un anneau
local noeth\'erien, $M$ et $N$ deux $A$-modules, $k\geq 0$ un
entier, alors $Ext^{j}_{A}(M, N) = 0$ pour tout $j<k$ si et
seulement si ${\rm Prof}_{A_{\goth p}}N_{\goth p}\geq k$ pour tout
id\'eal premier ${\goth p}\in {\rm Supp}(M)$}}.\par\noindent Dans la
situation relative, on adopte la m\^eme strat\'egie que dans le
lemme pr\'ec\'edent en supposant, sans enfreindre la
g\'en\'eralit\'e,  ${\cal F}={\cal O}_{X}$. Comme les propri\'et\'es
d'\^etre \`a fibres de dimension pure constante et ouvert sont
stables par changement de base, le morphisme induit $X'\rightarrow
S'$ est encore universellement \'equidimensionnel; de plus, $X'$
s'identifie \`a un sous espace d'int\'erieur vide dans $X$.  Comme
dans le {\it lemme 1} et  gr\^ace \`a une suite de Mayer-Vietoris,
on se ram\`ene facilement  \`a l'\'etude du faisceau coh\'erent
$${\cal E}xt^{k}({\cal O}_{X'}, \Omega^{n+p}_{Z'/S'})$$ avec
$X':=\ds{X{\times_{\atop S}} S'}$ et $Z':= \ds{Z{\times_{\atop S}}
S'\simeq S'\times U\times B}$, dont on d\'eduit le r\'esultat
gr\^ace aux hypoth\`eses de r\'ecurrences.
$\,\,\blacksquare$\smallskip\noindent En conservant  les notations
et hypoth\`eses du {\it lemme 1}  et posant $\pi:X\rightarrow S$ le
morphisme induit sur $X$ par $q$, on a le
 \Lemme{3}{} Soit $F$  un
ferm\'e analytique de $Z$ de codimension au moins deux fibre par
fibre (i.e ${\rm Codim}_{\pi^{-1}(s)}(F\cap \pi^{-1}(s))\geq
2,\,\,\,\forall\,s\in S$). Alors
$${\cal H}^{0}_{F}({\cal E}xt^{p}(\sigma_{*}{\cal O}_{X},
\Omega^{n+p}_{Z/S}))= {\cal H}^{1}_{F}({\cal
E}xt^{p}(\sigma_{*}{\cal O}_{X}, \Omega^{n+p}_{Z/S})) = 0$$ \dem
 Le cas absolu d\'ecoule du r\'esultat
 classique d'alg\`ebre ([Gr1])  disant que si
{\it{$A$ est  un anneau local noeth\'erien, $M$ et $N$ deux
$A$-modules de type fini avec ${\rm Prof}(M)\geq 2$ alors ${\rm
Prof}{\cal H}om(N, M)\geq 2$.}}\smallskip\noindent Comme le
probl\`eme est de nature locale, on peut supposer $Z:=S\times
U\times B$ avec $U$ (resp. $B$) un polydisque ouvert de ${\Bbb
C}^{n}$ (resp. ${\Bbb C}^{p}$). Par ailleurs, $\pi$ \'etant  \`a
fibres de dimension pure $n$, on peut se donner une installation
locale correspondant au  diagramme commutatif
$$\xymatrix{X\ar[d]^{\pi} \ar [dr]_{f}\ar [r]^{\sigma}
 & S\times U\times B\ar[d]^{p_{1}}\\
 S&\ar[l]^{p_{2}}S\times U}$$
dans lequel  $f$ est un morphisme fini et surjectif, $\sigma$ est un
plongement local, $p_{1}$ et $p_{2}$ \'etant les projections
canoniques.\par\noindent On se ram\`ene ainsi au cas fini et
utiliser le fait bien connu (cf [B.S] par exemple) que pour tout
morphisme fini et surjectif d'espaces complexes, $f:X\rightarrow Y$,
on a, ${\rm Prof}({\cal F})= {\rm Prof}(f_{*}{\cal F})$ pour tout
faisceau coh\'erent ${\cal F}$ sur $X$. Posons $m=n+p$ et
$Y:=S\times U$.  Comme $f_{*}{\cal E}xt^{m-n}({\cal O}_{X},
\Omega^{m}_{Z/S})\simeq {\cal H}om(f_{*}{\cal O}_{X},
\Omega^{n}_{Y/S})$\footnote{$^{(1)}$}{Dans cette situation
particuli\`ere, la dualit\'e analytique relative de [RRV], donne
$$f_{*}{\rm I}\!{\rm R}{\cal H}om({\cal O}_{X},
\Omega^{n+p}_{Z/S}[p])\simeq {\rm I}\!{\rm R}{\cal H}om(f_{*}{\cal
O}_{X}, \Omega^{n}_{Y/S})$$ dont on prend la cohomologie de degr\'e
0.; ce qui traduit une d\'eg\'en\'erescence de suite spectrale
classique.},  il suffit simplement de constater que
$\Omega^{n}_{Y/S}$ \'etant localement libre (donc de profondeur au
moins deux fibre par fibre), le faisceau  ${\cal H}om(f_{*}{\cal
O}_{X}, \Omega^{n}_{Y/S})$ est aussi de profondeur au moins
deux.\par\noindent Pour s'en convaincre, on consid\`ere un sous
ensemble analytique $F$ de $X$ tel que $F\cap \pi^{-1}(s)$ soit
d'int\'erieur vide dans $X_{s}$ et l'on regarde les suites
spectrales
$$E^{i,j}_{2}= {\cal E}xt^{i}({\cal
O}_{T}, {\cal H}^{j}_{F}(\Omega^{n}_{Y/S}))\,\,\, {\rm et}\,\,
'E^{i,j}_{2}= {\cal H}^{i}_{F}({\cal E}xt^{i}({\cal O}_{T},
\Omega^{n}_{Y/S}))$$ de m\^eme aboutissement $\ds{{\cal
E}xt^{i+j}_{F}({\cal O}_{T}, \Omega^{n}_{Y/S})}$ dans laquelle on a
pos\'e $ {\cal O}_{T}:=f_{*}{\cal O}_{X}$. Alors, les termes de bas
degr\'e donnent l'isomorphisme
$${\cal H}^{i}_{F}({\cal H}om({\cal O}_{T},  \Omega^{n}_{Y/S})\simeq {\cal H}om({\cal O}_{T},{\cal H}^{i}_{F}(\Omega^{n}_{Y/S}))$$
qui prouve notre assertion puisque le {\it lemme 2}  donne  ${\cal
H}^{i}_{F}(\Omega^{n}_{Y/S}))=0\,\,{\rm
pour}\,\,i=0,1$.$\,\,\blacksquare$
\smallskip\noindent
On en d\'eduit aussit\^ot le \Cor{1}{} Pour tout ouvert $U$ de $X$
dense fibre par fibre sur $S$, le morphisme naturel $ \ds{\Gamma(X,
\omega^{n}_{\pi})\rightarrow \Gamma(U, \omega^{n}_{\pi})}$ est
injectif.\rm\smallskip\noindent
 dont la cons\'equence est  l'unicit\'e de $\omega^{n}_{\pi}$ \`a
 isomorphisme canonique pr\`es.\smallskip\noindent
{\bf 1.0.1.  Remarques.}\par\noindent {\bf(i)} Le {\it lemme 2} se
d\'eduit ais\'ement du lemme 5.2.2 p.46 de [A.LJ] dont l'\'enonc\'e
se pr\'esente sous la forme :\par\noindent {\it{Soit $k$ un anneau
local noeth\'erien, $N$ une $k$-alg\`ebre locale
diff\'erentiellement lisse sur $k$, $M$ un $N$-module de type fini.
Soit ${\goth K}$ le corps r\'esiduel de $k$ et $m$ la dimension de
Krull de l'anneau r\'egulier $N\otimes_{k} {\goth K}$. Posons $n:=
{\rm dim}_{N\otimes_{k} {\goth K}}M\otimes_{k} {\goth K}$ et
$p=m-n$. Alors $Ext^{j}_{N}(M, N)=0$ pour tout $j<p$}} et qui peut
\^etre vu comme une version relative du r\'esultat d'alg\`ebre
cit\'e ci-dessus.\par\noindent On peut remarquer que si ${\cal F}$
est plat relativement \`a $S$, le th\'eor\`eme du changement de base
pour les ``${\cal E}xt$'' donne imm\'ediatement le
r\'esultat.\par\noindent {\bf(ii)} Pour le {\it lemme 3}, on peut se
r\'ef\'erer \`a [A.LJ], en adaptant sans difficult\'e le  lemme
5.2.3 \`a notre situation en utilisant, cette fois, les suites
spectrales $\ds{E^{i,j}_{2}= {\cal E}xt^{i}(\sigma_{*}({\cal
O}_{X}), {\cal H}^{j}_{Y}(\Omega^{n+p}_{Z/S}))}$ et
$\ds{'E^{i,j}_{2}= {\cal H}^{i}_{Y}({\cal E}xt^{i}(\sigma_{*}({\cal
O}_{X}), \Omega^{n+p}_{Z/S}))}$ de m\^eme aboutissement $\ds{{\cal
E}xt^{i+j}_{Y}(\sigma_{*}({\cal O}_{X}), \Omega^{n+p}_{Z/S})}$, on
voit que
$$'E^{i,j}_{2}= 0,\,\forall\,(i,j):\,i+j<p;\,\,\,\,'E^{i,j}_{2}= 0,\,\forall\,j:\,j<p$$
$$E^{i,j}_{2}= 0,\,\forall\,(i,j):\,i+j<p+2;\,\,\,\, E^{i,j}_{2}=
0,\,\forall\,j:\,j<p+2$$ et donc $\ds{{\cal H}^{0}_{Y}({\cal
E}xt^{p}(\sigma_{*}({\cal O}_{X}), \Omega^{n+p}_{Z/S}))=0}$.
\par\noindent L'annulation du terme $'E^{1,p}_{2}$ est moins
automatique, certes, mais facile \`a \'etablir puisqu'il suffit
simplement de constater que
$$'E^{1,p}_{2}= 'E^{1,p}_{3}=\cdots= 'E^{1,p}_{\infty}=0$$
sachant que $'E^{i,j}_{r+1}$ est constitu\'e des classes de cycles
de $'E^{i,j}_{r}$ muni de la diff\'erentielle $\ds{d^{i,j}_{r}:
{'E}^{i,j}_{r}\rightarrow {'E}^{i+r,j-r+1}_{r}}$ (cf [C.E]
p.326).\smallskip\noindent { \bf 1.1. Preuve du th\'eor\`eme
1.}\smallskip\noindent Dans un premier temps, on  va montrer que le
faisceau coh\'erent $\sigma^{*}{\cal E}xt^{p}(\sigma_{*}{\cal
O}_{X}, \Omega^{n+p}_{Z/S})$, qui est de profondeur au moins deux
fibre par fibre sur $S$ en vertu du {\it lemme 3}, est en fait
ind\'ependant de la factorisation locale choisie et d\'efini par
suite un faisceau intrins\`eque sur $X$. On montrera, ensuite, les
propri\'et\'es \'enonc\'ees.\smallskip\noindent {\bf{(a) Recollement
des faisceaux}} ${\cal E}xt^{p}(\sigma_{*}{\cal O}_{X},
\Omega^{n+p}_{Z/S})$\smallskip\noindent
 Supposons donn\'e un diagramme commutatif
$$\xymatrix{Z_{1}\ar[rd]_{p_{1}}&X\ar[d]^{\pi}\ar@{_{(}->}[l]_{\sigma_{1}}\ar@{^{(}->}[r]^{\sigma_{2}}&Z_{2}
\ar[ld]^{p_{2}}\\
&S&}$$ dans lequel $\sigma_{1}$, ${\sigma_{2}}$ sont des plongements
locaux dans des espaces lisses sur $S$ munis de projections
canoniques $p_{1}$ et $p_{2}$. \par\noindent Il est, alors, possible
de le compl\'eter naturellement en le diagramme
$$\xymatrix{&Z_{1}{\times_{S}}Z_{2}\ar[ld]_{q_{1}}\ar[rd]^{q_{2}}&\\
Z_{1}\ar[rd]_{p_{1}}&X\ar[d]^{\pi}
\ar@{^{(}->}[u]^{\sigma_{1,2}}\ar@{_{(}->}[l]^{\sigma_{1}}\ar@{^{(}->}[r]_{\sigma_{2}}
&Z_{2}\ar[ld]^{p_{2}}\\
&S&}$$ o\`u $\sigma_{1,2}$ est encore un plongement dans un $S$-
espace complexe lisse sur $Z_{1}$ et $Z_{2}$.\par\noindent Ces
consid\'erations nous ram\`enent \`a la situation donn\'ee par le
diagramme commutatif de $S$-espaces complexes
$$\xymatrix{&Z\times V\ar@<2pt>[d]^{q}\\
X\ar[ru]^{\rho}\ar[r]_{\sigma}&Z\ar@<2pt>[u]^{h}}$$ dans lequel  $V$
 est un ouvert d'un certain espace num\'erique ${\Bbb C}^{N}$ et
$h:Z\rightarrow Z\times V$ une section de $q$ (i.e $ho\sigma=\rho$,
$qo\rho=\sigma$) faisant de $Z$ un r\'etract de $Z_{1}:= Z\times V$.
Notons $m_{1}$ (resp. $m$) la dimension $S$-relative de $Z_{1}$
(resp. $Z$) et $r:=m_{1}-m$.\par\noindent Comme $Z$ et $Z_{1}$ sont des espaces $S$- plats
et  \`a fibres Gorenstein (donc Cohen-Macaulay),  le  {\it lemme 2} (qui est d'une application classique dans cette situation)  assure la d\'eg\'en\'erescence de  la suite spectrale
$${\cal E}xt^{i}({\cal O}_{X}, {\cal E}xt^{j}({\cal O}_{Z}, \Omega^{m_{1}}_{Z_{1}/S}))\Longrightarrow {\cal E}xt^{i+j}({\cal O}_{X}, \Omega^{m_{1}}_{Z_{1}/S})$$
donnant, en particulier,  l'isomorphisme
$${\cal E}xt^{p}({\cal O}_{X}, {\cal E}xt^{m_{1}-m}({\cal O}_{Z}, \Omega^{m_{1}}_{Z_{1}/S}))\simeq {\cal E}xt^{p+r}({\cal O}_{X}, \Omega^{m_{1}}_{Z_{1}/S})$$
Mais, gr\^ace \`a la classe fondamentale relative de $Z/S$ dans
$Z_{1}$, on a un isomorphisme (donn\'e par cup produit par la classe
fondamentale relative dont l'inverse est le morphisme r\'esidu dans
ce cas!)
$$\Omega^{m}_{Z/S}\simeq {\cal E}xt^{m_{1}-m}({\cal O}_{Z}, \Omega^{m_{1}}_{Z_{1}/S})$$
On peut aussi invoquer le fait que, dans cette situation
particuli\`ere, le complexe dualisant
 relatif ${\cal D}^{\bullet}_{Z/S}$ est quasi-isomorphe au complexe $\Omega^{m}_{Z/S}[m]$ puisqu'il
  ne poss\`ede qu'un faisceau d'homologie en degr\'e $-m$ qui est par d\'efinition $\ds{{\cal E}xt^{m_{1}-m}({\cal O}_{Z},
  \Omega^{m_{1}}_{Z_{1}/S})}$.  On en d\'eduit l'ind\'ependance vis-\`a-vis du $S$-plongement
choisi et l'on pose, alors,
 $$\omega^{n}_{\pi}:= \sigma^{*}{\cal
E}xt^{p}(\sigma_{*}{\cal O}_{X}, \Omega^{n+p}_{Z/S})$$ Ce faisceau
coh\'erent v\'erifie les propri\'et\'es  annonc\'ees dans le
th\'eor\`eme. En effet, il est clairement  coh\'erent et de
profondeur au moins deux fibre par fibre en vertu de la construction
et du {\it lemme 3}.  Aux points r\'eguliers de $\pi$,il coincide
naturellement, avec le faisceau usuel des formes holomorphes
$S$-relatives. L'isomorphisme \'etant donn\'e par la classe
fondamentale relative en ces points. Par construction, il  est
compatible aux inclusions ouvertes dans le sens de Verdier puisque,
si   $\pi_{i}:X_{i}\rightarrow S,\,\,i=1,2$ sont deux morphismes
analytiques  universellement $n$-\'equidimensionnels  et $U$ un
ouvert de $X_{1}$ muni
  de deux inclusions ouvertes $j_{i}:U\rightarrow X_{i}$
  tels  que le diagramme
  $$\xymatrix{&U\ar[ld]_{j_{1}}\ar[rd]^{j_{2}}&\\
  X_{1}\ar[rd]_{\pi_{1}}&&X_{2}\ar[ld]^{\pi_{2}}\\
  &S&}$$
  soit commutatif, il est facile de voir que
  $$j_{1}^{*}(\omega^{n}_{\pi_{1}}) =
  j_{2}^{*}(\omega^{n}_{\pi_{2}})$$
Pour s'en convaincre, on se ram\`ene au diagramme commutatif
 $$\xymatrix{U\times_{S} U\ar[rd]_{\pi_{U}}\ar[r]&X_{1}\times_{S} X_{2}\ar[d]^{\pi}\\
&S}$$
 en remarquant que  $\pi$ est encore universellement $n$-\'equidimensionnel
 puisque les  propri\'et\'es  d'\^etre ouvert et d'avoir des fibres de
 dimension constante sont stables par changement de base. De plus,
 comme $X_{1}$ et $X_{2}$ sont de dimension pure, il en va de m\^eme
  de leur produit fibr\'e qui est universellement $n$-\'equidimensionnel  sur $S$.\par\noindent
On doit principalement v\'erifier que   pour tout ouvert  $U$ de $X$ muni de
l'injection naturelle $j:U\rightarrow X$ et de la restriction de
$\pi$ \`a $U$ que l'on note $\pi_{U}$, on a
$$\omega^{n}_{\pi}|_{U}=\omega^{n}_{\pi_{U}}$$
Mais cela d\'ecoule de la construction m\^eme de ce
faisceau puiqu'il est bien connu que les faisceaux {\it Ext} poss\`ede de
 bonnes propri\'et\'es de variances relativement aux inclusions ouvertes;
  ce qui n'est qu'un cas particulier de {\it la stabilit\'e par changement
   de base plat}. Pour s'en assurer, il suffit, gr\^ace \`a la nature locale du probl\`eme, de  factoriser  localement  $\pi$ en un
plongement dans un espace complexe $Z$ muni d'une projection lisse  sur $S$. On consid\`ere, alors, un morphisme d'espaces complexes r\'eduits (de dimension pure) $\eta:S_{1}\rightarrow S$ et le diagramme commutatif (cart\'esien) de changement de base \par
\centerline{$\xymatrix{&Z_{1}\ar[ldd]_{q_{1}}\ar[r]^{\phi}&Z\ar[ldd]^{q}&\\
X_{1}\ar[d]_{\pi_{1}}\ar[ru]^{\sigma_{1}}\ar[r]_{\Theta}&X\ar[ru]^{\sigma_{1}}\ar[d]_{\pi}&\\
S_{1}\ar[r]_{\eta}&S&}$} \noindent
 Les morphismes lisses, plats  et plus g\'en\'eralement
universellement ouverts \'etant stables par changement de base, on
se ram\`ene \`a \'etablir l'existence d'une fl\`eche
$\Theta^{*}\omega^{n}_{\pi}\rightarrow
\omega^{n}_{\pi_{1}}$.\par\noindent Mais les \'egalit\'es
fonctorielles ${\rm I}\!{\rm L}\Theta^{*}{\rm I}\!{\rm R}{\cal
H}om(A,B)= {\rm I}\!
 {\rm R}{\cal H}om({\rm I}\!{\rm L}\Theta^{*}A,{\rm I}\!{\rm
 L}\Theta^{*}B)$ et la platitude de $\Theta$,
  traduite par l'\'egalit\'e ${\rm I}\!{\rm L}\Theta^{*}=\Theta^{*}$,
  donnent, en vertu des annulations du {\it
lemme 2}, un morphisme canonique de changement de base
$\Theta^{*}{\cal E}xt^{p}({\cal O}_{X},
\Omega^{n+p}_{Z/S})\rightarrow {\cal E}xt^{p}({\cal O}_{X_{1}},
\Omega^{n+p}_{Z_{1}/S_{1}}) $.\smallskip\noindent Il nous reste \`a
montrer le point  {\bf(v)} disant que ce faisceau munit $\pi$ d'un
morphisme canonique
 $${\rm
I}\!{\rm R}^{n}\pi_{!}\omega^{n}_{\pi}\rightarrow {\cal O}_{S}$$
compatible aux localisations sur $X$ (resp. $S$) et aux changements
de bases plats. Pour cela, nous allons d'abord exhiber une telle
fl\`eche relativement \`a une installation locale de $\pi$ puis
proc\'eder \`a la globalisation gr\^ace au lemme de
Reiffen.\smallskip $\bullet$ {\bf Construction locale.
}\par\noindent Soient  $s_{0}$ un point de $S$ et  $x$ un point de
la fibre $\pi^{-1}(s_{0})$.  Consid\'erons une factorisation locale
de $\pi$ en $x$  donn\'ee (abusivement) par
$\xymatrix{X\ar@/_/[rr]_{\pi}\ar[r]^{f}&Y\ar[r]^{q}&S}$ dans
laquelle $f$ est fini, ouvert  et surjectif sur $Y=S\times U$ avec
$U$ polydisque ouvert relativement compact de ${\Bbb C}^{n}$, $S$ un
ouvert de Stein et $q$ la projection canonique.\par\noindent Comme
on l'a d\'ej\`a vu dans la preuve du {\it lemme 3}, on a
l'isomorphisme
$$f_{*}{\cal E}xt^{p}({\cal O}_{X}, \Omega^{n+p}_{Z/S})\simeq {\cal H}om(f_{*}{\cal O}_{X}, \Omega^{n}_{Y/S})$$
qui, compos\'e  avec la fl\`eche canonique ${\cal
O}_{Y}\rightarrow f_{*}{\cal O}_{X}$, donne  un morphisme ``naturel''
$$f_{*}{\cal E}xt^{p}({\cal O}_{X}, \Omega^{n+p}_{Z/S})\rightarrow \Omega^{n}_{Y/S}$$
duquel r\'esulte
$${\rm I}\!{\rm R}^{n}q_{!}f_{*}{\cal E}xt^{p}({\cal O}_{X}, \Omega^{n+p}_{Z/S})\rightarrow
{\rm I}\!{\rm R}^{n}q_{!}\Omega^{n}_{Y/S}.$$ Mais, pour tout ouvert
$S'$ de Stein dans $S$, une formule de type K\"unneth donne
$${\rm I}\!{\rm R}^{n}q_{!}\Omega^{n}_{Y/S}(S')\simeq {\cal
O}_{S}(S') {{\widehat{\otimes}}{\atop_{{\Bbb C}}}}{\rm H}^{n}_{c}(U,
\Omega^{n}_{U})$$ qui, eu \'egard \`a l'isomorphisme ${\rm
H}^{n}_{c}(U, \Omega^{n}_{U})\simeq {\Bbb C}$ donn\'e par la
dualit\'e de Serre, assure l'isomorphisme ${\rm I}\!{\rm
R}^{n}q_{!}\Omega^{n}_{Y/S}\simeq {\cal O}_{S}$. Alors $f$ \'etant
fini, la d\'eg\'en\'erescence de la suite spectrale de Leray fournit
l'isomorphisme de foncteurs ${\rm I}\!{\rm R}^{n}\pi_{!} \simeq {\rm
I}\!{\rm R}^{n}q_{!}f_{*}$ duquel r\'esulte la fl\`eche
$${\rm I}\!{\rm R}^{n}\pi_{!}{\cal E}xt^{p}({\cal O}_{X}, \Omega^{n+p}_{Z/S})\rightarrow
{\cal O}_{S}$$
\indent $\bullet$ {\bf  Construction globale.}\par\noindent
Commen\c cons par remarquer que la construction d'une telle
fl\`eche, qui est toujours de nature locale sur $S$, est aussi  de
nature locale sur $X$ en vertu du  lemme de Reiffen donnant
l'exactitude \`a droite du foncteur ${\rm I}\!{\rm R}^{n}\pi_{!}$. Par ailleurs, $\omega^{n}_{\pi}$ \'etant de nature locale d'apr\`es le point {\bf(iv)}, on a,  pour tout ouvert  $U$ de $X$ muni de l'injection naturelle $j:U\rightarrow X$ et de
la restriction de $\pi_{U}$ (restriction de $\pi$ \`a $U$), $$\omega^{n}_{\pi}|_{U}=\omega^{n}_{\pi_{U}}$$
Cela \'etant dit, supposons $X$ paracompact et compl\`etement
paracompact\footnote{$^{(2)}$}{Tout ouvert de $X$ est un espace
paracompact; en particulier, pour toute famille paracompactifiante
de supports, l'\'etendue $\bigcup_{F\in \Phi}F$ est un ouvert
paracompact.}  et consid\'erons un ouvert de Stein (que l'on notera
encore $S$) de $S$ et un recouvrement ouvert localement fini
$(X_{\alpha})_{\alpha\in A}$ de $X$ qui soit $S$-adapt\'e et muni
d'installations locales
 $$\xymatrix{X_{\alpha}\ar[rdd]_{\pi_{\alpha}}\ar[rr]^{\sigma_{\alpha}}\ar[rd]^{f_{\alpha}}&&Z_{\alpha}\ar[ld]_{p_{\alpha}}\ar[ldd]^{r_{\alpha}}\\
&Y_{\alpha}\ar[d]^{q_{\alpha}}&\\&S&}$$ avec $\sigma_{\alpha}$  un
plongement local, $\pi_{\alpha}$ la restriction de $\pi$ \`a
$X_{\alpha}$, $Y_{\alpha}$ et $Z_{\alpha}$ sont lisses sur $S$,
$p_{\alpha}$,$q_{\alpha}$ et $r_{\alpha}$ lisses.\par\noindent Comme
$$\omega^{n}_{\pi}|_{X_{\alpha}} = \omega^{n}_{\pi_{\alpha}}\simeq
{\sigma_{\alpha}}^{*}{\cal E}xt^{p}({\cal O}_{X_{\alpha}},
\Omega^{n+p}_{{Z_{\alpha}}/S}),$$ on d\'eduit,  du cas local pr\'ec\'edent, une collection de morphismes
$${\rm I}\!{\rm R}^{n}{\pi_{\alpha}}_{!}\omega^{n}_{{\pi}_{\alpha}}\rightarrow {\cal O}_{S}.$$
Mais l'exactitude \`a droite du foncteur ${\rm I}\!{\rm
R}^{n}\pi_{!}$ (d\^ue au lemme de Reiffen) qui assure la
surjectivit\'e de la fl\`eche naturelle
$$ \bigoplus_{\atop\alpha\in A}{\rm I}\!{\rm
R}^{n}{\pi_{\alpha}}_{!} \omega^{n}_{{\pi_{\alpha}}} \rightarrow
{\rm I}\!{\rm R}^{n}{\pi_{!}}{\omega^{n}_{\pi}}$$ permet de
produire, par recollement sur $X$, le morphisme d\'esir\'e. La
compatibilit\'e aux localisations sur $X$ et $S$ est ais\'ee \`a
v\'erifi\'ee.\smallskip\noindent {\bf(vi)} Si $\pi$ est propre, le
{\it corollaire 3} de [KI], {\bf (3.2.1.2)} p.72, montre bien que le
faisceau coh\'erent $\omega^{n}_{X/S}:= {\cal H}^{-n}(\pi^{!}({\cal
O}_{S}))$ est localement isomorphe au faisceau $\omega^{n}_{\pi}$.
Par ailleurs, au vu des propri\'et\'es de $\omega^{n}_{X/S}$ , il
nous suffit simplement de montrer qu'il existe une fl\`eche non
triviale $\omega^{n}_{\pi}\rightarrow\omega^{n}_{X/S}$ pour en
d\'eduire qu'ils sont isomorphes globalement. Or ceci d\'ecoule de
la propri\'et\'e ``dualisante''  universelle de $\omega^{n}_{X/S}$
et du  point  {\bf(v)} pr\'ec\'edent; puisqu'alors la donn\'ee du
morphisme $ {\rm I}\!{\rm
 R}^{n}{\pi}_{*}\omega^{n}_{\pi}\rightarrow {\cal O}_{S}$
est strictement \'equivalente, par dualit\'e,  \`a la donn\'ee d'un
morphisme (d'ailleurs injectif) $\omega^{n}_{\pi}\rightarrow
\omega^{n}_{X/S}$. Mais ces deux faisceaux (sans torsion), \'etant
localement isomorphes, sont automatiquement globalement
isomorphes.\smallskip\noindent {\bf(vii)} Si dans le  diagramme
commutatif d'espaces analytiques complexes
$$\xymatrix{X_{2}\ar[rr]^{\Psi}\ar[rd]_{\pi_{2}}&&X_{1}\ar[ld]^{\pi_{1}}\\
&S&}$$ tous les morphismes sont propres, le r\'esultat est une
cons\'equence imm\'ediate de  {\bf(vi)} et du {\it
 corollaire 3} de [KI], {\bf (3.2.1.2)} p.73.\par\noindent
 Dans la situation de {\bf(vii)}, on remarque que la construction d'une fl\`eche ${\rm I}\!{\rm
 R}^{d}{\Psi}_{*}\omega^{n_{2}}_{\pi_{2}}\rightarrow
 \omega^{n_{1}}_{\pi_{1}}$ est de nature locale sur $X_{2}$ et $X_{1}$ en
 vertu du lemme de Reiffen. Il nous suffit, alors, de v\'erifier que
 pour tout germe de param\'etrisation locale $f:X_{1}\rightarrow S\times
 U$, on a un morphisme (bien d\'efini) $f_{*}{\rm I}\!{\rm
 R}^{d}{\Psi}_{*}\omega^{n_{1}+d}_{\pi_{2}}\rightarrow
 \Omega^{n_{1}}_{{S\times U}/S}$ et donc se ramener au cas o\`u
 $X_{1}:=S\times U$ et $\pi_{1}$ la projection canonique $q:S\times
 U\rightarrow S$. Alors, en vertu de l'hypoth\`ese sur la dimension
 relative des morphismes, les annulations d\^ues au lemme de
 Reiffen donnent l'isomorphisme $${\rm I}\!{\rm
 R}^{n_{1}}{q_{!}}{\rm I}\!{\rm
 R}^{d}{\Psi}_{*}\omega^{n_{1}+d}_{\pi_{2}}\simeq {\rm I}\!{\rm
 R}^{n_{1}+d}{\pi_{2}}_{*}\omega^{n_{1}+d}_{\pi_{2}}$$
qui permet d'obtenir, gr\^ace \`a {\bf (v)}, le morphisme ${\rm
I}\!{\rm
 R}^{d}{\Psi}_{*}\omega^{n_{1}+d}_{\pi_{2}}\rightarrow
 \Omega^{n_{1}}_{{S\times U}/S}$. Comme les images directes sup\'erieures
  ${\rm I}\!{\rm
 R}^{d}{\Psi}_{*}\omega^{n_{1}+d}_{\pi_{2}}$ sont coh\'erentes (puisque
  $\Psi$ est propre), la caract\'erisation du faisceau
  $\omega^{n_{1}}_{\pi_{1}} $ suffit pour conclure. En effet, comme
  $f_{*}\omega^{n_{1}}_{\pi_{1}}\simeq {\cal H}om(f_{*}{\cal O}_{X},
  \Omega^{n_{1}}_{{S\times U}/S})$, tout morphisme
  $f_{*}{\cal F}\rightarrow \Omega^{n_{1}}_{{S\times U}/S}$,
   avec ${\cal F}$ coh\'erent sur $X$, induira naturellement un
   morphisme $ {\cal H}om(f_{*}{\cal O}_{X}, f_{*}{\cal
   F})\rightarrow {\cal H}om(f_{*}{\cal O}_{X},
  \Omega^{n_{1}}_{{S\times U}/S})$. Mais $f$ \'etant fini, cela
  correspond exactement \`a un morphisme $f_{*}{\cal F}\rightarrow
  f_{*}\omega^{n_{1}}_{\pi_{1}}$ qui, induira, \`a son tour, puisque
  $f$ est arbitraire (et fini) un morphisme global
  ${\cal F}\rightarrow\omega^{n_{1}}_{\pi_{1}}$.
\smallskip\noindent
{\bf(viii)}  Dans le cas absolu, il coincide naturellement avec les
faisceaux dualisants de Andr\'eotti-Kas- Golovin- Grothendieck (cf
[KI], {\bf (3.0.1)}, p.61):
$${\cal D}^{n}({\cal O}_{X}) = {\cal H}_{-n}({\cal
O}_{X})=\omega^{n}_{X} \,\,\,\,\,\,\blacksquare$$
\bigskip\noindent {\bf 1.32 Preuve du
corollaire 1.1}\smallskip\noindent Les m\^emes arguments que ceux
avanc\'es pour montrer le caract\`ere intrins\`eque de la
caract\'erisation locale du faisceau $\omega^{n}_{\pi}$, permettent
de voir  que, pour chaque entier $k$, le faisceau coh\'erent
$\sigma^{*}({\cal E}xt^{p}(\sigma_{*}\Omega^{n-k}_{X/S},
  \Omega^{n+p}_{Z/S}))$ est ind\'ependant du plongement choisi. On notera $\omega^{k}_{\pi}$ le faisceau coh\'erent sur
  $X$ d\'efini par recollement de ces donn\'ees locales.\par\noindent
Comme le faisceau coh\'erent $\sigma_{*}\Omega^{n-k}_{X/S}$ est \`a support dans $X$ qui est fibre par fibre de codimension $p$ dans $Z$ (avec les notations pr\'ec\'edentes), le {\it lemme 2} assure les annulations
$${\cal E}xt^{j}(\sigma_{*}\Omega^{n-k}_{X/S},
  \Omega^{n+p}_{Z/S})=0,\,\,\,\forall\,j< p,\,\,\forall\, k\in\{0,\cdots,n\}$$
Alors,  utilisant la d\'eg\'en\'erescence de suites spectrales
classiques (cf les lemmes {\bf 2} et {\bf 3}, p.7-8), il est facile de
 produire un isomorphisme canonique
$${\cal H}om_{{\cal O}_{Z}}(\sigma_{*}\Omega^{n-k}_{X/S},{\cal E}xt^{p}(\sigma_{*}{\cal O}_{X},
\Omega^{n+p}_{Z/S}))\simeq {\cal
E}xt^{p}(\sigma_{*}\Omega^{n-k}_{X/S},
  \Omega^{n+p}_{Z/S})$$
On constitue, ainsi, en rempla\c cant $X$, $Z$, $\sigma$
respectivement par  $X_{\alpha}$, $Z_{\alpha}$, $\sigma_{\alpha}$,
une donn\'ee simplicielle (param\'etr\'ee par les installations
locales de $\pi$) qui, gr\^ace \`a l'existence de $\omega^{n}_{\pi}$
et ses propri\'et\'es fonctorielles (cf le {\it th\'eor\`eme 1}), se
recolle sur $X$ pour fournir l'expression intrins\`eque
 $$\omega^{k}_{\pi}={\cal H}om(\Omega^{n-k}_{X/S}, \omega^{n}_{\pi})$$
En reprenant les id\'ees du {\it lemme 3}, il apparait clairement
que $\omega^{k}_{\pi}$ est  de profondeur au moins deux fibre par
fibre puisque on obtient, pour tout ferm\'e $F$ de $X$,
l'isomorphisme
$${\cal H}^{i}_{F}({\cal H}om(\Omega^{n-k}_{X/S}, \omega^{n}_{\pi}))\simeq {\cal H}om(\Omega^{n-k}_{X/S},
{\cal H}^{i}_{F}( \omega^{n}_{\pi} )),\,\,\,\,{\rm
pour}\,\,i=0,\,1.$$ Il suffit alors d'invoquer le {\it th\'eor\`eme
1}, pour avoir  les annulations d\'esir\'ees pour tout ferm\'e $F$
de codimension au moins deux fibre par fibre sur $S$.
\par\noindent Dans le cas absolu, ceux sont bien les
dualis\'es de Andr\'eotti- Kas- Golovin des faisceaux
$\Omega^{n-k}_{X}$ (cf [KI], {\bf (3.0.1)},
p.61).$\,\,\,\,\,\blacksquare$\bigskip\bigskip\noindent
\centerline{{\tite  II. La preuve du th\'eor\`eme
2.}}\smallskip\bigskip\noindent {\bf 2.0. Platitude g\'eom\'etrique
analytique
 et int\'egration globale.}\smallskip\noindent La preuve de ce th\'eor\`eme
   va se faire en scindant l'\'etude en deux parties donnant,
 chacune, une caract\'erisation relativement simple de la platitude
 g\'eom\'etrique analytique (ou continue). Dans la premi\`ere, elle
  s'\'enoncera en terme de morphisme d'int\'egration sur les fibres avec
  de bonnes propri\'et\'es fonctorielles en les arguments; la seconde
  \'etablira une connexion \'etroite avec le faisceau du
  {\it th\'eor\`eme 1} par le biais d'un morphisme canonique que l'on appelle
   {\it morphisme classe fondamentale relative  } de $\pi$
   poss\'edant un certain nombre de propri\'et\'es
   fonctorielles parmi lesquelles la fameuse propri\'et\'e
   de la trace relative. \bigskip\noindent
{\bf 2.0.1. Les propositions fondamentales.} \Prop{1}{} {\it{Cas
global.}}\par\noindent Soit  $\pi$ un \'el\'ement de ${\cal E}(S,n)$
muni d'une pond\'eration ${\goth X}$. Alors, on a les
\'equivalences\par\noindent {\bf(i)} $\pi$ analytiquement
g\'eom\'etriquement plat\smallskip\noindent {\bf(ii)} il existe un
unique morphisme ${\cal O}_{S}$-lin\'eaire continu
\footnote{$^{(3)}$}{On entend par l\`a, une continuit\'e topologique
au niveau des  groupes  des  sections globales munis de leurs
structures d'espaces vectoriels topologiques {\bf{F.S}} ou
{\bf{Q.F.S}}} $\ds{{\int}_{\pi, {\goth X}}:{\rm I}\!{\rm
R}^{n}\pi_{!}\Omega^{n}_{X/S}\rightarrow {\cal O}_{S}}$
v\'erifiant:\par $\bullet$ il est stable par  changement de base
entre espace complexes r\'eduits (donnant, {\indent en particulier},
l'int\'egration usuelle dans le cas o\`u $S$ est un point).\par
$\bullet$ il est compatible  \`a l'additivit\'e des pond\'erations
sur $\pi$ et aux inclusions ouvertes sur $X$.\rm\smallskip\noindent
Cette assertion  repose essentiellement sur une version locale
donn\'ee par la \Prop{2}{}{\it{Cas local.}}\par\noindent
 Soient $X$, $Z$ et $S$ trois espaces analytiques
complexes r\'eduits munis d'un diagramme commutatif
$\xymatrix{X\ar@/_/[rr]_{\pi}\ar[r]^{\sigma}&Z\ar[r]^{q}&S}$ dans
lequel $\sigma$ est un plongement  dans $Z$ lisse et de dimension
relative $n+p$ sur $S$, $\pi$ un morphisme
 universellement $n$-\'equidimensionnel muni d'une pond\'eration ${\goth
 X}$. Alors $\pi$ est analytiquement g\'eom\'etriquement plat si et
seulement si il existe un unique morphisme ${\cal O}_{S}$-lin\'eaire
continue $\ds{{\int}_{\pi,{\goth X}}:{\rm I}\!{\rm
R}^{n}\pi_{!}\Omega^{n}_{X/S}\rightarrow {\cal O}_{S}}$ v\'erifiant
les propri\'et\'es de la {\it proposition
1}.\rm\smallskip\bigskip\noindent {\bf 2.0.2. Remarque.} On peut attirer l'attention du lecteur sur le fait que la correspondance donn\'ee par ${\int}_{\pi, {\goth X}}$ n'est pas
compatible \`a la composition des morphismes $\pi$ puisque la
compos\'ee de deux morphismes analytiquement g\'eom\'etriquement
plats n'est pas analytiquement g\'eom\'etriquement plate en
g\'en\'eral.\smallskip\par\par\noindent
{\bf 2.0.3. Fibres de dimension
nulle.}\smallskip\noindent
Pour mieux \'eclairer le lecteur, nous allons commencer par le cas particulier
  d'un morphisme fini. La situation est, certes,
restrictive (du moins en apparence!) mais suffisamment instructive
pour que l'on s'y'attarde un peu pour mettre en \'evidence la
strat\'egie et les ingr\'edients fondamentaux dont l'adaptation au
cas g\'en\'eral est \`a la base de la preuve du {\it th\'eor\`eme
2}. Pour cela, rappelons rapidement  quelques notions et faits
rencontr\'ees dans le \S1 de  [KI].\par $\star$ un morphisme propre
universellement \'equidimensionnel
  $\pi:X\rightarrow S$ est
  {\it analytiquement g\'eom\'etriquement plat} si et seulement s'il
existe une application analytique (qui est un plongement)
$\ds{\eta:S\rightarrow {\cal B}(X)}$ qui \`a $s\in S$ associe un
cycle de $X$ dont le support coincide avec la fibre ensembliste
$\pi^{-1}(s)$ et que, par construction de l'espace des cycles
analytiques effectifs ([B1]), ${\cal B}(X):= {\cal  B}(X_{red})$; de
plus, dans le cas particulier (et fondamental!) des $0$-cycles, on a
$\ds{{\cal B}_{0}(X)=\coprod_{r\geq 0}{\rm Sym}^{r}(X)}$ (cf [KI],
{\bf (1.0.2.3)}, p.21) .
\par
$\star$  un morphisme fini et surjectif $\pi:X\rightarrow S$
d\'efinit, g\'en\'eriquement sur $S$, un v\'eritable rev\^etement
ramifi\'e. Plus pr\'ecisemment, il existe un ouvert dense $S_{0}$,
que l'on peut supposer \^etre l'ensemble des points normaux de $S$,
au dessus duquel $\pi$ est  un rev\^etement ramifi\'e d'un certain
degr\'e $k$ dont on notera $(f_{j}(s))_{1\leq j\leq k}$ les branches
locales et $\ds{f: S_{0}\rightarrow{\rm Sym}^{k}(X)}$ l'application
analytique classifiante associ\'ee envoyant chaque point $s$ de
$S_{0}$ sur les fonctions sym\'etriques de ces branches locales
$({\rm s}_{j}(f_{1}(s),\cdots,f_{k}(s)))_{1\leq j\leq k}$. Alors, le
morphisme $\pi_{0}:X_{0}\rightarrow S_{0}$ d\'eduit de $\pi$ dans le
changement de base $S_{0}\rightarrow S$, d\'efinit une famille
continue de $0$-cycles puisque $S_{0}$ est localement irr\'eductible
(rappelons que les notions de finitude et de surjectivit\'e sont
stables par changement de base).  Mais  $S_{0}$ \'etant normal (donc
faiblement normal), cette famille est analytique.  Tout comme il a
\'et\'e dit dans [KI], {\bf (1.0.3.1)}, p.22, on dispose, pour
$\pi$, d'une image directe d\'efinie au sens des courants et d'un
op\'erateur {\it trace} associant \`a toute fonction holomorphe $g$
sur $X$, l'expression $\ds{{\cal
T}^{0}_{\pi}(g):=\sum^{j=k}_{j=1}g(s, f_{j}(s))}$ dont on notera
$\ds{{\cal T}^{0}_{\pi_{0}}(g)}$ sa restriction \`a $X_{0}$. Mais
$S_{0}$ \'etant normal, un courant $\bar\partial$-ferm\'e de type
$(0,0)$ d\'efinit naturellement une fonction holomorphe en vertu du
principe de prolongement de Hartogs ou du lemme de Dolbeault-
Grothendieck et par cons\'equent $\ds{{\cal T}^{0}_{\pi_{0}}(g)}$
est holomorphe sur $S_{0}$. A noter qu'en pr\'esence de
singularit\'es arbitraires, ce lemme n'est pas valable et il nous
est pas possible de statuer sur l'\'eventuelle variation analytique
de cette expression sans condition sur $\pi$.\Prop{3}{}  Soit
$\pi:X\rightarrow S$ un morphisme universellement
$0$-\'equidimensionnel muni d'une pond\'eration $\goth X$. Alors
$\pi$ est analytiquement g\'eom\'etriquement plat si et seulement si
il existe un morphisme trace ${\cal O}_{S}$-lin\'eaire continue
$\ds{{\cal T}^{0}_{\pi, \goth X}: \pi_{*}{\cal O}_{X}\rightarrow
{\cal O}_{S}}$, compatible aux changement de bases entre espaces
complexes r\'eduits quelconques et \`a l'additivit\'e des
pond\'erations.\rm\dem Avec les notations et remarques
pr\'ec\'edentes, on constate que  la platitude g\'eom\'etrique
analytique de $\pi$ est strictement \'equivalente au prolongement
analytique de  $f$ \`a $S$ tout entier ou, ce qui revient au m\^eme,
au fait que  le courant $\bar\partial$-ferm\'e sur $S$ donn\'e par
$\ds{\sum^{j=k}_{j=1}g(s, f_{j}(s))}$ est en fait une fonction
holomorphe sur $S$ ! Il est clair que  la finitude de $\pi$, la
notion de platitude g\'eom\'etrique et la construction d'une telle
trace conf\`erent \`a notre  probl\`eme une  nature locale sur $S$
et $X$. En se fixant un point $s_{0}$ de $S$ et un voisinage ouvert
de Stein (que l'on r\'etr\'ecit \`a souhait!), on se ram\`ene \`a la
factorisation locale
$\xymatrix{X\ar@/_/[rr]_{\pi}\ar[r]^{\sigma}&S\times V\ar[r]^{q}&S}$
dans laquelle $X$ et $S$ d\'esignent abusivement des ouverts de
localisations,  $\sigma$ un  $S$-plongement local avec $V$
polydisque  ouvert d'un espace num\'erique  ${\Bbb C}^{p}$ et $q$ la
projection canonique. Enfin, on consid\'erera
$\ds{\tilde{\sigma}:{\rm Sym}^{k}(V)\rightarrow
W:=\bigoplus^{k}_{h=1} S_{h}(V)}$ le plongement naturel (cf [KI],
{\bf (1.0.1.1)}, p.14) induit par passage au quotient de la somme
directe des applications $\ds{s_{h}:({\Bbb C}^{p})^{k}\rightarrow
S_{h}({\Bbb C}^{p})}$ d\'efinies par $\ds{s_{h}(X_{1},\cdots,
X_{k}):=\sum_{1\leq i_{1}<\cdots<i_{h}\leq h}X_{i_{1}}\cdots
X_{i_{h}}}$ pour $h\in\{1,\cdots,k\}$ et permettant d'identifier
$\ds{\rm Sym}^{k}_{S}(X)$ \`a un sous espace de $W$.\par
$\blacklozenge$ Supposons  $\pi$ analytiquement g\'eom\'etriquement
plat. Alors,  il s'installe naturellement dans le diagramme
commutatif
$$\xymatrix{&S\times V\ar[r]^{F}\ar[dd]^{q}&{\rm Sym}^{k}(V)\times V\ar[dd]^{\tilde{q}}\\
X\ar[ru]\ar[rd]_{\pi}\ar[rrr]&&&{\rm Sym}^{k}(V)\# V\ar[ld]^{\pi_{\#}}\ar[lu]\\
&S\ar[r]_{f}\,\,\,\,\,&{\rm Sym}^{k}(V)}$$ dans lequel  $F:=f\times
Id_{S}$, $q$ (resp. $\tilde{q}$) est la projection canonique, ${\rm
Sym}^{k}(V)\# V$ d\'esigne le sous espace d'incidence et  $\pi_{\#}$
le morphisme  naturellement d\'eduit de $\pi$.\par\noindent
Remarquons que $X$ s'identifie au graphe de $\pi$ quotient\'e par le
groupe
 sym\'etrique d'ordre $k$.\par\noindent
Convenons de noter ${\goth  T}^{0}_{\# }$ la trace universelle
relative au
 morphisme universel ${\rm Sym}^{k}(V)\# V\rightarrow {\rm Sym}^{k}(V)$ qui
 d\'efinit une famille analytique de points puisqu'il est \'equidimensionnel
  sur une base normale !\par\noindent
Par d\'efinition (cf [KI], {\bf (1.0.3.1)}, p.23), cette trace
universelle g\'en\`ere la trace recherch\'ee puisqu'elles sont
li\'ees par la formule de changement de base $f^{*}({\goth
T}^{0}_{\# })(-)= {\cal T}^{0}_{\pi,\goth{X}}F^{*}(-)$ qui nous
montre que l'analyticit\'e de $f$ entraine automatiquement
l'analyticit\'e de notre trace sur les fonctions
${\sigma}_{k}$-invariantes. Comme  de telles  fonctions sont
engendr\'ees par les fonctions de Newton, le probl\`eme est r\'eduit
au seul fait de savoir si  l'action de cette trace sur les fonctions
de Newton suffit \`a la caract\'eriser compl\`etement pour pouvoir
conclure.\par Il en est effectivement ainsi. En effet, d\'esignons
par $x:=(x_{1},\cdots,x_{p})$ le point courant de $V$ et $s$ une
coordonn\'ee sur $S$ que l'on peut supposer \^etre un ouvert de
Stein. Alors,  il est bien connu que, par d\'ecomposition
nucl\'eaire, toute fonction holomorphe sur $S\times V$ s'\'ecrit
sous la forme $\ds{\sum_{m}\sum_{|I|=m; m\geq 0}a_{I}(s)x^{I}}$ avec
$a_{I}(s)$ holomorphes sur $S$. Mais la ${\cal O}_{S}$-lin\'earit\'e
de cette trace nous ram\`ene principalement \`a \'etudier son action
sur les fonctions de Newton $x^{I}$. Alors, $f$ \'etant  holomorphe,
essentiellement  par d\'efinition, les fonctions sym\'etriques des
branches locales sont holomorphes et par suite, leurs  fonctions de
Newton aussi. On en d\'eduit, alors, l'analyticit\'e en $s$  de
toutes les expressions du type $\ds{{\rm
N}_{m}(f(s)):=\sum^{j=k}_{|I|=m; j=1}f^{I}_{j}(s)}$  et par suite
l'existence d'un unique morphisme ${\cal O}_{S}$-lin\'eaire
$\ds{{\cal T}^{0}_{\pi,\goth{X} }: \pi_{*}{\cal O}_{X}\rightarrow
{\cal O}_{S}}$. Ce morphisme  prolonge naturellement $\ds{{\cal
T}^{0}_{\pi_{0}}}$ puisque  $\ds{{\cal T}^{0}_{\pi_{0}}(x^I):= {\cal
T}^{0}_{\pi,\goth{X}}(x^I)|_{X_{0}}:={\rm N}_{m}(f(s))}$ et il est
de plus continue relativement aux structures d'espaces vectoriels
topologiques {\bf F.S}  dont on munit naturellement $\Gamma(X,{\cal
O}_{X})$ et $\Gamma(S,{\cal O}_{S})$. On se ram\`ene essentiellement
\`a un probl\`eme sur des alg\`ebres analytiques complexes munis de
leurs topologies canoniques (cf [Ju]) et \`a utiliser des
in\'egalit\'es de type Cauchy (cf [B1], lemme 5).
\par\noindent
La compatibilit\'e avec l'additivit\'e des pond\'erations est
une donn\'ee inh\'erente \`a un morphisme trace agissant sur les
familles continues de cycles!
 \par $\blacklozenge$
R\'eciproquement si, pour toute fonction holomorphe $g$ sur $X$,
 $\ds{{\cal T}^{0}_{\pi,\goth{X}}(g)}$ est holomorphe sur $S$, on a, en
particulier,$\ds{{\cal T}^{0}_{\pi,\goth{X}}(x^{I})}$ holomorphe sur $S$.
Bien s\^ur, le proc\'ed\'e de sym\'etrisation nous permet de nous
concentrer sur les fonctions ${\sigma}_{k}$-invariantes m\^eme si ce
n'est pas explicitement dit.  Or sur $S_{0}$, les expressions de
Newton $\ds{{\rm N}_{m}(f(s))}$ sont holomorphes et coincident
naturellement avec $\ds{{\cal T}^{0}_{\pi,\goth{X}}(x^{I})}$. On en d\'eduit
que les fonctions de Newton se prolongent analytiquement sur $S$
tout entier, et bien entendu, il en va de m\^eme pour les fonctions
sym\'etriques des branches locales puisqu' il est bien connu  qu'en
caract\'eristique nulle les fonctions de Newton engendrent
l'alg\`ebre des fonctions sym\'etriques. D'o\`u une application
holomorphe classifiante $\ds{\Psi_{\pi}:S\rightarrow {\cal
B}_{0}(X)}$ et donc la platitude g\'eom\'etrique analytique  de
$\pi$. $\,\,\blacksquare$
\smallskip\par\noindent
{\bf 2.0.4. Remarques.}\par\noindent {\bf(i)} Les  exemples simples donn\'es  dans [KI],{\bf (1.3.6)},p.41) \`a d\'efaut d'un mode de variance
raisonnable sur les fibres, il faut se garder de croire
 que l'on dispose naturellement d'une  {\it{trace}} $\pi^{*}\pi_{*}{\cal O}_{S}\rightarrow
{\cal O}_{S}$.  \par\noindent {\bf(ii)} Si $\pi$ est fini et
g\'eom\'etriquement plat, la trace d\'esir\'ee est, en fait, un cas
particulier de l'int\'egration sur les familles analytiques de
cycles ([BV], [K2]); la famille de $0$-cycles est induite par les
fibres de $\pi$ munies des multiplicit\'es
convenables.\smallskip\noindent {\bf 2.0.5. Fibres de dimension
strictement positive.} \par\noindent Dans cette situation, l'id\'ee
fondamentale est de se raccrocher au cas pr\'ec\'edent par le biais
des  factorisations locales du morphisme. Mais \`a l'instar du cas
des fibres $0$-dimensionnelles (ou de bases particuli\`eres telles
que faiblement normales), on ne peut se contenter de tester au moyen
de traces de fonctions holomorphes cette platitude g\'eom\'etrique
dans {\bf
 une} seule installation locale (hors-mis des bases particuli\`eres telles que faiblement
normales, par exemple). Par ailleurs, pour des raisons que l'on peut
facilement comprendre, il est tout-\`a-fait naturel de penser le
test de la trace sur les formes diff\'erentielles relatives  de
degr\'e maximum sur les fibres. Le contre exemple de [B1] propos\'e
dans [KI], {\bf (1.1)}, p.30,  montre bien que, relativement \`a
{\bf une} factorisation locale, on peut tr\`es bien avoir une trace
$f_{*}{\cal O}_{X}\rightarrow {\cal O}_{Y}$ sans avoir la trace
fondamentale $f_{*}\Omega^{n}_{X/S}\rightarrow \Omega^{n}_{Y/S}$. Le
passage au cas g\'en\'eral ne requiert pas de strat\'egie
particuli\`ere mais seulement d'adapter les outils en rempla\c cant,
par exemple, les fonctions de Newton par les formes de Newton (cf
[KI],{\bf (1.0.1.3)}, p.16).
\par\noindent
\Prop{4}{}\par\noindent Soit $\pi\in
{\cal E}(S,n)$ pond\'er\'e par un cycle $\goth{X}$. Alors $\pi\in
{\cal G}_{a}(S,n)$ si et seulement si
 pour tout
espace complexe $S$-lisse $Y$ et tout morphisme fini et surjectif
$f:X\rightarrow Y$, il existe un unique morphisme {\it{trace}}
${\cal O}_{S}$-lin\'eaire continu
\par\noindent \centerline{${\cal
T}^{n}_{f,\goth X}:f_{*}\Omega^{n}_{X/S}\rightarrow
\Omega^{n}_{Y/S}$}\smallskip\noindent prolongeant  la trace
usuelle\par\noindent
\centerline{$f_{*}f^{*}\Omega^{n}_{Y/S}\rightarrow
\Omega^{n}_{Y/S}$}\smallskip\noindent compatible \`a tout changement
de base sur $S$, aux localisations sur $X$, \`a l'additivit\'e des
pond\'erations sur $\pi$ et  induisant (et est induit par)  le
morphisme d'alg\`ebres diff\'erentielles gradu\'ees\par\noindent
\centerline{${\cal T}^{\bullet}_{f,\goth X}
:f_{*}\Omega^{\bullet}_{X/S}\rightarrow
\Omega^{\bullet}_{Y/S}$}.\rm\smallskip\noindent \dem \par
$\blacklozenge$ Supposons  $\pi$  analytiquement
g\'eom\'etriquement plat. Alors, \'etant universellement $n$-\'equidimensionnel, pour tout point $s_{0}$ de
$S$, il existe voisinage ouvert $S_{0}$ de $s_{0}$ dans $S$ et une
famille d'\'ecailles ou de cartes $(E_{\alpha}:=(U_{\alpha},
B_{\alpha}, \sigma_{\alpha}))_{\alpha\in A}$ $S_{0}$- adapt\'ees \`a
$X_{s_{0}}$ de sorte que, pour tout $\alpha\in A$,  on ait une
installation (locale) ($\clubsuit$)
$$\xymatrix{X_{\alpha}\ar[d]^{\pi_{\alpha}} \ar [dr]_{f_{\alpha} }\ar [r]^{\sigma_{\alpha}}
 & S_{0}\times U_{\alpha}\times B_{\alpha} \ar[d]^{p_{\alpha}}\\
 S_{0}&\ar[l]^{q_{\alpha}} S_{0}\times U_{\alpha}}$$
dans laquelle $f_{\alpha}$ est un morphisme ouvert, fini et
surjectif, $\sigma_{\alpha}$ un plongement local, $p_{\alpha}$ et
$q_{\alpha}$ \'etant les projections canoniques. De telles cartes,
que l'on peut choisir centr\'ees par un rapport \`a un point
privil\'egi\'e sur la fibre, existent toujours (cf [A.S]). $\pi$
\'etant analytiquement g\'eom\'etriquement plat, les morphismes
 $f_{\alpha}$ le seront aussi  et d\'efiniront donc des familles analytiques de points. De plus,  cette factorisation locale est toujours munie de deux applications  analytiques\par
$\bullet$ $\ds{F_{\alpha}: S_{0}\times U_{\alpha}\rightarrow {\rm
Sym}^{k}(B_{\alpha})},\,\,{\rm et}$\par $\bullet$
$\ds{T_{i,j}(F_{\alpha}(s,t)):(s, t)\rightarrow
\sum^{k}_{l=1}(f_{l}(s,t))^{i}\otimes\Lambda^{j}(D_{t}f_{l}(s,t))}$\smallskip\noindent
Cette derni\`ere expression \'etant d\'efinie localement et
seulement g\'en\'eriquement holomorphe sans la condition de
platitude g\'eom\'etrique!\par\noindent
 Les morphismes $F_{\alpha}$
appel\'es {\it{classifiants}} font la connexion avec la situation
universelle d\'ecrite par le diagramme commutatif
$$\xymatrix{&S_{0}\times U_{\alpha}\times B_{\alpha}\ar[r]^{F_{\alpha}\times Id_{B_{\alpha}}}\ar[dd]&{\rm Sym}^{k}(B_{\alpha})
 \times B_{\alpha}\ar[dd]&\\
X_{\alpha}\ar[ur]\ar[rd]_{f_{\alpha}}\ar[rrr]&&& {\rm
Sym}^{k}(B_{\alpha}) \# B_{\alpha}\ar[ld]^{\pi_{\#}}\ar[ul]
\\&S_{0}\times U_{\alpha}\ar[r]_{F_{\alpha}}&{\rm Sym}^{k}(B_{\alpha})&}$$
dans lequel, on a ensemblistement $\ds{(F_{\alpha}\times
Id_{B_{\alpha}})^{-1}({\rm Sym}^{k}(B_{\alpha})\# B_{\alpha}) =
X}$.\par\noindent On va montrer que, pour chaque $\alpha$ fix\'e, la
platitude g\'eom\'etrique analytique  donne naissance \`a une
collection de morphismes traces ${\cal T}^{n}_{f_{\alpha},
\goth{X_{\alpha}}}:{f_{\alpha}}_{*}\Omega^{n}_{X_{\alpha}/S}\rightarrow
\Omega^{n}_{{S\times U_{\alpha}}/S}$ v\'erifiant les propri\'et\'es
requises avec $\goth{X_{\alpha}}:={\goth {X}}\bullet(S\times
X_{\alpha})$. On en d\'eduira ais\'ement qu'il en est ainsi pour
tout morphisme fini et surjectif $f:X\rightarrow Y$ de $X$ sur un
espace complexe $Y$ lisse sur $S$ et de dimension relative $n$.
\par\noindent
 D'apr\`es ce qui pr\'ec\`ede, il nous est facile de voir que l'assertion du
lemme se ram\`ene \`a l'\'equivalence:
$$\left [\matrix{\forall\,\alpha,\,\,\forall \,\,{\rm j}\in\{1,\cdots,p\}\,\,\,T_{i,j}(F_{\alpha}(s,t))\,\,\,\,{\rm  S-
holomorphe}\hfill\cr\noalign{\vskip 0,2cm}
\,\,\,\,\,\,\,\,\,\,\,\,\,\,\,\,\,\,\,\,\,{\hskip
2cm}{\Updownarrow}\hfill\cr {\cal T}^{n}_{f_{\alpha},
\goth{X_{\alpha}}}(\xi):={\ds
\sum^{l=k}_{l=1}}f^{*}_{\alpha,l}(\xi)\in
\Gamma(V,\Omega^{n}_{Y/S}),\,\,\,\forall\,\,\xi\in\Gamma(U,\Omega^{n}_{X/S})}\right
]$$\par\noindent
 Afin de ne pas  trop alourdir le texte, nous omettrons
d\'elib\'eremment les indices en nous fixant une \'ecaille et on
pose $\goth{X}:=\goth{X_{\alpha}}$, $Y:= S\times U$, $Z=S\times
U\times B$, $f:=f_{\alpha}$, $x:=(x_{1},\cdots, x_{p})$ le point
courant de $B$ dans ${\Bbb C}^{p}$ et  $w_{i,j}:= x^{I}dx^{J}$ , $I$
et $J$ \'etant des ensembles d'indices de longueur $i$ et $j$
respectivement. Par ailleurs, nous utiliserons, sans le mentionner,
le plongement  $\ds{{\rm Sym}^{k}(B)\rightarrow
\bigoplus_{\atop{1\leq m\leq k}}{\rm S}_{m}({\Bbb C}^{p}),}$
permettant  un "rep\'erage" des fonctions sym\'etriques en termes de
fonctions de Newton.\par\noindent Comme il a \'et\'e mentionner dans
[KI],{\bf (1.0.3.1.1)}, p.23, la situation universelle est idyllique
puisque le morphisme $\pi_{\#}:\#\rightarrow {\rm Sym}^{k}(B)$ est
analytiquement $0$-g\'eom\'etriquement plat et dot\'e  d'un
morphisme  d'alg\`ebres diff\'erentielles gradu\'ees
$${\goth T}_{\#}^{\bullet}:({\pi_{\#}})_{*}\omega^{\bullet}_{\#}\rightarrow
\omega^{\bullet}_ {{\rm Sym}^{k}({\Bbb C}^{p})}$$ dont les
composantes sont des morphismes traces. D'apr\`es
[KI],{\bf(1.0.3.1.1)}, p.24 et  {\bf(1.0.3.2.1)}),p.27,  ${\cal T}^{\bullet}_{f, \goth{X}}$ est
d\'efini gr\^ace \`a la relation
$${\cal T}^{\bullet}_{f, \goth{X}}(w_{i,j}):= F^{*}{\goth
T}^{\bullet}_{\#}(\tilde{w}_{i,j}),$$ dans laquelle
$\tilde{w}_{i,j}$ est une forme m\'eromorphe, section du faisceau
$\omega^{\bullet}_ {{\rm Sym}^{k}(B)\times B}$, naturellement
induite par la forme $\sigma_{k}$-invariante image r\'eciproque de
$w_{i,j}$ par les projections canoniques de $B^{k}$ sur ses
diff\'erents facteurs. Il est absolument \'evident que
pour les  formes de Newton \`a valeurs vectorielles
$\ds{W_{i,j}:=\sum_{|I|=i,\,|J|=j}x^{I}\otimes dx^{J}}$, on a
 $\ds{T_{i,j}(F(s,t))
=(F)^{*}{\goth T}^{\bullet}_{\#}(\tilde{W}_{i,j}) }$ et qu'il n'est
pas difficile de voir que ces expressions d\'efinissent (et sont
d\'efinies par) les formules pr\'ec\'edentes.
 On d\'efinit, ainsi, une trace
naturelle pour $f$ dont l'action sur les formes de Newton est $S$-
holomorphe sur $S_{0}\times U$. En d'autres termes $\ds{{\cal
T}^{n}_{f, \goth{X}}(x^{I}dx^{J})}$ est holomorphe sur $S_{0}\times
U$.\par\noindent Pour passer au cas d'une $n$-forme holomorphe
$S$-relative  g\'en\'erale sur $X$, il suffit d'invoquer les faits
standards suivants :\smallskip {\bf(a)} la
$\Omega^{\bullet}_{S\times U}$- lin\'earit\'e du morphisme trace
${\cal T}^{n}_{f,\goth{X}}$\par {\bf(b)} l'isomorphisme
$\ds{\Gamma(S'\times U\times B, \Omega^{r}_{S\times U\times
B/S})\simeq {\cal O}_{S}(S')\widehat{\bigotimes}_{{\Bbb
C}}\Omega^{r}_{U\times B} }$, donn\'e pour tout ouvert de Stein $S'$
de $S$, et  la d\'ecomposition standard $\ds{\Omega^{r}_{U\times
B}=\sum_{i+j=r}\Omega^{i}_{U}\otimes_{{\cal O}_{U\times
B}}\Omega^{j}_{B}}$ \par {\bf(c)} le fait qu'il y'ait suffisamment
de projections lin\'eaires permettant de d\'ecrire, en chaque point
de $X$, le faisceau des $n$- formes diff\'erentielles relatives sur
$X$.\par{\bf(d)} En tout degr\'e $k$, le faisceau $\omega^{\bullet}_
{{\rm Sym}^{k}({\Bbb C}^{p})}$ est engendr\'e par les formes
m\'eromorphes dont l'image r\'eciproque par l'application quotient
$q: ({\Bbb C}^{p})^{k}\rightarrow{\rm Sym}^{k}({\Bbb C}^{p})$ ( qui
se comporte d'ailleurs comme une r\'esolution \'equivariante) est
une forme de Newton (cf [KI],{\bf(1.0.1.3)}, p.16).
\smallskip\noindent {\bf(a)} et {\bf(b)} montrent, de toute
\'evidence, que seules les formes diff\'erentielles provenant de $B$
jouent un r\^ole essentiel dans la description de cette trace.
{\bf(c)} signifie qu'en  tout point $x\in X$, pour tout plongement
local $\sigma$ dans $Z$ lisse sur $Y$ (et donc toute projection $S$-
lin\'eaire $q:Z\rightarrow Y$) telle que $q o \sigma$ soit encore
fini en $x$,  la famille $(q^{*}(\Omega^{n}_{Y/S}))_{q}$ engendrent
$\Omega^{n}_{Z/S}$ au voisinage de $x$ et donc $(q o
\sigma)^{*}(\Omega^{n}_{Y/S})$ engendre $\Omega^{n}_{X/S}$ au
voisinage de $x$. C'est  un probl\`eme d'alg\`ebre lin\'eaire simple
auquel on renvoie le lecteur au  {\it lemme (5.1)} de [KI],\S{\bf
5},p.86.\par\noindent  On termine en constatant que  le passage des
projections lin\'eaires aux projections quelconques se fait sans
douleurs gr\^ace au {\it lemme (5.2)} de [KI],
p.86.\smallskip\noindent La d\'ecomposition nucl\'eaire ci-dessus
d\'ecoule d'arguments standards (Kiehl, Verdier, Forster, Knorr)
assurant que  toute $r$ forme $S$- relative $\xi$  sur $S\times
U\times B$ se d\'ecompose sous forme de s\'erie convergente
$$\xi:=\sum _{I,K;|K|+|I|=r}\alpha_{K}\wedge w_{i,j}$$
$\alpha_{K}$ \'etant une $(r-j)$- forme holomorphe $S$-relative sur
$S\times U$ ``born\'ee'' sur tout compact  de $U$, on en d\'eduit
{\bf(ii)}.\par\noindent La continuit\'e est facile \`a \'etablir une
fois que l'on munit $\Gamma(X, \Omega^{n}_{X/S})$ et $\Gamma(Y,
\Omega^{n}_{Y/S})$ de structure topologique de type {\bf F.S} (cf
la {\it proposition 3}, p.18).\par\noindent
 On peut remarquer que l'existence de
la d\'ecomposition (nucl\'eaire) de $\xi$ cit\'ee plus haut, et la
$\Omega^{\bullet}_{S\times U}$- lin\'earit\'e du morphisme trace
donne imm\'ediatement l'analyticit\'e de la trace ${\cal
T}^{n}_{f,\goth{X}}(\xi)$.\smallskip\noindent $\bullet$ {\bf Compatibilit\'e
avec  la trace usuelle} $\ds{f_{*}f^{*}\Omega^{n}_{Y/S}\rightarrow
\Omega^{n}_{Y/S}}$.\par\noindent
Tout morphisme fini et analytiquement g\'eom\'etriquement plat  $\ds{f:X\rightarrow
Y}$, o\`u  $Y:=S\times U$, provenant d'une factorisation locale de $\pi$, induit, en vertu de la {\it proposition 3}, un morphisme trace  $\ds{{\cal T}^{0}_{f,\goth{X}}:f_{*}{\cal
O}_{X}\rightarrow {\cal O}_{Y}}$ . Mais l'isomorphisme $\ds{f_{*}f^{*}\Omega^{n}_{Y/S}\simeq
f_{*}{\cal O}_{X}\otimes_{{\cal O}_{Y}}\Omega^{n}_{Y/S}},$ permet de d\'efinir
 un morphisme trace (que l'on appelle ``naturel'')  ${\cal O}_{Y}$-lin\'eaire\par\noindent
$\ds{\tilde{{\cal T}}^{n}_{f,\goth{X}}:={\cal T}^{0}_{f}\otimes Id_{\Omega^{n}_{Y}}:f_{*}f^{*}\Omega^{n}_{Y/S}\rightarrow
\Omega^{n}_{Y/S}}$ (donn\'e  $deg(f).Id$).\par\noindent Il est
\'evident que les morphismes  traces $\ds{{\cal
T}^{n}_{f,\goth{X}}:f_{*}\Omega^{n}_{X/S}\rightarrow \Omega^{n}_{Y/S}}$
d\'efinis pr\'ec\'edemment prolongent naturellement les morphismes $\tilde{{\cal T}}^{n}_{f,\goth{X}}$  en s'ins\'erant dans le  diagramme commutatif
$$\xymatrix{f_{*}f^{*}\Omega^{n}_{Y/S}\ar[r]\ar[d]_{\tilde{{\cal T}}^{n}_{f,\goth{X}}}&f_{*}
\Omega^{n}_{X/S}\ar[d]^{{\cal T}^{n}_{f,\goth{X}}}\\
\Omega^{n}_{Y/S}\ar[r]_{Id}&\Omega^{n}_{Y/S}}$$ \par\noindent
$\bullet$ {\bf Compatibilit\'e aux changements  de bases, aux localisations sur $X$  et passage
au cas} $Y$ {\bf  lisse sur} $S$.\par\noindent
La compatibilit\'e \`a la localisation sur $X$ revient \`a voir dans quelle
mesure des  installations locales (d'intersection non vide!) autour d'une fibre donn\'ee produiront elles des traces qui se correspondent par changement de base dans toutes les directions de la param\'etrisation locales et non seulement dans celle de $S$. Cela nous ram\`ene principalement \`a examiner le cas de deux param\'etrisations $f:X\rightarrow Y:=S\times U$ et $f':X'\rightarrow Y':=S\times U'$ provenant de factorisations locales de $\pi$ et dans lesquelles  $U'\Subset U\subset {\Bbb C}^{n}$. Le recollement d\'esir\'e se traduit, alors, en exhibant ${\cal T}^{n}_{f,\goth{X}}$ comme un prolongement naturel de ${\cal T}^{n}_{f',\goth{X}}$. Mais ceci est la cl\'e de vo\^ute de la construction de [B1].\par\noindent
La compatibilit\'e aux changement de base en l'argument  $S$ est facile \`a \'etablir. En effet,  comme les notions de finitude, surjectivit\'e, lissit\'e et platitude g\'eom\'etrique analytique sont stables  par changement de base
quelconques entre espaces complexes r\'eduits, il est facile de voir qu\'etant donn\'e un morphisme d'espaces complexes r\'eduits $g:S'\rightarrow S$ et
$$\xymatrix{X'\ar[r]^{h}\ar[d]_{f'}&X\ar[d]^{f}\\
S'\times U\ar[r]_{G:=g\otimes Id_{U}}&S\times U}$$
le  diagramme cart\'esien de changement de base  qui lui est associ\'e, on obtient un diagramme commutatif
 $$\xymatrix{f'_{*}\Omega^{n}_{X'/S'}\ar[d]_{{\cal T}^{n}_{f'}}&f'_{*}h^{*}\Omega^{n}_{X/S}\simeq g^{*}f_{*}\Omega^{n}_{X/S}\ar[l]^{\!\!\!\!\!\simeq}\ar[d]^{g^{*}{\cal T}^{n}_{f}}\\
\Omega^{n}_{S'\times U/S'}&g^{*}\Omega^{n}_{S\times
U/S}\ar[l]_{\simeq}}$$
Le passage au cas d'un espace $S$-lisse quelconque
$Y$  ne pr\'esente aucune difficult\'e puisqu'en vertu de la nature
locale du probl\`eme, on se ram\`ene \`a la situation pr\'ec\'edente
dans laquelle  il est toujours possible d'installer $f:X\rightarrow
Y$ dans le diagramme
$$\xymatrix{X\ar[r]\ar[d]_{f}&S\times U\times B\ar[d]\\
Y\ar[r]_{h}&S\times U}$$ o\`u  $h$ est un morphisme  \'etale
c'est-\`a-dire un rev\^etement non ramifi\'e et  donc un
isomorphisme local. Alors, tenant compte de la surjectivit\'e du
morphisme  image r\'eciproque $h^{*}\Omega^{n}_{{S\times
U}/S}\rightarrow \Omega^{n}_{Y/S}$, il est clair que la trace
donn\'ee par le morphisme  $hof$
 donnera, sans difficult\'e, une trace pour $f$.
\smallskip\noindent
$\bullet$ {\bf Compatibilit\'e avec la structure  d'alg\`ebre
gradu\'ee.}\par\noindent
 Pour tout entier naturel  $k$ fix\'e, le morphisme canonique $\Omega^{k}_{X/S}\rightarrow {\cal
H}om(\Omega^{n-k}_{X/S}, \Omega^{n}_{X/S})$  auquel on applique le
morphisme fini $f$, donne  la  fl\`eche
$$f_{*}\Omega^{k}_{X/S}\rightarrow f_{*}{\cal H}om(\Omega^{n-k}_{X/S}, \Omega^{n}_{X/S})\simeq {\cal H}om(f_{*}\Omega^{n-k}_{X/S}, f_{*}\Omega^{n}_{X/S})$$
Il suffit alors d'utiliser l'image r\'eciproque
$\ds{\Omega^{n-k}_{Y/S}\rightarrow f_{*}\Omega^{n-k}_{X/S}}$ et la
trace pr\'ec\'edente $\ds{{\cal T}^{n}_{f,\goth{X}}:
f_{*}\Omega^{n}_{X/S}\rightarrow \Omega^{n}_{Y/S}}$ pour en
d\'eduire le morphisme trace $\ds{{\cal T}^{k}_{f,\goth{X}}:
f_{*}\Omega^{k}_{X/S}\rightarrow \Omega^{k}_{Y/S}}$. Alors, en d\'esignant par $d_{Y/S}$ la diff\'erentielle ext\'erieure relative usuelle  sur $Y$, il est facile d'en d\'eduire un  morphisme de $(\Omega^{\bullet}_{Y/S}, d)$- alg\`ebres diff\'erentielles  gradu\'ees \smallskip \centerline{${\cal
T}^{\bullet}_{f,\goth{X}}: f_{*}\Omega^{\bullet}_{X/S}\rightarrow
\Omega^{\bullet}_{Y/S}$}\par\noindent
Une autre mani\`ere d'exprimer ce fait serait de dire, qu'au niveau des morphismes de degr\'e $0$, on a \par
 \centerline{$ {\rm
H}om^{0}_{\Omega^{\bullet}_{Y/S}}(f_{*}\Omega^{\bullet}_{X/S},
\Omega^{\bullet}_{Y/S}) \simeq {\rm H}om_{{\cal O}_{Y}}(
f_{*}\Omega^{n}_{X/S},\Omega^{n}_{Y/S}).$}
\noindent
puisque, $Y$ \'etant lisse sur $S$,
l'homomorphisme $\ds{\Omega^{k}_{Y/S}\rightarrow {\rm H}om_{{\cal
O}_{Y}}(\Omega^{n-k}_{Y/S}, \Omega^{n}_{Y/S})}$ qui, \`a $\xi$,
associe le morphisme $\ds{h_{\xi}(\phi)= \phi\wedge\xi}$, donne les
isomorphismes
$$\ds{\Omega^{\bullet}_{Y/S}\simeq {\rm H}om_{{\cal
O}_{Y}}(\Omega^{\bullet}_{Y/S}, \Omega^{n}_{Y/S})},\,\,{\rm et}\,\,
\ds{\Omega^{\bullet}_{Y/S}\simeq {\rm H}om_{{\cal O}_{Y}}({\rm
H}om_{{\cal O}_{Y}}(\Omega^{\bullet}_{Y/S}, \Omega^{n}_{Y/S}),
\Omega^{n}_{Y/S})}$$ desquels r\'esulte
$$ {\rm H}om_{\Omega^{\bullet}_{Y/S}}(f_{*}\Omega^{\bullet}_{X/S}, \Omega^{\bullet}_{Y/S})\simeq  {\rm H}om_{\Omega^{\bullet}_{Y/S}}(f_{*}\Omega^{\bullet}_{X/S},Hom_{{\cal O}_{Y}}( \Omega^{\bullet}_{Y/S},\Omega^{n}_{Y/S}))\simeq Hom_{{\cal O}_{Y}}( f_{*}\Omega^{\bullet}_{X/S},\Omega^{n}_{Y/S})$$
$\bullet${\bf Compatibilit\'e avec  l'addition des
pond\'erations.}\par\noindent Soient  ${\goth X}_{1}$ et ${\goth
X}_{2}$ deux pond\'erations (pas n\'ecessairement analytiquement
g\'eom\'etriquement plate!)  de $\pi$. Si  l'on d\'efinit   la somme
de deux familles de cycles en prenant les sommes point par point des
cycles correspondants et l'on utilise l'additivit\'e de la trace
dans le cas absolu\footnote{$^{(5)}$}{Si $X_{1}$ et $X_{2}$ sont
deux cycles install\'es localement dans une m\^eme \'ecaille
adapt\'ee avec $f_{1}$ et $f_{2}$ pour morphismes finis respectifs,
alors $\ds{{\cal T}^{\bullet}_{f_{1}\oplus f_{2}}= {\cal
T}^{\bullet}_{f_{1}} + {\cal T}^{\bullet}_{f_{2}}}$} il est facile
de constater que dans toute \'ecaille  $(U,B,\sigma)$-$S$-adapt\'ee,
on a
$${\cal
T}^{\bullet}_{f,{\goth X}_{1}\oplus{\goth X}_{2} }= {\cal
T}^{\bullet}_{f,{\goth X}_{1}}+ {\cal T}^{\bullet}_{f,{\goth
X}_{2}}$$
\par\noindent
$\blacklozenge$ La r\'eciproque est facile \`a \'etablir, vu ce qui
a \'et\'e dit pr\'ec\'edemment sur le lien entre les morphismes
traces et les $\ds{T_{i,j}(F_{E}(s,t))}$. En effet, il est clair que
l'analyticit\'e des traces de formes holomorphes $S$-relatives
quelconques sur $X$ entraine automatiquement l'analyticit\'e des
traces de formes de Newton , donc l'analyticit\'e des applications
$\ds{T_{i,j}(F_{E}(s,t))}$ dans toute \'ecaille $E$ et par suite
l'isotropie de $F_{E}(s,t)$. Il suffit simplement d'utiliser le
proc\'ed\'e de sym\'etrisation des formes pour obtenir des sections
du faisceau des formes ``sym\'etriques'' sur ${\rm Sym}^{k}(B)$ qui,
via le plongement ``en Newton'' $\ds{{\rm Sym}^{k}(B)\rightarrow
\bigoplus_{1\leq m\leq k}{\rm S}_{m}({\Bbb C}^{p}),}$ se d\'ecrivent
naturellement en fonction de formes de Newton
$\,\,\blacksquare$\smallskip\par\noindent
 {\bf 2.0.6. Preuve de la proposition
2.}
\smallskip
$\blacklozenge$ Supposons $\pi$ est analytiquement  g\'eom\'etriquement plat. Comme il est, en particulier, universellement $n$-
\'equidimensionnel, on peut toujours ins\'erer ces donn\'ees locales
dans un  diagramme commutatif ($\spadesuit$)
$$\xymatrix{X\ar[rdd]_{\pi}\ar[rr]^{\sigma}\ar[rd]^{f}&&Z\ar[ld]_{q'}\ar[ldd]^{r}\\
&Y\ar[d]_{q}&\\
&S&}$$ avec $f$ fini et surjectif, $Z:=S\times U\times B$,
$Y:=S\times U$, $U$ (resp. $B$) \'etant un polydisque relativement
compact de ${\Bbb C}^{n}$ (resp.${\Bbb C}^{p}$ ), $q$, $q'$ et $r$
\'etant les projections canoniques.\smallskip\noindent La {\it  proposition 4} assure l'existence de morphisme trace ${\cal
T}^{n}_{f,\goth{X}}:f_{*}\Omega^{n}_{X/S}\rightarrow
\Omega^{n}_{Y/S}$ ayant les propri\'et\'es que l'on sait. Comme pour
tout faisceau de groupes ab\'eliens ${\cal F}$ sur $X$, on dispose
d'une suite spectrale de Leray $\ds{{\rm I}\!{\rm R}^{i}q_{!}{\rm
I}\!{\rm R}^{j}f_{!}{\cal F}\Longrightarrow{\rm I}\!{\rm
R}^{i+j}{\pi}_{!}{\cal F}}$, qui, au vu de la finitude de $f$,
d\'eg\'en\`ere et donne, en  particulier, l'isomorphisme
 $\ds{{\rm I}\!{\rm R}^{n}{\pi}_{!}\Omega^{n}_{X/S}\simeq {{\rm I}\!{\rm
 R}^{n}q_{!}}f_{*}\Omega^{n}_{X/S}}$, on obtient, en appliquant
 \`a ${\cal T}^{n}_{f}$, le foncteur ${\rm I}\!{\rm
R}^{n}q_{!}$ (qui est exact \`a droite), un morphisme canonique
$$\Phi:={\rm I}\!{\rm R}^{n}q_{!}({\cal T}^{n}_{f,\goth{X}}): {\rm I}\!{\rm
R}^{n}{{\pi}_{!}}\Omega^{n}_{X/S}\rightarrow{\rm I}\!{\rm
R}^{n}q_{!}\Omega^{n}_{Y/S}.$$ Par ailleurs,  pour tout ouvert de
Stein $S'$ de $S$, une simple application d'une formule de type
K\"unneth et du th\'eor\`eme {\bf B} de Cartan,  donne l'
isomorphisme ${\rm I}\!{\rm R}^{n}q_{!}\Omega^{n}_{Y/S}(S')\simeq
{\rm H}^{n}_{c}(U,\Omega^{n}_{U})\hat{\otimes}_{{\Bbb C}} {\cal
O}_{S}(S')$, dont la  compos\'e avec la trace absolue  (ou la
dualit\'e de Serre) ${\rm H}^{n}_{c}(U,\Omega^{n}_{U})\rightarrow
{\Bbb C}$  donn\'ee par l'int\'egration de repr\'esentants de
Dolbeault de classes de cohomologie, donne l'isomorphisme de
faisceaux $\phi: {\rm I}\!{\rm R}^{n}q_{!}\Omega^{n}_{Y/S}\simeq
{\cal O}_{S}$. On en d\'eduit, alors, le diagramme  naturellement
commutatif de faisceaux analytiques
$$\xymatrix{{\rm I}\!{\rm R}^{n}{{\pi}_{!}}\Omega^{n}_{X/S}\ar[rd]_{\int_{\pi,\goth{X}}}\ar[rr]^{\Phi}&& {\rm I}\!{\rm R}^{n}q_{!}\Omega^{n}_{Y/S}\ar[ld]^{\phi}\\
&{\cal O}_{S}&}$$
{\bf Propri\'et\'es de $\int_{\pi,\goth{X}}$}\par
  $\bullet$ {\bf Continuit\'e topologique.} Montrons que ce
morphisme d'int\'egration  est continu dans un sens que
l'on va pr\'eciser. Pour cela, rappelons d'abord que,\par
{\bf(a)}  selon [A.K], [A.B1]), si $X$ est un espace topologique et
$\Phi$ est une famille paracompactifiante de supports de $X$, un
recouvrement ouvert ${\goth U}$ est dit {\bf adapt\'e} \`a $\Phi$ si
pour tout $F\in \Phi$, $\ds{\overline{\bigcup_{{U_{i}\cap
F\not=\emptyset}}U_{i}}}$ appartient \`a $\Phi$. [A.B1] nous
explique bien comment munir les groupes de cohomologie ${\rm
H}^{k}_{\Phi}(X, {\cal F})$ d'une topologie de type limite inductive
d'espaces de Banach et de Souslin.\par {\bf(b)} Si $\pi:X\rightarrow
S$ est un morphisme d'espaces complexes, il est connu que la famille
des ferm\'es de $X$ sur lesquels la restriction de $\pi$ est propre
est  une famille paracompactifiante de supports de $X$ admettant
toujours des recouvrements ouverts adapt\'es.\smallskip\noindent De plus, on sait que  pour tout ouvert de Stein $W$ de $S$, le  groupe
de cohomologie  $\Gamma(W, {\rm I}\!{\rm
R}^{n}{\pi_{!}}{\Omega^{n}_{X/S}})$ (resp. $\Gamma(W, {\cal
O}_{S})$) peut \^etre muni d'une topologie d'espace de type {\bf
Q.F.S} (resp. {\bf F.S}) (cf par exemple  [A.K], [A.B1], [B.S]). Si $\Phi(W)$ est la famille de supports
dont les ferm\'es sont propres sur $W$,  le faisceau analytique
${\rm I}\!{\rm R}^{n}{\pi_{!}}\Omega^{n}_{X/S}$ est engendr\'e, sur
$S$, par le pr\'efaisceau
$$W\rightarrow {\rm H}^{n}_{\Phi(W)}(\pi^{-1}(W), \Omega^{n}_{X/S})$$
et chaque ouvert  $\pi^{-1}(W)$ admet un recouvrement ouvert ${\goth
V}$ dont les \'el\'ements sont $\Phi(W)$-adapt\'es et
``plongeables'' dans des ouverts de type $S\times U\times B$. On
peut remarquer que, gr\^ace \`a  l'ouverture et la surjectivit\'e de
$\pi$,  on peut simplifier la construction en prenant $W:=\pi(V)$ et
nous   ramener  principalement \`a l'\'etude de  la topologie sur
${\rm H}^{n}_{\Phi}(V, \Omega^{n}_{X/W})$. Dans la suite, on note
encore par $X$ cet ouvert $V$.\par\noindent Alors, dans cette
configuration locale, on a les isomorphismes essentiellement d\^us aux d\'eg\'en\'erescence de suites spectrales de Leray
$${\rm H}^{n}_{\Phi}(X, \Omega^{n}_{X/S})\simeq {\rm H}^{n}_{S\times c\times B}(S\times U\times B, \sigma_{*}\Omega^{n}_{X/S})\simeq  {\rm H}^{n}_{S\times c}(S\times U, f_{*}\Omega^{n}_{X/S})$$
On en d\'eduit par des arguments classiques la topologie ad\'equate
(cf [A.K], [A.B1], [Do] etc).\smallskip\noindent Comme on l'a
d\'ej\`a dit dans [KI], {\bf (2.2.1)}, p.56,  la trace est
intimement li\'ee \`a l'int\'egration locale par la formule
$$\int_{|X_{s}|\cap (U'\times B)}\xi\wedge f^{*}\alpha = \int_{\{s\}\times U'}\sum^{j=k}_{j=1}f^{*}_{j}(\xi)\wedge
 \alpha =\int_{\{s\}\times U'}{\cal T}^{q}_{f,\goth{X}}(\xi)\wedge \alpha$$
 o\`u  $s$ est un point fix\'e de $S_{0}$, $U$ un polydisque de ${\Bbb C}^{n}$ d\'efini comme de coutume, $U'$ un polydisque de  $U$ ne rencontrant pas l'ensemble de ramification $R(|X_{s}|)$, $(f_{j})_{1\leq j\leq k}$ d\'esignant les branches locales de $|X_{s}|$ sur $U'$, $\xi\in \Gamma(U\times B,
\Omega^{q}_{U\times B})$ et $\alpha$ une section de $ {\cal
A}^{n-q,n}_{c}(U')$, c'est-\`a-dire une  $(n-q,n)$ forme  \`a
coefficients ${\cal C}^{\infty}$ et \`a support compact dans
$U'$.\par\noindent La convergence des ces int\'egrales \'etant
assur\'ee  par le th\'eor\`eme d'int\'egration sur un sous ensemble
analytique  de Lelong ([L]). Comme $\pi$ est analytiquement
g\'eom\'etriquement plat, les formes m\'eromorphes  $\ds{{\cal
T}_{f}(\xi)= \sum^{j=k}_{j=1}f^{*}_{j}(\xi)}$ se  prolongent
analytiquement  \`a $S\times U$ ; ce qui nous d\'efinit  une
application lin\'eaire \par \centerline{$\ds{\sum^{q=n}_{q=0}{\cal
A}^{n-q,n}_{c}(U)\otimes \Omega^{q}(B)\rightarrow {\cal
O}_{S}(S)}$}\noindent La continuit\'e d'une telle application
r\'esulte soit d' arguments classiques reposant sur les formules
int\'egrales et in\'egalit\'es de Cauchy soit en regardant
simplement la trace d\'efinie dans la situation universelle (cf
[KI], {\bf(1.0.3.2)}, p.23) nous ramenant \`a des majorations de
fonctions ou formes de Newton  (cf {\it lemme 5} de [B1]). Elle se
prolonge donc en une application lin\'eaire continue sur les
compl\'et\'es topologiques $\ds{\Theta:\sum^{q=n}_{q=0}{\cal
A}^{n-q,n}_{c}(U){\widehat{\otimes}}_{{\Bbb C}}
\Omega^{q}(B)\rightarrow {\cal O}_{S}(S)}$. Par ailleurs, la formule
de Stokes nous montre qu'elle est  nulle sur les formes \`a support
$B$-propre et $\bar\partial$-exacte. D'o\`u, gr\^ace \`a une formule
de type K\"unneth ``nucl\'eaire'' et une interpr\'etation connue de
la cohomologie \`a support $B$-propre,
 une (unique) application
lin\'eaire continue
$$\Gamma(S,{\rm I}\!{\rm
R}^{n}{\pi}_{!}(\Omega^{n}_{{S\times U\times B}/S})|_{X})\rightarrow
\Gamma(S, {\cal O}_{S})$$ donnant, de fa\c con \'evidente, gr\^ace
\`a la fl\`eche de restriction (qui est un morphisme nucl\'eaire de
faisceaux  de Frechet)  $\ds({\Omega^{n}_{{S\times U\times
B}/S})|_{X}\rightarrow \Omega^{n}_{X/S}}$, un diagramme commutatif
d'applications continues\par
$$\xymatrix{\Gamma(S,{\rm I}\!{\rm
R}^{n}{\pi}_{!}(\Omega^{n}_{{S\times U\times
B}/S})|_{X}\ar[rr]\ar[rd]&&\Gamma(S,{\rm I}\!{\rm
R}^{n}{\pi}_{!}(\Omega^{n}_{X/S}))\ar[ld]^{\int_{\pi,{\goth X}}}\\
&\Gamma(S, {\cal O}_{S})&}$$ Dans cette situation locale, on dispose
d'une classe fondamentale relative pond\'er\'ee, ${\rm
C}_{\pi,\goth{X}}$,
 pour $\pi$ ([KI], \S {\bf 4},{\it th\'eor\`eme 0}, p.82
 ). On peut supposer, si on le d\'esire, que ${\goth X}$
 est la pond\'eration standard. D\`es lors, on
peut  ``raccrocher''  ce morphisme d'int\'egration $\int_{\pi,{\goth
X}}$  \`a la fl\`eche d'int\'egration le long des fibres de la
projection $p_{1}: S\times U\times B\rightarrow S$.  En effet, en
notant $Z:=U\times B$, il est connu que  la r\'esolution de
Dolbeault $\Omega^{\bullet}_{Z}\rightarrow ({\cal A}^{\bullet,
\star}_{Z},\bar\partial)$ fournit, par tensorisation compl\'et\'ee,
une r\'esolution acyclique  de ${\cal
O}_{S}(V){{\widehat{\otimes}}_{\Bbb C}}\Omega^{\bullet}_{Z}$ par des
faisceaux fins, et ce pour tout ouvert de Stein $V$ de $S$. Il est,
alors,  facile d'en d\'eduire, gr\^ace \`a la trace absolue,  une
trace relative\smallskip \centerline{$\ds{\int_{p_{1}}:{\rm I}\!{\rm
R}^{n+p}{p_{1}}_{!}\Omega^{n+p}_{S\times Z/S}\rightarrow {\cal
O}_{S}}$} \noindent de formation compatible \`a tout changement  de
base. Alors, le cup produit par  la classe fondamentale relative
${\cal C}_{X/S}$ de [B4] ou plus g\'en\'eralement par ${\rm
C}_{\pi,\goth{X}}$, induit le diagramme commutatif\smallskip
\centerline{$\ds{\xymatrix{ {\rm I}\!{\rm R}^{n}\pi_{!}
(\Omega^{n}_{X/S})\ar[rr]^{\wedge{\rm C}_{\pi,\goth{X}}} \ar[rd]_{\int_{\pi,{\goth X}}} &&{\rm I}\!{\rm R}^{n+p}{p_{1}}_{!}\Omega^{n+p}_{S\times Z/S}\ar[ld]^{\int_{p_{1}}}\\
&{\cal O}_{S}&}}$}\par\noindent \par\noindent La compatibilit\'e de
la trace vis-\`a-vis de l'additivit\'e des pond\'erations entraine
\'evidemment la m\^eme propri\'et\'e pour le morphisme
$\int_{\pi,\goth X}$.\smallskip\par $\blacklozenge$
R\'eciproquement, si $\int_{\pi,\goth{X}}:{\rm I}\!{\rm
R}^{n}{\pi}_{!}\Omega^{n}_{X/S}\rightarrow {\cal O}_{S}$ est un
morphisme ${\cal O}_{S}$- lin\'eaire continu et compatible aux
changements de bases r\'eduits, alors $\pi$ est analytiquement
g\'eom\'etriquement plat.\par\noindent En effet, la dualit\'e
analytique relative pour une projection et un morphisme propre de
[R.R.V] et l'isomorphisme de Verdier [V] donne l'isomorphisme
$${\rm I}\! {\rm H}om({\rm I}\!{\rm
R}^{n}\pi_{!}\Omega^{n}_{X/S},  {\cal O}_{S})\simeq {\rm I}\! {\rm
H}om(f_{*}\Omega^{n}_{X/S}, \Omega^{n}_{Y/S})$$ Cela nous dit que
cette int\'egration sur les fibres de $\pi$, accouche naturellement
d' un morphisme $f_{*}\Omega^{n}_{X/S}\rightarrow \Omega^{n}_{Y/S}$
ayant de toute \'evidence les m\^emes propri\'et\'es que
$\int_{\pi}$. On prendra garde que le morphisme ainsi obtenu par
dualit\'e n'est pas automatiquement un morphisme trace!\par\noindent
Mais, par construction m\^eme  de $\int_{\pi,\goth X}$, la fl\`eche
pr\'ec\'edente est un morphisme trace ${\cal T}^{n}_{f}:
f_{*}\Omega^{n}_{X/S}\rightarrow \Omega^{n}_{Y/S}$ prolongeant
naturellement la trace usuelle
$f_{*}f^{*}\Omega^{n}_{Y/S}\rightarrow \Omega^{n}_{Y/S}$ ,
compatible aux changements de bases r\'eduits et aux localisations
sur $X$ et $Y$. Ceci \'etant vrai pour tout $f$ fini et surjectif,
on conclut gr\^ace a la {\it proposition 4}.
$\,\,\blacksquare$\bigskip\noindent Arm\'e de la {\it proposition 2} traitant la situation locale que l'on va globaliser,  nous
sommes en mesure \`a pr\'esent de donner la \smallskip\noindent {\bf
2.0.7. Preuve de la proposition fondamentale
{\bf 1}.}\smallskip\noindent Dans toute la suite on supposera
$X$  paracompact et compl\`etement paracompact.\par $\blacklozenge$
Supposons $\pi$ analytiquement g\'eom\'etriquement plat et donnons
nous un recouvrement ouvert localement fini $(X_{\alpha})_{\alpha\in
A}$ de $X$  muni d'installations locales
 $$\xymatrix{X_{\alpha}\ar[rdd]_{\pi_{\alpha}}\ar[rr]^{\sigma_{\alpha}}\ar[rd]^{f_{\alpha}}&&Z_{\alpha}\ar[ld]_{p_{\alpha}}\ar[ldd]^{r_{\alpha}}\\
&Y_{\alpha}\ar[d]^{q_{\alpha}}&\\&S&}$$ avec $\sigma_{\alpha}$  un
plongement local, $\pi_{\alpha}$ la restriction de $\pi$ \`a
$X_{\alpha}$, $Y_{\alpha}$ et $Z_{\alpha}$ sont lisses sur $S$,
$p_{\alpha}$,$q_{\alpha}$ et $r_{\alpha}$ lisses. De tels
recouvrements existent puisque $\pi$ est universellement
$n$-\'equidimensionnel. De plus, on peut les choisir convenablement
adapt\'es \`a la famille des supports $S$-propres dans $X$ (cf [A.K],
[A.B1]\footnote{$^{(6)}$}{La {\it proposition (1.1)} de [A.B1] montre d'ailleurs que, si $X$ est paracompact et compl\`etement paracompact, toute famille paracompactifiante \`a base d\'enombrable admet toujours un recouvrement ouvert localement fini adapt\'e.} Alors, $\pi$ \'etant analytiquement g\'eom\'etriquement
plat et cette notion \'etant compatible aux localisations ouvertes
sur $X$ et $S$, chaque morphisme  $\pi_{\alpha}$ l'est aussi. Il
r\'esulte, alors, de la {\it proposition 2}, une collection de
morphismes d'int\'egration
$$\int_{\pi_{\alpha}, \goth{X}_{\alpha}}:{\rm I}\!{\rm R}^{n}{\pi_{\alpha}}_{!}\Omega^{n}_{{X}_{\alpha}/S}\rightarrow {\cal O}_{S}$$
$\goth{X}_{\alpha}$ \'etant la restriction \`a l'ouvert $X_{\alpha}$
de la pond\'eration analytiquement g\'eom\'etriquement plate $\goth
X$.\par\noindent Mais $\pi$ \'etant \`a fibres de dimension pure $n$,
 le lemme de Reiffen assure l'exactitude \`a droite du foncteur ${\rm
I}\!{\rm R}^{n}\pi_{!}$ de laquelle r\'esulte la suite exacte \`a
droite
$$ \bigoplus_{\alpha\in A}{\rm I}\!{\rm R}^{n}{\pi_{\alpha}}_{!}
\Omega^{n}_{{X_{\alpha}}/S} \rightarrow {\rm I}\!{\rm
R}^{n}{\pi_{!}}{\Omega^{n}_{X/S}}\rightarrow 0$$
\par\noindent Ainsi, obtient-on, par recollement, une fl\`eche globale
$$\int_{\pi,\goth{X}}:{\rm I}\!{\rm R}^{n}{\pi}_{!}\Omega^{n}_{X/S}\rightarrow {\cal O}_{S}$$
dont il est facile de v\'erifier les propri\'et\'es annonc\'ees
puisqu'elles sont h\'erit\'ees des propri\'et\'es de
$\int_{\pi_{\alpha}, \goth{X_{\alpha}}}$ donn\'ees par la {\it proposition 2}.\par\noindent
 L'ind\'ependance de la
construction vis-\`a-vis  du recouvrement ouvert choisi d\'ecoule de
la notion m\^eme de platitude g\'eom\'etrique contenue dans la {\it proposition 4} ou du lemme des recouvrements de Varouchas [Va] en
prenant en compte la conservation de la $S$-analyticit\'e par
r\'etr\'ecissement en $U_{\alpha}$ et $B_{\alpha}$.\par\noindent Sa
compatibilit\'e globale  avec tout changement de base entre espaces
complexes r\'eduits s'impose par d\'efinition. En effet, si
$\eta:T\rightarrow S$ est  un morphisme d'espaces complexes
r\'eduits et
$$\xymatrix{X_{T}\ar[r]^{\Theta}\ar[d]_{\pi_{T}}&X\ar[d]^{\pi}\\
T\ar[r]_{\eta}&S}$$ le diagramme de changement de base associ\'e.
Alors $\pi_{T}$ est analytiquement  g\'eom\'etriquement plat et
 la compatibilit\'e aux changements de bases r\'esulte de la
 formule (classique) $$\int_{X_{t}}\Theta^{*}\xi :=
 \eta^{*}\int_{X_{\eta(t)}}\xi$$
$\blacklozenge$ R\'eciproquement, supposons donn\'e un morphisme
$\pi:X\rightarrow S$ universellement $n$-\'equidimensionnel muni
d'une certaine pond\'eration $\goth{X}$ et d'un morphisme
d'int\'egration global
$$\int_{\pi,\goth{X}}: {\rm I}\!{\rm R}^{n}{\pi_{!}}\Omega^{n}_{X/S}\rightarrow {\cal O}_{S}$$
${\cal O}_{S}$-lin\'eaire, continu, compatible \`a tout changement
de base entre espaces complexes r\'eduits de dimension finie et
donnant l'int\'egration usuelle dans le cas absolu. Montrons, alors,
que $\pi$ est analytiquement  g\'eom\'etriquement plat.\par\noindent
On suppoera $X$ d\'enombrable \`a l'infini, paracompact et
compl\`etement paracompact.\par\noindent
 Comme le faisceau
analytique ${\rm I}\!{\rm R}^{n}{\pi_{!}}\Omega^{n}_{X/S}$ est
engendr\'e, sur $S$, par le pr\'efaisceau
$$W\rightarrow {\rm H}^{n}_{\Phi(W)}(\pi^{-1}(W), \Omega^{n}_{X/S})$$
$\Phi(W)$ \'etant une famille de supports dont les ferm\'es sont
propres sur $W$ ouvert de Stein de $S$, il est facile d'adapter
convenablement  [B.V] pour d\'ecrire explicitement le morphisme
d'int\'egration au moyen d'un calcul en termes  de cochaines  de  \v
Cech (cf [KI], {\bf (2.1)}, p.53). La seule difficult\'e technique consiste
\`a produire des recouvrements ouverts trivialisants pour $\pi$, en
prenant des cartes adapt\'ees \`a chaque fibre et \`a la famille des
ferm\'es $S$-propres dans $X$.
\par\noindent  Comme $\int_{\pi,\goth{X}}$  est global et
compatible aux changement de bases, il ne d\'epend pas du
recouvrement ouvert choisi pour le d\'efinir et, est compatible en
un sens \'evident aux raffinements entre recouvrements ouverts sur
$X$ qui est paracompact. On peut donc se ramener \`a la situation
locale donn\'ee par la {\it proposition 2}. Les propri\'et\'es de
ce morphisme ont pour effet imm\'ediat de rendre invariante, par
raffinement ou r\'etr\'ecissement des ouverts,  la $S$-analyticit\'e
des traces obtenues par d\'ecoupages dans des donn\'ees locales. En
effet, $\pi$ \'etant universellement $n$- \'equidimensionnel, on
peut toujours le factoriser localement en un diagramme commutatif
$\ds{\xymatrix{X\ar@/_/[rr]_{\pi}\ar[r]^{f}&Y\ar[r]^{q}&S}}$
  avec $f$ fini et surjectif sur $Y$ lisse sur $S$ et de
dimension relative $n$. Mais comme on l'a d\'ej\`a vu, dans  cette
situation locale, on a  un morphisme d'int\'egration locale
$\ds{{\rm I}\!{\rm R}^{n}{{q}_{!}}f_{*}\Omega^{n}_{X/S}\rightarrow
{\cal O}_{S}}$, induisant, par construction m\^eme de
l'int\'egration sur les cycles, un morphisme de faisceaux ${\cal
O}_{X}$-coh\'erents $\ds{f_{*}\Omega^{n}_{X/S}\rightarrow
\Omega^{n}_{Y/S}}$ qui est n\'ecessairement un morphisme
{\it{trace}} de formation naturellement  stable par n'importe quel
changement de base entre espaces complexes r\'eduits et rendant
commutatif le diagramme
$\ds{\xymatrix{f_{*}f^{*}\Omega^{n}_{Y/S}\ar[r]\ar@/_/[rr]&f_{*}\Omega^{n}_{X/S}\ar[r]&\Omega^{n}_{Y/S}}}$
 qui impose la platitude
g\'eom\'etrique de $f$.\par\noindent De plus, l'invariance par
r\'etr\'ecissement ou raffinement des ouverts  dont il a \'et\'e
question plus
 haut \`a propos du morphisme d'int\'egration globale montre que ces traces locales h\'eritent de cette propri\'et\'e.\par\noindent
 Il apparait, ainsi, que le morphisme global d'int\'egration induit  une collection de
morphismes {\it{trace}} compatibles aux  changement de base entre
espaces complexes r\'eduits, prolongeant  naturellement les traces
usuelles, compatibles aux localisations sur $X$  (et par suite
insensible aux changement de projections quelconques); ce qui prouve
la platitude g\'eom\'etrique analytique  de $\pi$ en vertu de la
{\it proposition 2}. $\,\,\blacksquare$\rm\bigskip\noindent
{\bf 2.0.8. Quelques corollaires.}\smallskip\noindent
 \Cor{2}{}{\it{Equidimensionnalit\'e  universelle et
trace.}}\par\noindent
 Soit $\pi\in {\cal E}(S,n)$ pond\'er\'e par un
cycle $\goth{X}$. Alors, pour toute projection finie $f:X\rightarrow
S\times U$ o\`u $U$, comme de coutume, d\'esigne un ouvert
relativement compact de ${\Bbb C}^{n}$, il existe un unique
  morphisme trace  ${\cal O}_{S}$-lin\'eaire continue\par
  \centerline{${\goth
T}^{n}_{f,\goth{X}}:f_{*}\Omega^{n}_{X/S}\rightarrow
{\cal L}^{0}_{S}{\hat\otimes}_{\Bbb C}\Omega^{n}_{U}$}
\noindent compatible aux changements de bases et \`a
 l'additivit\'e des pond\'erations. De plus, si la pond\'eration
 $\goth{X}$ est contin\^ument g\'eom\'etriquement plate, ${\Im}({\goth
T}^{n}_{f,\goth{X}})\subset {\cal O}^{c}_{S}{\hat\otimes}_{\Bbb
C}\Omega^{n}_{U}$. \rm \dem $\pi:X\rightarrow S$ \'etant
universellement $n$-\'equidimensionnel, il se factorise localement
en un morphisme fini suivi d'une projection ou plus pr\'ecisemment
s'installe localement, par rapport \`a l'une quelconque de ses
fibres, dans  le diagramme commutatif
$$\xymatrix{X\ar[r]^{\sigma}\ar[d]_{\pi}\ar[rd]^{f}&S\times U\times B\ar[d]\\
S&S\times U\ar[l]}$$, dans lesquelles $f$ est un morphisme fini
surjectif et
 ouvert, $\sigma$ un plongement local et o\`u
 les fl\`eches non pr\'ecis\'ees sont les projections canoniques.\par\noindent
Alors, l'image directe d\'efinie au sens des courants (cf
[KI],{\bf(1.0.3.1)}, p.22) donn\'ee par le morphisme\par\noindent
$\ds{{f}_{*}: {f}_{*}{\cal D}^{n,0}_{X}\rightarrow {\cal
D}^{n,0}_{S_{0}\times U}}$ transforme, en particulier, toute
$n$-forme holomorphe sur $X$ en un courant $\bar\partial$-ferm\'e de
type $(n,0)$ sur  $S_{0}\times U$. Mais, modulo ${\cal
O}_{S_{0}\times U}$- torsion, tout courant $\bar\partial$-ferm\'e de
type $(n,0)$ s'identifie naturellement, d'apr\`es [B2],  \`a une
section du faisceau $\omega^{n}_{S_{0}\times U}$ qui, d'apr\`es
[K2], se `` d\'ecompose nucl\'eairement `` en
$\ds{\bigoplus_{i+j=n}\omega^{i}_{S_{0}}{\hat{\otimes}_{\Bbb
C}}\Omega^{j}_{U}}$, gr\^ace au lemme de Dolbeault-Grothendieck et
aux r\'esulats g\'en\'eraux de ce dernier sur les espaces
nucl\'eaires. Il apparait clairement que ce morphisme pr\'ecedemment
cit\'e agit sur les formes $S$-relatives en donnant la fl\`eche
$\ds{{f_{*}}:f_{*}\Omega^{n}_{X/S_{0}}\rightarrow
\omega^{0}_{S_{0}}{\hat{\otimes}_{{\Bbb
C}}}\Omega^{n}_{U}}$.\par\noindent Par ailleurs, on dispose, d'un
op\'erateur {\it trace} associant \`a toute $n$- forme  holomorphe
$\xi$ sur $X$ l'expression $\ds{{\cal
T}^{n}_{f}(\xi):=\sum^{l=k}_{l=1}f^{*}_{l}(\xi)}$, dans laquelle
$(f_{l})_{1\leq l\leq k}$ d\'esigne  les branches locales du
rev\^etement ramifi\'e g\'en\'erique d\'efini par $f$ (cf [KI],
{\bf( 1.0.3.1)}, p.23). Il va sans dire que
 les op\'erateurs  ${\cal T}^{n}_{f}$ et ${f}_{*}$ agissent
 de la m\^eme fa\c con sur les formes $S$-relatives. \par\noindent
On peut pr\'eciser un peu plus l'image de cette trace en utilisant
la stabilit\'e par changement de base de la notion d'universalit\'e
\'equidimensionnelle. En effet, en consid\'erant une
d\'esingularisation de $S$ ou tout simplement sa normalisation
$\eta:\hat{S}\rightarrow S$, on voit,  gr\^ace a la {\it proposition
4},  que le morphisme $\hat{\pi}:\hat{X}\rightarrow {\hat{S}}$
d\'eduit de $\pi$ dans un tel changement de base est encore
$n$-universellement \'equidimensionnel sur une base normale (ou
lisse). Les factorisations locales produisent, alors,  des traces
locales $\ds{{\hat{\cal
T}}^{n}_{\hat{f}}:{\hat{f}}_{*}\Omega^{n}_{\hat{X}/\hat{S}}\rightarrow
\Omega^{n}_{\hat{S}\times U/\hat{S}}}$ puisque sur un espace normal
toute fonction g\'en\'eriquement holomorphe est automatiquement
holomorphe. Par cons\'equent,  $\hat{\pi}$ est analytiquement
g\'eom\'etriquement plat en vertu de la {\it proposition 4} . Alors,
utilisant des arguments standards de transversalit\'e d'
applications nucl\'eaires, on en d\'eduit  que l'image des  traces
locales relatives \`a $\pi$ sont \`a  valeurs dans le faisceau
$\ds{\eta_{*}{\cal O}_{\hat{S}}{\hat\otimes}_{\Bbb
C}\Omega^{n}_{U}}$.
 Mais $\eta_{*}{\cal O}_{\hat{S}}= {\cal O}^{b}_{S}={\cal L}^{0}_{S_{0}}$\footnote{$^{(7)}$}{ Si $Z$ est un espace complexe r\'eduit et de dimension pure, on note ${\cal L}^{k}_{Z}$ le faisceau des $k$-formes m\'eromorphes sur $Z$ dont l'image r\'eciproque par tout morphisme de d\'esingularisation se prolonge holomorphiquement sur la d\'esingularis\'ee de $Z$. Pour plus de pr\'ecisions, on renvoie le lecteur \`a   [K2].}  est  le faisceau des fonctions
m\'eromorphes se prolongeant holomorphiquement sur toute
d\'esingularis\'ee de $S_{0}$ ou  simplement le faisceau des
fonctions m\'eromorphes localement born\'ees que l'on avait not\'e
${\cal O}^{b}_{S}$ dans [KI], {\bf(1.2.0)}, p.32.\par \noindent La
compatiblit\'e aux changements de bases est une cons\'equence
conjointe de la stabilit\'e du faisceau   ${\cal L}^{0}_{S_{0}}$ par
image r\'eciproque de morphisme d'espaces complexes et des
propri\'et\'es de la trace donn\'es par la {\it proposition 4}. La
compatibilit\'e \`a l'additivit\'e des pond\'erations est
\'evidente.\par\noindent Si la pond\'eration est contin\^ument
g\'eom\'etriquement plate, le
 th\'eor\`eme du changement de projection (cf [B1], {\it
 th\'eor\`eme 2}, p.42) disant que le passage d'une carte \`a une
 autre est toujours m\'eromorphe continue, est strictement
 \'equivalent au fait que, pour toute param\'erisation locale
 $S$-relative $f:X\rightarrow S$, le morphisme trace est en fait \`a
 valeurs dans le faisceau ${\cal O}^{c}_{S}{\widehat\otimes}_{\Bbb
 C}\Omega^{n}_{U}$ (i.e la variation en $s$ est m\'eromorphe
 continue).$\blacksquare$\smallskip\noindent \Cor {3} Soit $\pi\in {\cal
E}(S,n)$ muni d'une pond\'eration $\goth{X}$. Alors $\pi$ appartient
\`a ${\cal G}_{a}(S,n)$ si et seulement si pour toute factorisation
locale $\ds{\xymatrix{X\ar@/_/[rr]_{\pi}\ar[r]^{f}&Y\ar[r]^{q}&S}}$
de $\pi$ avec $f$ fini et surjectif et $q$ lisse de dimension
relative $n$, il existe un (unique) morphisme trace ${\cal
T}^{0}_{f,\goth{X}}:f_{*}{\cal O}_{X}\rightarrow {\cal O}_{Y}$
compatible aux changements de bases entre espaces complexes
r\'eduits et de nature local sur $X$ et $Y$.\rm \dem C'est un
exercice simple d'alg\`ebre lin\'eaire. Pour \^etre plus pr\'ecis,
supposons  $Y=S\times U$ et $X$ plong\'e dans $S\times U\times B$ et
notons $t:=(t_{1},\cdots,t_{n})$ (resp  $x:=(x_{1},\cdots,x_{p})$)
le point courant de $U$ (resp. $B$) et d\'esignons par
$q_{0}:U\times B\rightarrow U$ la projection lin\'eaire canonique
envoyant $(t,x)$ sur $t$. En vertu du {\it lemme (5.3)} de [KI], \S{\bf 5}, p.86,  toute $n$-forme holomorphe $S$ relative
sur $X$ au voisinage d'un point fix\'e s'\'ecrit  comme la somme
$\ds{\sum_{|I|+|J|=n}g_{I,J}(s,t,x)dt^{I}\wedge dx^{J}}$. Mais
chacun de ces termes s'ecrira  sous la forme
$\ds{g_{I,J}(s,t,x)H_{I,J}(s,t,x)dt}$ avec $H_{I,J}(s,t,x)$
s'exprimant naturellement en fonction des mineurs de la jacobienne
du changement de projection associ\'e. Alors, la lin\'earit\'e de la
trace permet, gr\^ace \`a l'hypoth\`ese, d'en d\'eduire un morphisme
trace $\ds{{\cal
T}^{n}_{f,\goth{X}}:f_{*}\Omega^{n}_{X/S}\rightarrow
\Omega^{n}_{S\times U/S}}$ et ce pour tout $f$.D'o\`u la platitude
g\'eom\'etrique analytique en vertu de la {\it proposition 4}
$\,\,\blacksquare$  \smallskip\smallskip\noindent \Cor{4}{} Soit
$\pi\in {\cal E}(S,n)$ pond\'er\'e par un cycle $\goth{X}$. Alors
\par\noindent {\bf(i)} si $S$ est localement irr\'eductible, $\pi$
est contin\^ument g\'eom\'etriquement plat.\par\noindent {\bf(ii)}
si $S$ est normal $\pi$ est analytiquement g\'eom\'etriquement
plat.\par\noindent {\bf(iii)}si $S$ est faiblement normal et  $\pi$
contin\^ument g\'eom\'etriquement plat, $\pi$ est analytiquement
g\'eom\'etriquement plat.\rm
\dem Elle repose essentiellement sur le {\it corollaire 2}. Soit $s_{0}$ dans $S$. Alors,  il
existe voisinage ouvert $S_{0}$ de $s_{0}$ dans $S$ et une famille
d'\'ecailles ou de cartes  $(E_{\alpha}:=(U_{\alpha}, B_{\alpha},
\sigma_{\alpha}))_{\alpha\in A}$ $S_{0}$- adapt\'ees  \`a
$X_{s_{0}}$ de sorte que, pour tout $\alpha\in A$,  on ait une
installation (locale) ($\clubsuit$)
$$\xymatrix{X_{\alpha}\ar[d]^{\pi_{\alpha}} \ar [dr]_{f_{\alpha} }\ar [r]^{\sigma_{\alpha}}
 & S_{0}\times U_{\alpha}\times B_{\alpha} \ar[d]^{p_{\alpha}}\\
 S_{0}&\ar[l]^{q_{\alpha}} S_{0}\times U_{\alpha}}$$
dans laquelle $f_{\alpha}$ est un morphisme fini et surjectif,
$\sigma_{\alpha}$ est un plongement local, $p_{\alpha}$ et
$q_{\alpha}$ \'etant les projections canoniques. De telles cartes, que l'on peut choisir centr\'ees par un rapport \`a un point privil\'egi\'e sur la fibre, existent toujours (cf [A.S]).\par\noindent
D'apr\`es le {\it corollaire 2}, pour chaque indice $\alpha$, on dispose d'un morphisme trace\par
\centerline{$ {\goth
T}^{n}_{f_{\alpha}}:{f_{\alpha}}_{*}\Omega^{n}_{X_{\alpha}/S_{0}}\rightarrow
{\cal L}^{0}_{S_{0}}{\hat\otimes}_{\Bbb C}\Omega^{n}_{U_{\alpha}}$}
\noindent
Il nous suffit simplement d'utiliser la {\it proposition 4} et  la caract\'erisation de l'irr\'eductibilit\'e locale, de  la
normalit\'e et de la  faible normalit\'e  \'enonc\'ee dans [KI], {\bf(1.2.0)}, p.32,  disant que \par\noindent
$\bullet$ ${\cal O}^{b}_{S}\simeq {\cal O}^{c}_{S}$ si et seulement si $S$ est localement irr\'eductible,\par\noindent
$\bullet$ ${\cal O}_{S}\simeq {\cal O}^{c}_{S}$ si et seulement si $S$ est faiblement normal,\par\noindent
$\bullet$ ${\cal O}_{S}\simeq {\cal O}^{b}_{S}$ si et seulement si $S$ est normal.$\,\,\,\blacksquare$\smallskip\par\noindent
\Cor{5}{} Dans la cat\'egorie des espaces analytiques complexes
r\'eduits  de dimension finie, les morphismes  de {\it Tor-dimension
finie} \`a fibres de dimension pure constante  sont analytiquement g\'eom\'etriquement plats.\rm \dem
Comme les fibres sont de dimension pure constante $n$, le lemme de factorisation ([Fi]) montre que pour  tout $x$ de $X$ il existe une factorisation locale en $x$ c'est-\`a-dire un diagramme coomutatif  $\ds{\xymatrix{X\ar@/_/[rr]_{\pi}\ar[r]^{f}&Y\ar[r]^{q}&S}}$ o\`u $X$ d\'esigne en fait un  voisinage  ouvert $U$ de $x$, $Y$ lisse sur $S$ que l'on peut \'ecrire sous la forme $S\times V$ avec $V$ ouvert  de ${\Bbb C}^{n}$,  $f$ fini et surjectif et $q$ lisse de dimension relative
$n$.\par\noindent
 Il est facile de voir que $f$ est n\'ecessairement de {\it Tor-dimension} finie. La preuve donn\'ee dans [M.F.K] et que nous rappelons pour la commodit\'e du lecteur  s'applique bien dans le contexte de la g\'eom\'etrie analytique locale. En effet, consid\'erons le diagramme cart\'esien
$$\xymatrix{X\times_{S}Y\ar[r]^{q'}\ar[d]_{\pi'}&X\ar@/^/[l]^{\sigma}\ar[d]^{\pi}\\
Y\ar[r]_{q}&S}$$
dans lequel $\sigma:=(id_{X}\times f)$  est une section de $q'$.
Comme $q$ est lisse, $q'$ l'est aussi. Par ailleurs, le changement de
 base donn\'e par $q$ \'etant cohomologiquement transverse, $\pi'$ est aussi de {\it Tor-dimension finie} (en fait, on peut se ramener au cas o\`u $Y=S\times U$ permettant de voir rapidement que $\pi'$ poss\`ede cette propri\'et\'e). Alors $X$ s'identifie \`a un sous espace ferm\'e  localement intersection compl\`ete dans $ X\times_{S}Y$ et donc de {\it Tor-dimension finie}. Mais
 $f$ \'etant la compos\'ee de deux morphismes de  {\it Tor-dimension finie}
  l'est aussi\footnote{$^{(8)}$}{R\'eciproquement, si un morphisme se factorise localement en un morphisme fini de {\it Tor-dimension finie} suivi d'un morphisme lisse, il est n\'ecessairement de {\it Tor-dimension finie}. Par ailleurs, d'apr\`es [Ry], un morphisme fini de {\it Tor-dimension finie} d\'efinit toujours une famille analytique de $0$-cycles et il est donc ouvert.  } \par\noindent
Il nous reste alors  \`a voir que chaque morphisme fini d'une telle
factorisation est naturellement muni d'un morphisme trace. Or, pour
cela, il suffit simplement de se rappeler que  la construction de
telles traces \'etant de nature locale, on se ram\`ene, apr\`es
``germification,  \`a un probl\`eme  d'alg\`ebre  analytique  locale
que l'on surmonte en adaptant sans efforts cons\'equents  le lemme
de construction de la trace pour les morphismes de {\it
Tor-dimension finie} (alg\'ebriques) (cf [KI], {\it lemme (5.0)},
p.85). On produit, ainsi,  des morphismes traces au sens
pr\'ec\'edemment d\'ecrits  et,  on  termine  en  invoquant la {\it
proposition 4}.$\,\,\blacksquare$ \bigskip\noindent {\bf 2.0.9.
Remarques.}\par\noindent {\bf(i)} les hypoth\`eses du {\it
corollaire 3} nous donnent, si $X$ est plong\'e dans $Z$ lisse et de
dimension relative $p$ sur $Y:=S\times U$, un morphisme ${\cal
O}_{X}\rightarrow {\cal E}xt^{p}({\cal O}_{X}, \Omega^{p}_{Z/Y})$.
Il est facile de voir que la platitude g\'eom\'etrique analytique
revient \`a montrer l'existence d'un unique morphisme ${\cal
O}_{X}\rightarrow {\cal E}xt^{p}({\cal O}_{X}, \Omega^{p}_{Z/S})$
envoyant la fonction $1$ sur une classe particuli\`erement
interessante que l'on aura l'occasion d'\'evoquer dans
{\bf(2.1)}.\par\noindent {\bf(ii)} Le {\it corollaire 4} donne des
conditions suffisantes (et non n\'ecessaires) portant sur l'espace
des param\`etres. Par construction, les morphismes  contin\^ument
  g\'eom\'etriquement plat (resp. analytiquement g\'eom\'etriquement plat)
  produisent toujours des fonctions m\'eromorphes continues sur $S$
   c'est-\`a-dire des sections du faisceau ${\cal O}^{c}_{S}$
   (resp. holomorphes sur $S$). \par\noindent
{\bf(iii)}  Comme on le voit, la nature du morphisme est
enti\`erement
 d\'etermin\'ee par la possibilit\'e ou pas de prolonger
 analytiquement les sections de ce faisceau. Ainsi,  le coeur du
probl\`eme  consiste \`a trouver une raison assurant ce prolongement
holomorphe sur $S_{0}\times U_{\alpha}$ tout entier, en ayant la
stabilit\'e par changement de base quelconque entre espaces
complexes r\'eduits et la compatibilit\'e aux restrictions ouvertes
sur $U_{\alpha}$\smallskip\bigskip\noindent {\bf 2.1. Platitude
g\'eom\'etrique analytique et morphisme classe fondamentale
 relative pond\'er\'e.}\bigskip\noindent
{\bf 2.1.0.}  Avant d'entrer dans le vif du sujet, rappelons rapidement
 que:\par\noindent
$\bullet$  si $\pi:X\rightarrow S$ est un morphisme universellement $n$-\'equidimensionnel
 pond\'er\'e par un certain cycle $\goth{X}$, un
faisceau coh\'erent ${\cal F}$ sur $X$ est dit v\'erifi\'e {\bf la
propri\'et\'e de la trace relative} si  pour tout diagramme
commutatif d'installation locale ($\spadesuit$)
$$\xymatrix{X\ar[rdd]_{\pi}\ar[rr]^{\sigma}\ar[rd]^{f}&&Z\ar[ld]_{q'}\ar[ldd]^{q}\\
&Y\ar[d]^{q''}&\\
&S&}$$ dans lequel  $\sigma$ est un  plongement local de codimension
$p$ dans $Z$ lisse sur $S$ de dimension relative $m:=n+p$, $Y$ lisse
sur $S$ de dimension relative $n$, $q'$, $q''$ et $q$ sont des
morphismes  lisses, $f$  fini, surjectif et o\`u l'on d\'esigne
encore par $\pi$ la restriction de $\pi$ \`a  une telle donn\'ee
locale,  est naturellement associ\'e \`a  un diagramme commutatif
$$\xymatrix{f_{*}f^{*}\Omega^{n}_{Y/S}[n]
\simeq f_{*}({\cal O}_{X}\otimes
f^{*}\Omega^{n}_{Y/S}[n])\ar[d]_{{\cal T}^{0}_{f,\goth{X}}\otimes
Id[n]} \ar[r]&\!\!\!f_{*}\Omega^{n}_{X/S}[n]
\ar[r]^{\indent{f_{*}{\cal C}_{\pi}}}&f_{*}{\cal F}
\ar[d]^{{\goth T}^{n}_{f,\goth{X}}}\\
\Omega^{n}_{Y/S}[n]\ar[rr]_{Id}&&\Omega^{n}_{Y/S}[n]}$$\par\noindent
dans lequel ${\cal T}^{0}_{f,\goth{X}}: f_{*}{\cal O}_{X}\rightarrow
{\cal O}_{Y}$ est la trace usuelle et ${\goth T}^{n}_{f,\goth{X}}$
est un morphisme {\it{trace}} d\'etermin\'e par
${\goth{X}}$.\smallskip\noindent
$\bullet$  Si $\pi:X\rightarrow S$ est un
morphisme universellement $n$-\'equidimensionnel et
$(X_{\alpha})_{\alpha\in A}$ un recouvrement ouvert localement fini
de $X$ adapt\'e \`a $\pi$,  on d\'esignera par
$(\spadesuit)_{\alpha}$ tout diagramme d'installation locale
donn\'ee par
 $$\xymatrix{X_{\alpha}\ar[rdd]_{\pi_{\alpha}}\ar[rr]^{\sigma_{\alpha}}\ar[rd]^{f_{\alpha}}&&Z_{\alpha}\ar[ld]_{q'_{\alpha}}\ar[ldd]^{q_{\alpha}}\\
&Y_{\alpha}\ar[d]^{q''_{\alpha}}&\\&S&}$$ o\`u $\sigma_{\alpha}$ est
un plongement local, $\pi_{\alpha}$ la restriction de $\pi$ \`a
$X_{\alpha}$, $Y_{\alpha}$ et $Z_{\alpha}$ des espaces complexes
lisses sur $S$ de dimension relative $n+p$ et $n$ respectivement,
$q_{\alpha}$,$q'_{\alpha}$ et $q''_{\alpha}$ lisses. De tels
recouvrements existent puisque $\pi$ est universellement
$n$-\'equidimensionnel. De plus, on peut les choisir convenablement
adapt\'es \`a la famille des supports $S$-propres dans $X$(cf [A.K],
[A.B1]). Si $\pi$ est analytiquement g\'eom\'etriquement plat, il en
est de m\^eme de $\pi_{\alpha}$ et $f_{\alpha}$ puisque cette notion
 est  compatible aux localisations ouvertes sur $X$ et
$S$.\smallskip\noindent Enfin, rappelons (cf [KI], {\bf(4.1.1)}, p.82)
 \Th{0}{} Soient $S$, $Z$ et $X$ des espaces
analytiques complexes r\'eduits install\'es dans le diagramme
commutatif
$\xymatrix{X\ar[r]^{\sigma}\ar@/_/[rr]_{\pi}&Z\ar[r]^{q}&S}$ dans
lequel $q$ est lisse de dimension relative $n+p$, $\sigma$ est un
plongement local et $\pi$ un morphisme contin\^ument
$n$-g\'eom\'etriquement plat (donc  universellement
$n$-\'equidimensionnel)  muni d'une certaine pond\'eration
$\goth{X}$.
\par\noindent Alors $\pi$ est  analytiquement
$n$-g\'eom\'etriquement plat si et seulement si il existe  une
classe de cohomologie  ${\rm C}_{\pi,\goth{X}}$ de ${\rm
H}^{p}_{X}(Z, \Omega^{p}_{Z/S})$ v\'erifiant les propri\'et\'es
suivantes:\par\noindent elle est de nature locale sur $X$ et $S$, de
formation compatible \`a l'additivit\'e des pond\'erations et \`a
tout changement de base entre espaces complexes r\'eduits donnant,
en particulier, pour chaque $s$ fix\'e la classe fondamentale
absolue du cycle $[\pi^{-1}(s)]$. De plus, si $\goth X$ est la
pond\'eration standard, ${\rm C}_{\pi,\goth{X}}={\cal C}_{X/S}$
construite dans [B4]. \rm\smallskip\noindent L'\'el\'ement ${\rm
C}_{\pi,\goth{X}}$ a \'et\'e appel\'e {\it{classe fondamentale
relative pond\'er\'ee}} associ\'ee \`a $\pi$. \smallskip\noindent
Dans ce qui suit, nous allons \'etablir les deux propositions
suivantes :
 \Prop{5}{} {\it Cas global}.\par\noindent
 Soit  $\pi$ un
\'el\'ement de ${\cal E}(S,n)$ muni d'une pond\'eration ${\goth X}$.
Alors, $\pi$ est analytiquement g\'eom\'etriquement plat si et
seulement si il existe un unique morphisme canonique de faisceaux
coh\'erents
 ${\cal C}_{\pi,\goth{X}}:\Omega^{n}_{X/S}\rightarrow \omega^{n}_{\pi}$
 v\'erifiant les propri\'et\'es suivantes:\par\noindent {\bf(i)} il est de
formation compatible \`a l'additivit\'e des pond\'erations\footnote{$^{(9)}$}{En un sens que l'on pr\'ecisera.}, aux
 changement de base r\'eduit et satisfait
  {\it la propri\'et\'e de la trace relative},\par\noindent
{\bf(ii)} il induit  un morphisme de complexes  diff\'erentiels
gradu\'es\par
 \centerline{${\cal C}^{\bullet}_{\pi, \goth{X}}:\Omega^{\bullet}_{X/S}\rightarrow \omega^{\bullet}_{\pi}$}
\rm\smallskip \Prop{6}{}{\it Cas local}.\par\noindent
  Soient $X$, $Z$ et $S$ trois espaces analytiques
complexes r\'eduits munis du diagramme commutatif $\xymatrix{X\ar@/_/[rr]_{\pi}\ar[r]^{\sigma}&Z\ar[r]^{q}&S}$ dans lequel $\pi$ est un morphisme universellement $n$-\'equidimensionnel pond\'er\'e par un certain cycle  ${\goth X}$, $\sigma$ est un plongement dans $Z$ lisse sur $S$ et
 de dimension relative $n+p$ et $q$ la  projection canonique.\smallskip\noindent
 Alors $\pi$ est analytiquement
g\'eom\'etriquement plat  si et seulement si il existe un unique
morphisme canonique de faisceaux coh\'erents ${\cal C}_{\pi,\goth{X}}:\Omega^{n}_{X/S}\rightarrow
{\cal E}xt^{p}({\cal O}_{X}, \Omega^{n+p}_{Z/S})$ de formation compatible \`a l'additivit\'e des pond\'erations et \`a
tout changement de base entre espaces complexes r\'eduits de
dimension finie et satisfaisant  {\it la propri\'et\'e de la trace
relative}.\rm\smallskip\noindent Il nous est, alors,  facile de connecter ce morphisme aux
classes fondamentales relatives de [B4] et [A.E] gr\^ace au
\Cor{6}{} Sous les hypoth\`eses de la {\it proposition 6}, on a:\par\noindent {\bf(i)}  ${\cal C}_{\pi,\goth{X}}$
d\'efinit un unique \'el\'ement de ${\rm E}xt^{p}({\cal O}_{X},
\Omega^{p}_{Z/S})$ dont l'image par le morphisme canonique ${\rm
E}xt^{p}({\cal O}_{X}, \Omega^{p}_{Z/S})\rightarrow {\rm
H}^{p}_{|X|}(Z, \Omega^{p}_{Z/S})$ est la classe fondamentale
relative pond\'er\'ee  ${\rm C}_{\pi,\goth{X}}$, coincidant, si
$\goth{X}$ est la pond\'eration standard, avec la classe de [B4] (cf
le {\it th\'eor\`eme 0)).\par\noindent {\bf(ii)} Si $\goth{X}$ est la pond\'eration
alg\'ebrique (ce qui correspond au cas o\`u $\pi$ est un morphisme
plat ou plus g\'en\'eralement  de {\it Tor-dimension finie}), ce
morphisme coincide, dans le cadre alg\'ebrique, avec  le morphisme classe fondamentale relative de
[A.E]{\footnote{$^{(10)}$}{Cela signifie exactement que, si
$\pi:X\rightarrow S$ est un morphisme analytiquement g\'eom\'etriquement plat
d'espaces complexes r\'eduits de dimension finie avec $X$ muni d'un
$S$-plongement, $\sigma$,  de codimension $p$ fibre par fibre  dans
$Z$ lisse sur $S$, la construction de [A.E] s'adapte parfaitement
\`a ce cadre et fournit un rel\`evement naturelle de la classe
fondamentale relative pond\'er\'ee \`a valeurs dans $\ds{{\rm
E}xt^{p}({\cal O}_{X}, \Omega^{n+p}_{Z/S})}$.}} \'etendu \`a ce
cadre. On a, de plus, le diagramme commutatif
$$\xymatrix{\sigma_{*}\Omega^{n}_{X/S}\ar[rr]^{{\cal C}_{\pi,
\goth{X}}}\ar[rd]_{\wedge{\rm C}_{\pi,\goth{X}}}&&{\cal
E}xt^{p}({\cal O}_{X}, \Omega^{n+p}_{Z/S})\ar[ld]\\ &{\cal
H}^{p}_{X}(\Omega^{p}_{Z/S})&}$$\rm\bigskip\noindent {\bf 2.1.1.
Preuve de la proposition {\bf 6}.}\smallskip $\blacklozenge$
Supposons $\pi$ analytiquement g\'eom\'etriquement plat. Alors, la {\it proposition 2},  nous conduit \`a  montrer que la
donn\'ee du morphisme d'int\'egration,
$$\int_{\pi,\goth{X}}: {\rm I}\!{\rm R}^{n}{\pi_{!}}\Omega^{n}_{X/S}\rightarrow {\cal O}_{S}$$
dot\'e des propri\'et\'es qui le caract\'erisent, est
enti\`erement conditionn\'ee par la donn\'ee d'un unique morphisme
global de faisceaux ${\cal O}_{X}$-coh\'erents ${\cal C}_{\pi,\goth{X}}:\Omega^{n}_{X/S}\rightarrow {\cal E}xt^{p}({\cal
O}_{X}, \Omega^{p}_{Z/S})$.
\noindent v\'erifiant les propri\'et\'es requises.\par\noindent
Comme il a \'et\'e dit \`a maintes reprises, les donn\'ees relatives \`a $X$,
 $Z$ et $S$ peuvent \^etre compl\'et\'ees et consign\'ees dans un
 diagramme d'installation locale du type $(\spadesuit)$ de la page 37. La finitude de
  $f$ et $\sigma$, l'annulation des faisceaux de cohomologie
  ${\rm I}\!{\rm R}^{j}{\pi_{!}}\Omega^{n}_{X/S}$,   pour tout $j>n$ ($\pi$
  \'etant  \`a fibres de dimension pure $n$) donnent, gr\^ace aux suites spectrales de Leray,  les
isomorphismes $${\rm I}\!{\rm
R}^{j}{\pi_{!}}\Omega^{n}_{X/S}\simeq {\rm I}\!{\rm
R}^{n}{q''_{!}}f_{*}\Omega^{n}_{X/S}\simeq {\rm I}\!{\rm
R}^{n}{q_{!}}\sigma_{*}\Omega^{n}_{X/S}$$ \noindent Alors, la
dualit\'e analytique relative appliqu\'ees aux  projections
canoniques $q$ (resp. $q''$) et l'isomorphisme de Verdier (cf [V],
th\'eor\`eme 3, p.397) ${q}^{!}({\cal O}_{S})\simeq
\Omega^{n+p}_{Z/S}[n+p]$ (resp. ${q''}^{!}({\cal O}_{S})\simeq
\Omega^{n}_{Y/S}[n]$),  donnent la succession d'isomorphismes
$$\xymatrix{{\rm I}\!{\rm H}om({\rm I}\!{\rm R}^{n}{q''_{!}}f_{*}\Omega^{n}_{X/S}, {\cal O}_{S})\ar[r]\ar[d]&{\rm I}\!{\rm H}om({\rm I}\!{\rm R}^{n}{q_{!}}\sigma_{*}\Omega^{n}_{X/S}, {\cal O}_{S})\ar[d]\\
{\rm I}\!{\rm H}om({\rm I}\!{\rm R}{q''_{!}}f_{*}\Omega^{n}_{X/S}[n], {\cal O}_{S})\ar[r]\ar[d]&{\rm I}\!{\rm H}om({\rm I}\!{\rm R}{q_{!}}\sigma_{*}\Omega^{n}_{X/S}[n], {\cal O}_{S})\ar[d]\\
{\rm I}\!{\rm H}om(f_{*}\Omega^{n}_{X/S}[n], {q''}^{!}({\cal O}_{S}))\ar[r]\ar[d]&{\rm I}\!{\rm H}om(\sigma_{*}\Omega^{n}_{X/S}[n], {q}^{!}({\cal O}_{S}))\ar[d]\\
{\rm I}\!{\rm H}om(f_{*}\Omega^{n}_{X/S}[n],
\Omega^{n}_{Y/S}[n])\ar[r]&{\rm I}\!{\rm
H}om(\sigma_{*}\Omega^{n}_{X/S}[n],\Omega^{n+p}_{Z/S}[n+p])}$$
Appliquant, de nouveau,  la dualit\'e analytique pour  un morphisme propre (cf [RRV]) et donc, en particulier, pour un plongement et utilisant  la d\'efinition de $\sigma^{!}$\footnote{$^{(11)}$}{Pour tout faisceau coh\'erent ${\cal F}$ sur $Z$, on a $\sigma_{*}\sigma^{!}{\cal F}:={\rm I}\!{\rm R}hom(\sigma_{*}{\cal O}_{X}, {\cal F})$}, on aboutit, gr\^ace au  {\it lemme 2},
donnant l'annulation  des faisceaux de cohomologie ${\cal H}^{j}({\rm I}\!{\rm R}hom(\sigma_{*}{\cal O}_{X},\Omega^{n+p}_{Z/S}[p] )$ pour tout $j<0$, \`a l'isomorphisme
$$
{\rm I}\!{\rm
H}om(\sigma_{*}\Omega^{n}_{X/S}[n],\Omega^{n+p}_{Z/S}[n+p]
)\simeq{\rm I}\!{\rm H}om(\Omega^{n}_{X/S}, {\cal E}xt^{p}({\cal
O}_{X}, \Omega^{n+p}_{Z/S}))$$  Ainsi, la donn\'ee d'un
morphisme de
 ${\rm I}\!{\rm R}^{n}{\pi_{!}}\Omega^{n}_{X/S}$ \`a valeurs dans
  le faisceau coh\'erent  ${\cal O}_{S}$ est strictement \'equivalente \`a
   la donn\'ee des deux morphismes
    $f_{*}\Omega^{n}_{X/S}\rightarrow \Omega^{n}_{Y/S}$ et
     $\Omega^{n}_{X/S}\rightarrow {\cal E}xt^{p}({\cal O}_{X},
     \Omega^{n+p}_{Z/S})$. \par\noindent D\'esignons par  ${\cal T}^{n}_{f,{\goth
X}}$ et ${\cal C}_{\pi, \goth{X}}$
 les images  de  $\int_{\pi,\goth{X}}$ par ces isomorphismes.
 La {\it proposition 2} et la {\it proposition 4} montrent que
  ${\cal T}^{n}_{f,{\goth X}}$ est un v\'eritable morphisme trace h\'eritant
   de  toutes les propri\'et\'es fonctorielles de $\int_{\pi,\goth{X}}$ et en particulier de la stabilit\'e par changement de base et de la compatibilit\'e
avec les pond\'erations. Il nous reste \`a \'etudier les
propri\'et\'es du morphisme ${\cal C}_{\pi, \goth{X}}:
\Omega^{n}_{X/S}\rightarrow {\cal E}xt^{p}({\cal O}_{X},
 \Omega^{n+p}_{Z/S})$. \smallskip\noindent
 {\bf(i) Propri\'et\'e de la trace relative.}\par\noindent
Comme, d'une part,  ${\cal T}^{n}_{f,{\goth X}}$ est un morphisme
trace prolongeant la trace naturelle
$f_{*}f^{*}\Omega^{n}_{Y/S}\rightarrow
 \Omega^{n}_{Y/S}$ et que, d'autre part, la version locale du point {\bf(v)} du  {\it th\'eor\`eme 1} nous fournit une fl\`eche  ${\rm I}\!{\rm R}^{n}{\pi}_{!}{\cal
E}xt^{p}({\cal O}_{X}, \Omega^{n+p}_{Z/S})\rightarrow {\cal O}_{S}$,  d'ailleurs, strictement \'equivalente \`a la donn\'ee d'un morphisme $f_{*}{\cal
E}xt^{p}({\cal O}_{X}, \Omega^{n+p}_{Z/S})\rightarrow
\Omega^{n+p}_{Y/S}$, on en d\'eduit  un diagramme
$$\xymatrix{f_{*}f^{*}\Omega^{n}_{Y/S}
\simeq f_{*}({\cal O}_{X}\otimes f^{*}\Omega^{n}_{Y/S})\ar[d]_{{\cal
T}^{0}_{f,\goth{X}}\otimes Id} \ar[r]&\!\!\!f_{*}\Omega^{n}_{X/S}
\ar[r]^{\!\!\!\!\!\!\!\!{f_{*}{\cal C}_{\pi,\goth{X}}}}&f_{*}
{\cal E}xt^{p}({\cal O}_{X}, \Omega^{n+p}_{Z/S})\ar[d]^{{\goth T}^{n}_{f,\goth{X}}}\\
\Omega^{n}_{Y/S}\ar[rr]_{Id}&&\Omega^{n}_{Y/S}}$$ dont la
commutativit\'e se v\'erifie simplement aux points g\'en\'eriques de
$X$
 en choisissant une suite r\'eguli\`ere d\'efinissant $X$ autour d'un
  tel point et proc\'edant  \`a un calcul facile en terme de
  symboles. De plus, le faisceau ${\cal E}xt^{p}({\cal O}_{X},
  \Omega^{n+p}_{Z/S})$ est universel pour cette propri\'et\'e dans le sens o\`u
  tout faisceau coh\'erent ${\cal F}$ sur $X$ v\'erifiant
  la propri\'et\'e de la trace relative donne une factorisation
$$\xymatrix{f_{*}f^{*}\Omega^{n}_{Y/S}
\simeq f_{*}({\cal O}_{X}\otimes f^{*}\Omega^{n}_{Y/S})\ar[d]_{{\cal
T}^{0}_{f,\goth{X}}\otimes Id} \ar[r]&\!\!\!f_{*}{\cal F}
\ar[r]^{\!\!\!\!\!\!\!\!{\Phi}}&f_{*}
{\cal E}xt^{p}({\cal O}_{X}, \Omega^{n+p}_{Z/S})\ar[d]^{{\goth T}^{n}_{f,\goth{X}}}\\
\Omega^{n}_{Y/S}\ar[rr]_{Id}&&\Omega^{n}_{Y/S}}$$
{\bf(ii)
Stabilit\'e par changement de base.}\par\noindent Soulignons le fait
que la stabilit\'e par changement de base arbitraire  d'un
 tel morphisme est difficilement r\'ealisable (voir m\^eme irr\'ealisable!)  sans  hypoth\`ese sur $\pi$\footnote{$^{(12)}$}{Le probl\`eme vient du fait que les faisceaux $Ext$ se comportent mal par changement de
base!  Mais c'est le cas si $\pi$  est plat ou si le changement de base est {\it cohomologiquement transverse} \`a $\pi$ (cf [Ba])}. Il s'av\`ere que la
platitude g\'eom\'etrique analytique est suffisamment forte et souple pour permettre une telle propri\'et\'e sur le  faisceau  ${\cal E}xt^{p}({\cal O}_{X},
 \Omega^{n+p}_{Z/S})$. Comme la situation (tout comme le probl\`eme ) est
purement locale, tout morphisme d'espaces complexes r\'eduits de
dimension pure $\eta:S_{1}\rightarrow S$ induit un diagramme
commutatif de changement de base du type\par
\centerline{$\xymatrix{&Z_{1}\ar[ldd]_{q_{1}}\ar[r]^{\phi}&Z\ar[ldd]^{q}&\\
X_{1}\ar[d]_{\pi_{1}}\ar[ru]^{\sigma_{1}}\ar[r]_{\Theta}&X\ar[ru]^{\sigma_{1}}\ar[d]_{\pi}&\\
S_{1}\ar[r]_{\eta}&S&}$} \noindent avec $(Z_{1}, X_{1}, S_{1};
\sigma_{1}, q, \pi_{1})$  v\'erifiant les hypoth\`eses de la
 {\it proposition 6, p.40}. De plus, la notion de platitude
g\'eom\'etrique \'etant stable par changement de base, ${\pi_{1}}$
est aussi analytiquement g\'eom\'etriquement plat muni de  la
pond\'eration $\goth{X}_{1}:=(\Theta\times Id_{S})^{*}\goth{X}$, et
$q$  \'etant lisse, $q_{1}$ l'est aussi.\par\noindent Pour que les
morphismes classes fondamentales relatives pond\'er\'ees soient
stables par changement de base, il suffirait (et il est
n\'ecessaire) que le morphisme  $\Theta^{*}({\cal C}_{\pi,
\goth{X}}):  \Theta^{*}\Omega^{n}_{X/S}\rightarrow \Theta^{*}{\cal
E}xt^{p}({\cal O}_{X},
 \Omega^{n+p}_{Z/S})$ se factorise au travers de
 ${\cal E}xt^{p}({\cal O}_{X_{1}}, \Omega^{n+p}_{Z_{1}/S_{1}})$
 de sorte \`a avoir un  diagramme commutatif $$\xymatrix{\Theta^{*}\Omega^{n}_{X/S}
 \ar[d]\ar[r]^{\Theta^{*}({\cal C}_{\pi,
\goth{X}})}&\Theta^{*}{\cal E}xt^{p}({\cal O}_{X},
 \Omega^{n+p}_{Z/S})\ar[d]\\
 \Omega^{n}_{X_{1}/S_{1}}\ar[r]_{{\cal C}_{\pi_{1}, \goth{X}_{1}}}&{\cal E}xt^{p}({\cal O}_{X_{1}},
 \Omega^{n+p}_{Z_{1}/S_{1}})}$$
o\`u la premi\`ere fl\`eche verticale est un isomorphisme. Ainsi,
tout repose sur l'existence d'une image r\'eciproque
$\Theta^{*}{\cal E}xt^{p}({\cal O}_{X},
\Omega^{n+p}_{Z/S})\rightarrow {\cal E}xt^{p}({\cal O}_{X_{1}},
 \Omega^{n+p}_{Z_{1}/S_{1}})$ qui prolongerait l'image r\'eciproque des formes
 holomorphes $S$-relatives aux points r\'eguliers de $\pi$. Mais
 la succession d'isomorphismes pr\'ec\'edents, la formule de changement de base
  pour l'int\'egration  $\ds{\Theta^{*}(\int_{\pi,\goth X})=
  \int_{\pi_{1},\goth {X_{1}}}}$ et l'isomorphisme
$\Theta^{*}\Omega^{n}_{X/S}\simeq \Omega^{n}_{X_{1}/S_{1}}$ montrent
clairement que $\Theta^{*}$ fait correspondre aux morphismes ${\cal
T}^{n}_{f,\goth{X}}$ et ${\cal C}_{\pi, \goth{X}}$ les morphismes
${\cal T}^{n}_{f_{1},\goth{X}_{1}}$ et ${\cal C}_{\pi_{1},
\goth{X}_{1}}$. Mais la propri\'et\'e de la trace relative
caract\'erisant ce  dernier impose, par universalit\'e, l'existence
d'un unique morphisme $\Theta^{*}: \Theta^{*}{\cal E}xt^{p}({\cal
O}_{X}, \Omega^{n+p}_{Z/S}) \rightarrow {\cal E}xt^{p}({\cal
O}_{X_{1}},
 \Omega^{n+p}_{Z_{1}/S_{1}})$ prolongeant naturellement
  l'image r\'eciproque des formes holomorphes aux points r\'eguliers de $\pi$.
  \smallskip\noindent
En anticipant un peu, en  utilisant  le fait que ${\cal C}_{\pi,
\goth{X}}$ donne naturellement la classe  ${\rm C}_{\pi,\goth{X}}$ du {\it th\'eor\`eme 0}, on se ram\`ene \`a faire la v\'erification  dans le cas de  la pond\'eration standard et donc de la classe fondamentale relative ${\cal C}_{X/S}$ construite dans [B4]
qui est  enti\`erement caract\'eris\'ee par la propri\'et\'e de la
trace relative et stable par changement de base. Alors,
l'\'egalit\'e  $\ds{\Theta^{*}({\cal C}_{X/S})={\cal
C}_{X_{1}/S_{1}}}$ et le fait que
 $\Theta^{*}({\cal C}_{\pi,\goth{X}})$ coincide $S$-g\'en\'eriquement
 avec ${\cal C}_{\pi_{1},\goth{X}_{1} }$ montrent (par universalit\'e)
 qu'il existe n\'ecessairement un unique morphisme
$\Theta^{*}: \Theta^{*}{\cal E}xt^{p}({\cal O}_{X}, \Omega^{n+p}_{Z/S})
\rightarrow {\cal E}xt^{p}({\cal O}_{X_{1}},
 \Omega^{n+p}_{Z_{1}/S_{1}})$ relevant l'image r\'eciproque
 $\Theta^{*}:{\rm H}^{p}_{X}(Z, \Omega^{p}_{Z/S})\rightarrow {\rm H}^{p}_{X_{1}}(Z_{1},
 \Omega^{p}_{Z_{1}/S_{1}})$ et  rendant commutatif le diagramme pr\'ec\'edent.
\smallskip $\blacklozenge$  R\'eciproquement, supposons donn\'e un
morphisme $\pi:X\rightarrow S$ universellement
$n$-\'equidimensionnel, localement install\'e en codimension $p$
dans $Z$ lisse sur $S$, muni d'une pond\'eration $\goth{X}$ et
 d'un morphisme canonique de faisceaux coh\'erents ${\cal
C}_{\pi, \goth{X}}:\Omega^{n}_{X/S}\rightarrow {\cal E}xt^{p}({\cal
O}_{X}, \Omega^{n+p}_{Z/S})$ compatible aux changements de base
entre espaces complexes r\'eduits, \`a l'additivit\'e des
pond\'erations  et v\'erifiant la propri\'et\'e de la trace
relative. Montrons, alors, en vertu de la {\it proposition 2},
qu' il existe un unique morphisme $$\int_{\pi, \goth{X}}:{\rm
I}\!{\rm R}^{n}\pi_{!}\Omega^{n}_{X/S}\rightarrow {\cal O}_{S}$$
v\'erifiant les  propri\'et\'es cit\'ees ci-dessus.
\smallskip\noindent L'un des isomorphismes de la page 39, donn\'e par
 $${\rm I}\!{\rm H}om({\rm I}\!{\rm R}^{n}{\pi_{!}}\Omega^{n}_{X/S}, {\cal O}_{S})\simeq {\rm I}\!{\rm H}om(\Omega^{n}_{X/S}, {\cal E}xt^{p}({\cal O}_{X}, \Omega^{n+p}_{Z/S}))$$
permet d'associer \`a  ${\cal C}_{\pi, \goth{X}}$ un unique morphisme ${\rm I}\!{\rm R}^{n}{\pi_{!}}\Omega^{n}_{X/S}\rightarrow {\cal O}_{S}$ que nous allons noter $\int_{\pi, \goth{X}}$. Pour v\'erifier que cet \'el\'ement poss\`ede les propri\'et\'es requises, il suffit simplement, par exemple,  de se rappeler l'isomorphisme
$${\rm I}\!{\rm H}om({f_{*}}\Omega^{n}_{X/S}, {\Omega^{n}}_{Y/S})\simeq
 {\rm I}\!{\rm H}om(\Omega^{n}_{X/S}, {\cal E}xt^{p}({\cal O}_{X},
 \Omega^{n+p}_{Z/S}))$$
faisant correspondre \`a  ${\cal C}_{\pi, \goth{X}}$ un v\'eritable  morphisme trace  not\'e ${\cal T}^{n}_{f, \goth{X}}$. Il en r\'esulte, pour des raisons \'evidentes d\^ues  \`a la d\'efinition m\^eme de cette fl\`eche d'int\'egration et \`a  l'\'egalit\'e des foncteurs ${\rm I}\!{\rm
R}^{n}{\pi_{!}}\simeq {\rm I}\!{\rm R}^{n}q''_{!}f_{*}$, que  $\int_{\pi,
\goth{X}}$ v\'erifie les propri\'et\'es demand\'ees.
\par\noindent On peut aussi bien remarquer, comme il a \'et\'e dit
plus haut, que la version locale du {\it th\'eor\`eme 1} assure
l'existence d'une fl\`eche
$${\rm I}\!{\rm R}^{n}{\pi}_{!}{\cal
E}xt^{p}({\cal O}_{X}, \Omega^{n+p}_{Z/S})\rightarrow {\cal O}_{S}$$
\noindent
dont on d\'eduit le diagramme commutatif
$$\xymatrix{{\rm I}\!{\rm R}^{n}\pi_{!}\Omega^{n}_{X/S}\ar[rd]\ar[rr]^{{\rm R}^{n}\pi_{!}({\cal C}_{\pi, \goth{X}})}&&{\rm I}\!{\rm R}^{n}\pi_{!}{\cal
E}xt^{p}({\cal O}_{X}, \Omega^{n+p}_{Z/S})\ar[ld]\\
&{\cal O}_{S}&}$$ et les propri\'et\'es voulues pour le morphisme
$\int_{\pi, \goth{X}}: {\rm I}\!{\rm
R}^{n}\pi_{!}\Omega^{n}_{X/S}\rightarrow {\cal O}_{S}$ puisque le
faisceau ${\cal E}xt^{p}({\cal O}_{X}, \Omega^{n+p}_{Z/S})$
v\'erifie, par hypoth\`eses, la propri\'et\'e de la trace relative
et la stabilit\'e par changement de base. Il induira, donc, des
morphismes traces prolongeant les traces naturelles et, de ce fait,
 la compatibilit\'e  \`a l'additivit\'e des pond\'erations.
$\,\,\,\blacksquare$\bigskip\noindent {\bf 2.1.2. Preuve du
corollaire 6.}\smallskip\noindent
Nous allons v\'erifier que  ${\cal
C}_{\pi, \goth{X}}$ contient l'information port\'ee par la classe
fondamentale relative de Barlet (cf [B4]) et qu'il prolonge de fa\c
con naturelle, \`a ce cadre, le  morphisme classe fondamentale
relative  de Ang\'eniol-Elzein (cf [A.E]) construit, dans le cadre
alg\'ebrique, pour un morphisme de {\it Tor-dimension finie}.\par\noindent
Remarquons, d'abord,  que les isomorphismes  ${\rm H}^{p}_{X}(Z, \Omega^{p}_{Z/S})\simeq
\Gamma(Z, {\cal H}^{p}_{X}(\Omega^{p}_{Z/S}))$ et ${\rm  E}xt^{p}({\cal O}_{X},
\Omega^{p}_{Z/S})\simeq \Gamma(Z,{\cal E}xt^{p}({\cal O}_{X}, \Omega^{p}_{Z/S}))$ (cf les lemmes
{\bf 1} et {\bf 2}, p.6-7), permettent de v\'erifier localement, au voisinage de
$X$ dans $Z$, une quelconque assertion portant sur les groupes globaux.\par\noindent
Alors, le morphisme de restriction $\Omega^{n}_{Z/S}\rightarrow \Omega^{n}_{X/S}$, la locale libert\'e des faisceaux $\Omega^{\bullet}_{Z/S}$  et les d\'eg\'en\'erescences  des suites spectrales (cf {\bf (1.0.1)}, p.9)
   $$E^{i,j}_{2}={\cal E}xt^{i}(\Omega^{n}_{Z/S},  {\cal E}xt^{j}
({\cal O}_{X}, \Omega^{n+p}_{Z/S}))\,\,\,{\rm
avec}\,\,E^{i,j}_{2}=0,\,\,\forall\,j<p$$
$$E'^{i,j}_{2}={\cal E}xt^{i}({\cal O}_{X}, {\cal E}xt^{j}
(\Omega^{n}_{Z/S}, \Omega^{n+p}_{Z/S}))\,\,\,{\rm
avec}\,\,E'^{i,j}_{2}=0,\,\,\forall\,j>0$$ d'aboutissement ${\cal
E}xt^{i+j} ({\cal O}_{X}\otimes \Omega^{n}_{Z/S},
\Omega^{n+p}_{Z/S})$,  nous donnent imm\'ediatement un morphisme de faisceaux coh\'erents canonique
$${\cal H}om(\Omega^{n}_{X/S},  {\cal E}xt^{p} ({\cal
O}_{X}, \Omega^{n+p}_{Z/S}))\rightarrow {\cal  E}xt^{p}({\cal O}_{X},
\Omega^{p}_{Z/S})$$
\'equivalent \`a la donn\'ee d'une fl\`eche globalement d\'efinie
$${\rm I}\!{\rm H}om(\Omega^{n}_{X/S},  {\cal E}xt^{p}({\cal
O}_{X}, \Omega^{n+p}_{Z/S}))\rightarrow {\rm  E}xt^{p}({\cal O}_{X},
\Omega^{p}_{Z/S})$$
 \noindent que l'on compose avec le morphisme
naturel ${\rm E}xt^{p}({\cal O}_{X}, \Omega^{p}_{Z/S})\rightarrow
{\rm H}^{p}_{X}(Z, \Omega^{p}_{Z/S})$ envoyant, ainsi,  ${\cal
C}_{\pi,\goth{X}}$ sur un '\'el\'ement $\tilde{{\rm C}}_{\pi,\goth{X}}$. Or, par d\'efinition ou construction m\^eme, il  satisfait {\it la propri\'et\'e de la trace relative}. Cette derni\`ere se v\'erifiant localement, on se ram\`ene \`a un probl\`eme purement ``faisceautique'' en consid\'erant, pour tout germe de $S$-param\'etrisation locale  $f:X\rightarrow Y$, le morphisme d'{\it
int\'egration} de classes de cohomologie \`a support dans $X$ (qui est propre sur $Y$ puisque fini sur $Y$) sur les
fibres $p$-\'equidimensionnelles du morphisme lisse $q':Z\rightarrow
Y$. Plus pr\'ecisemment, on obtient, en tout degr\'e $k$, les fl\`eches
$$\phi:f_{*}{\cal
H}^{p}_{X}(\Omega^{k+p}_{Z/S}) \rightarrow f_{*}{\cal
H}^{p}_{Y-propre}(\Omega^{k+p}_{Z/S})\rightarrow \Omega^{k}_{Y/S}$$
prolongeant  naturellement les traces
usuelles $\ds{{f}_{*}{f}^{*}\Omega^{k}_{Y/S}\rightarrow
\Omega^{k}_{Y/S}}$ et rendent commutatif le diagramme
$$\xymatrix{{f}_{*}\Omega^{k}_{X/S}\ar[rd]^{{\cal T}^{k}_{{f,\goth{X}}}}
\ar[r]^{\!\!\!\!\!\!\!\!\!\wedge {f_{*}}{\rm C}_{\pi, \goth{X}}}&{f}_{*}
{\cal H}^{p}_{X}(\Omega^{p+k}_{Z/S})\ar[d]^{\phi}\\
{f}_{*}{f}^{*}\Omega^{k}_{Y/S}\ar[u]\ar[r]&\Omega^{k}_{Y/S}}$$
Le fait que $\tilde{{\rm C}}_{\pi,\goth{X}}$  soit stable par changement de base entre espaces complexes r\'eduits, compatible   aux localisations sur $X$ et \`a l'additivit\'e des pond\'erations, se d\'eduit directement des propri\'et\'es du morphisme  ${\cal C}_{\pi,\goth{X}}$ \`a l'aide d'une v\'erification locale faite aux points g\'en\'eriques de $X$ (ou sur un quelconque ouvert de $X$ dense fibre par fibre sur $S$!).\par\noindent
On en conclut que $\tilde{{\rm C}}_{\pi,\goth{X}}$ est exactement la classe ${\rm  C}_{\pi,\goth{X}}$ donn\'ee dans le {\it th\'eor\`eme 0} et que, par cons\'equent, elle coincide, dans le cas de la pond\'eration standard, avec  la classe fondamentale relative ${\cal C}_{X/S}$ construite dans [B4].
\smallskip\noindent
En vertu du {\it corollaire 5} (p.36), tout morphisme
$\pi:X\rightarrow S$  \`a fibres de dimension pure constante $n$ et
de {\it Tor-dimension finie} sur une base r\'eduite est
analytiquement g\'eom\'etriquement plat. Par cons\'equent, il lui
est associ\'e un unique morphisme de faisceaux coh\'erents ${\cal
C}_{\pi,\goth{X}}$ v\'erifiant les propri\'et\'es requises. Mais notre situation locale s'apparente \`a celle o\`u $S$ est un sch\'ema localement noeth\'erien excellent, $Z$ un $S$-sch\'ema localement noeth\'erien, localement de type fini et lisse sur $S$ et $\sigma$ est une immersion ferm\'ee induisant sur $X$ une projection de {\it Tor-dimension finie} \`a fibres de dimension constante  $n$ sur $S$ apparaissant, en un certain sens, comme un cas particulier de la situation g\'en\'erale d\'ecrite  par un morphisme localement de type fini, \'equidimensionnel et de  {\it Tor-dimension finie} entre  sch\'emas  localement noeth\'eriens (sans r\'ef\'erence \`a un quelconque plongement). Comparant les construction et invoquant le principe d'unicit\'e, on en d\'eduit que [A.E] qui s'adaptent parfaitement \`a cette situation {\bf locale} et que  ${\cal C}_{\pi,\goth{X}}$ est l'analogue analytique
local de la classe de [A.E].\par\noindent
Enfin, il est facile de faire  la connexion entre
  ces objets  gr\^ace au  diagramme
$$\xymatrix{\Omega^{n}_{X/S}\ar[rr]^{{\cal C}_{\pi,
\goth{X}}}\ar[rd]_{{\cal C}_{X/S}}&&{\cal E}xt^{p}({\cal O}_{X},
\Omega^{n+p}_{Z/S})\ar[ld]\\
&{\cal H}^{p}_{X}(\Omega^{n+p}_{Z/S})&}$$\noindent qui est
commutatif puisque ${\cal C}_{X/S}$ est  annul\'ee, au moins
g\'en\'eriquement sur $S$, par tout id\'eal  de d\'efinition de $X$
dans $Z$. Il apparait que le morphisme   ${\cal C}_{\pi,\goth{X}}$
est l'unique  rel\`evement du morphisme cup-produit
 $\wedge {\cal C}_{X/S}$ $\,\,\,\blacksquare$\bigskip\noindent
{\bf 2.1.3. Remarque.} \par\noindent
{\bf(i)} Signalons au lecteur  que   les arguments
globaux utilis\'es dans [A.E] ne conviennent pas du tout ou tout
simplement n'ont pas d'\'equivalents en g\'eom\'etrie analytique
complexe.\par\noindent
{\bf(ii)} Il est \'evident que ${\rm C}_{\pi, \goth{X}}$ induit, en tout degr\'e $k$, un morphisme de faisceaux de groupes ab\'eliens
${\rm C}^{k}_{\pi, \goth{X}}:\Omega^{k}_{X/S}\rightarrow {\cal H}^{p}_{X}(\Omega^{k+p}_{Z/S})$ produisant un morphisme de complexes diff\'erentiels 
${\rm C}^{\bullet}_{\pi, \goth{X}}:\Omega^{\bullet}_{X/S}\rightarrow {\cal H}^{p}_{X}(\Omega^{p+\bullet}_{Z/S})$ munis de la diff\'erentielle $S$-relative usuelle.
\bigskip\noindent
{\bf 2.1.4. Preuve
de la proposition 5.}\smallskip\noindent Tout comme dans la
 {\it proposition 6, p.40}, nous allons utiliser l'int\'egration
globale de  la {\it proposition 1} et montrer que la donn\'ee
du morphisme
$$\int_{\pi,\goth{X}}: {\rm I}\!{\rm R}^{n}{\pi_{!}}\Omega^{n}_{X/S}\rightarrow {\cal O}_{S}$$
dot\'e des propri\'et\'es \'enonc\'ees dans cette  proposition, est
enti\`erement conditionn\'ee par la donn\'ee d'un unique morphisme
global de faisceaux ${\cal O}_{X}$-coh\'erents\par
\centerline{${\cal C}_{\pi,\goth{X}}:\Omega^{n}_{X/S}\rightarrow
\omega^{n}_{\pi}$} \noindent v\'erifiant les propri\'et\'es
voulues.\smallskip $\blacklozenge$ Pour ce faire, on se ram\`ene \`a
la situation locale en  choisissant  un recouvrement ouvert
localement fini $(X_{\alpha})_{\alpha\in A}$ adapt\'e \`a la
localisation du morphisme $\pi$ autour de l'une de ses fibres
$\pi^{-1}(s_{0}):=X_{s_{0}}$, c'est-\`a-dire des diagrammes de
localisations commutatifs du type $(\spadesuit)_{\alpha}$ de la page
39. Posons $\omega^{n}_{\pi_{\alpha}}:= \sigma^{*}_{\alpha}{\cal
E}xt^{p}({\sigma_{\alpha}}_{*}{\cal O}_{X_{\alpha}},
\Omega^{m}_{Z_{\alpha}/S})$ et remarquons qu'en vertu de la nature
locale sur $X$ et $S$ de la notion de platitude g\'eom\'etrique
analytique, chaque morphisme $\pi_{\alpha}$ est analytiquement
g\'eom\'etriquement plat.\par\noindent La version faisceautique du
th\'eor\`eme de Reiffen donne la suite exacte \`a droite
$$ \bigoplus_{\alpha\in A}{\rm I}\!{\rm R}^{n}{\pi_{\alpha}}_{!}
\Omega^{n}_{{X_{\alpha}}/S} \rightarrow {\rm I}\!{\rm
R}^{n}{\pi_{!}}{\Omega^{n}_{X/S}}\rightarrow 0$$
 de laquelle on d\'eduit la fl\`eche injective
$${\rm H}om({\rm I}\!{\rm R}^{n}{\pi_{!}}{\Omega^{n}_{X/S}}, {\cal O}_{S})\rightarrow {\rm H}om(\bigoplus_{\alpha}{\rm I}\!{\rm R}^{n}{\pi_{\alpha}}_{!}\Omega^{n}_{{X_{\alpha}}/S},  {\cal O}_{S})$$
Par d\'efinition du produit et coproduit cat\'egoriel, on a
l'isomorphisme
$${\rm I}\!{\rm H}om(\bigoplus_{\alpha}{\rm I}\!{\rm R}^{n}{\pi_{\alpha}}_{!}\Omega^{n}_{{X_{\alpha}}/S},
  {\cal O}_{S})\simeq\prod_{\alpha}
{\rm I}\!{\rm H}om({\rm I}\!{\rm
R}^{n}{\pi_{\alpha}}_{!}\Omega^{n}_{{X_{\alpha}}/S},  {\cal O}_{S})
$$
Des arguments d\'ej\`a \'evoqu\'es dans la preuve de la {\it proposition 6}  (isomorphisme de Verdier et la dualit\'e
analytique relative pour une projection et un morphisme propre ) donnent, en particulier,  l'isomorphisme
$${\rm I}\!{\rm H}om({\rm I}\!{\rm
R}^{n}{\pi_{\alpha}}_{!}\Omega^{n}_{{X_{\alpha}}/S},  {\cal
O}_{S})\simeq {\rm I}\!{\rm H}om(\Omega^{n}_{{X_{\alpha}}/S},
\omega^{n}_{\pi_{\alpha}})$$ Cela nous conduit vers le diagramme commutatif \footnote{$^{(9)}$}{Il est tout \`a fait possible de "topologiser" la
situation en  consid\'erant  le foncteur ${\rm I}\!{\rm H}om()$ sur les
morphismes lin\'eaires continus et obtenir des isomorphismes. Mais
cela fait l'objet d'un article ult\'erieur sur la notion de paire
dualisante en g\'eom\'etrie analytique dont le d\'ecor naturel est
un aspect de la dualit\'e analytique relative.}
$$\xymatrix{0\ar[r]&{\rm I}\!{\rm H}om({\rm I}\!{\rm R}^{n}\pi_{!}\Omega^{n}_{X/S},{\cal
O}_{S})\ar[r]\ar[d]_{\Psi}&{\ds{\prod_{\alpha}}} {\rm I}\!{\rm
H}om({\rm I}\!{\rm
R}^{n}{\pi_{\alpha}}_{!}\Omega^{n}_{{X_{\alpha}}/S},  {\cal
O}_{S})\ar[d]^{\ds{\prod_{\alpha}}\Psi_{\alpha}}\\
&{\rm I}\!{\rm H}om(\Omega^{n}_{X/S}, \omega^{n}_{\pi})\ar[r] &
{\ds{\prod_{\alpha}}}{\rm I}\!{\rm H}om(\Omega^{n}_{{X_{\alpha}}/S},
\omega^{n}_{\pi_{\alpha}})&}$$  \par\noindent Ainsi,
$\int_{\pi, \goth{X}}$ donne une collection $(\int_{\pi_{\alpha}, \goth{X}_{\alpha}})_{\alpha}$ de morphismes d'int\'egration locales
 avec  $\int_{\pi_{\alpha},\goth{X}_{\alpha}}: {\rm I}\!{\rm
R}^{n}{\pi_{\alpha}}_{!}\Omega^{n}_{{X_{\alpha}}/S}\rightarrow {\cal
O}_{S}$ h\'eritant de toutes les propri\'et\'es de $\int_{\pi,
\goth{X}}$ et produisant naturellement un morphisme de faisceaux
${\cal O}_{X_{\alpha},\goth{X}_{\alpha}}$-coh\'erents ${\cal
C}_{\pi_{\alpha}}: \Omega^{n}_{X_{\alpha}/S}\rightarrow
\omega^{n}_{\pi_{\alpha}}$ v\'erifiant les propri\'et\'es de la
{\it proposition 6, p.40}. Mais le  {\it th\'eor\`eme 1} nous dit que cette
famille, $(\omega^{n}_{\pi_{\alpha}})_{\alpha\in A} $, de faisceaux
${\cal O}_{X_{\alpha}}$-coh\'erents,  se recollent, \`a isomorphisme
canonique pr\`es, en un unique  faisceau ${\cal O}_{X}$-coh\'erent
$\omega^{n}_{\pi}$. Il en r\'esulte aussit\^ot  que la famille de
morphismes $({\cal C}_{\pi_{\alpha},\goth{X}_{\alpha}})_{\alpha\in
A}$ se recollent globalement sur $X$, en un unique  morphisme ${\cal
C}_{\pi, \goth{X}}: \Omega^{n}_{X/S}\rightarrow \omega^{n}_{\pi}$.
Alors, la {\it proposition 6, p.40} et la nature locale du faisceau $\omega^{n}_{\pi}$
 donn\'e par le {\it th\'eor\`eme 1}  permettent de v\'erifier
 les principales propri\'et\'es \`a savoir:\par
$\bullet$ {\it la propri\'et\'e de la trace relative}, dont la v\'erification est locale et r\'esulte donc de la {\it proposition 6, p.40}.\par
$\bullet$ {\it  la compatibilit\'e avec les changements de bases},
traduisant le fait que si  $\eta:S_{1}\rightarrow S$ est  un
morphisme d'espaces complexes r\'eduits de dimension pure induisant le diagramme de changement de base \par
\centerline{$\xymatrix{X_{1}\ar[d]_{\pi_{1}}\ar[r]^{\Theta}&X\ar[d]_{\pi}\\
S_{1}\ar[r]_{\eta}&S}$} \noindent
 alors, il existe un morphisme  canonique
$\ds{\Theta^{*}\omega^{n}_{\pi}\rightarrow \omega^{n}_{\pi_{1}}}$
prolongeant celui des formes holomorphes relatives sur la partie
r\'eguli\`ere de $\pi$) et rendant commutatif le diagramme\par
\centerline{$\xymatrix{\Theta^{*}(\Omega^{n}_{X/S})\ar[r]^{\Theta^{*}({\cal C}_{\pi})}\ar[d]&\Theta^{*}\omega^{n}_{\pi}\ar[d]\\
 \Omega^{n}_{X_{1}/S_{1}}\ar[r]_{{\cal C}_{\pi_{1}}}&\omega^{n}_{\pi_{1}}}$}
\noindent
Insistons encore sur le fait que, en g\'en\'eral, le faisceau $\omega^{n}_{\pi}$ se comporte mal par changement de
base et pour qu'il en soit ainsi il faut une raison tr\`es s\'erieuse!\par\noindent
 Comme la notion de platitude g\'eom\'etrique analytique
est stable par changement de base entre espaces complexes r\'eduits,
$\pi_{1}$ est aussi analytiquement g\'eom\'etriquement plat. Par
cons\'equent, il lui est associ\'e un unique morphisme $\ds{{\cal
C}_{\pi_{1}, \goth{X}_{1}}:\Omega^{n}_{X_{1}/S_{1}}\rightarrow
\omega^{n}_{\pi_{1}}}$ v\'erifiant les propri\'et\'es {\bf(i)} et
{\bf(iii)} de la {\it proposition 1} en vertu de ce qui pr\'ec\`ede.\par\noindent La nature
locale du probl\`eme permet de se ramener \`a une situation
plong\'ee  compl\'etant le diagramme pr\'ec\'edent en un diagramme
commutatif faisant apparaitre un $S$-plongement local de $X$ et son
transform\'e par le changement de base donn\'e. Ainsi, sans perte de
g\'en\'eralit\'e, on se retrouve dans le contexte de [B4] dans
lequel on avait $\ds{\Theta^{*}({\cal C}_{X/S})={\cal
C}_{X_{1}/S_{1}}}$ traduisant la stabilit\'e par changement de base
entre espaces complexes r\'eduits  de la classe fondamentale
relative. \par\noindent De plus, l'\'egalit\'e
$\Theta^{*}(\int_{\pi, \goth{X}}):=\int_{\pi_{1}, \goth{X}_{1}}$, le fait que
$\Theta^{*}({\cal C}_{\pi,\goth{X}})$ coincide  $S$-g\'en\'eriquement
 avec ${\cal C}_{\pi_{1},\goth{X}_{1} }$ montrent (par universalit\'e) qu'il existe
 n\'ecessairement un unique morphisme
$\ds{\Theta^{*}\omega^{n}_{\pi}\rightarrow \omega^{n}_{\pi_{1}}}$
rendant commutatif le diagramme pr\'ec\'edent.\par\noindent En vertu
de l'abscence de torsion et vu qu'ils coincident sur un ouvert dense
fibre par fibre de $X$ (par exemple sur l'ouvert Cohen-Macaulay de
$\pi$), ce morphisme est injectif. \par\noindent Comme les
morphismes lisses sont stables par changement de base, il est clair
que cette image r\'eciproque prolonge  l'image r\'eciproque des
formes holomorphes $S$-relatives.\smallskip
 $\bullet$ {\it la compatibilit\'e avec la structure de complexes
diff\'erentiels.}\par\noindent
 On commence par remarquer que le morphisme ${\cal
C}_{\pi, \goth{X}}:\Omega^{n}_{X/S}\rightarrow\omega^{n}_{\pi}$
donne naturellement, en tout degr\'e $k$,  deux types de  morphismes canoniques , \`a savoir :\par {\bf(i)}
$\Omega^{k}_{X/S}\rightarrow\omega^{k}_{\pi}$\smallskip {\bf(ii)}
$\omega^{k}_{\pi}\rightarrow {\cal E}xt^{p}({\cal O}_{X},
\Omega^{k+p}_{Z/S})$\smallskip\noindent En effet, appliquant le
foncteur covariant ${\cal H}om(\Omega^{n-k}_{X/S}, -)$ \`a ${\cal
C}_{\pi, \goth{X}}$ et utilisant le morphisme naturel
$\Omega^{k}_{X/S}\rightarrow {\cal H}om(\Omega^{n-k}_{X/S},
\Omega^{n}_{X/S})$, on obtient une fl\`eche canonique
$$\Omega^{k}_{X/S}\rightarrow {\cal H}om(\Omega^{n-k}_{X/S},
\omega^{n}_{\pi})$$ Mais en vertu  du  {\it corollaire (1.1)}, on a  l'isomorphisme $\omega^{k}_{\pi}\simeq {\cal
H}om(\Omega^{n-k}_{X/S}, \omega^{n}_{\pi})$ et par suite le
morphisme donn\'e dans {\bf(i)}.\par\noindent Le second s'obtient en
utilisant les suites spectrales maintes fois utilis\'ees (cf par exemple p.46). Sans
d\'etailler, disons que l'on commence par \'ecrire (pour un
plongement local $\sigma$) la fl\`eche $${\cal
H}om(\sigma_{*}\Omega^{n-k}_{X/S}, {\cal E}xt^{p}(\sigma_{*}{\cal
O}_{X}, \Omega^{n+p}_{Z/S}))\rightarrow {\cal
H}om(\Omega^{n-k}_{Z/S}, {\cal E}xt^{p}(\sigma_{*}{\cal O}_{X},
\Omega^{n+p}_{Z/S})),$$ \noindent obtenue gr\^ace \`a l'incarnation
locale de $\omega^{n}_{\pi}$ et \`a la restriction naturelle
$\ds{\Omega^{n-k}_{Z/S}\rightarrow\sigma_{*}\Omega^{n-k}_{X/S}}$,
puis utilisant $\Omega^{\bullet}_{Z/S}\simeq {\cal
H}om(\Omega^{n+p-\bullet}_{Z/S}, \Omega^{n+p}_{Z/S})$ et
les annulations du {\it lemme 2, p.7}  (puisque  $X$
est de codimension pure $p$ fibre par fibre dans $Z$ qui est lisse
sur $S$), on obtient l'isomorphisme
$${\cal E}xt^{p}(\sigma_{*}{\cal
O}_{X}\otimes\Omega^{n-k}_{Z/S} , \Omega^{n+p}_{Z/S})\simeq {\cal
E}xt^{p}(\sigma_{*}{\cal O}_{X}, \Omega^{p+k}_{Z/S})$$
mettant en \'evidence le morphisme d\'esir\'e qui est manifestement
injectif et rendant  commutatif le diagramme
$$\xymatrix{\Omega^{k}_{X/S}\ar[rr]\ar[rd]&&\omega^{k}_{\pi}\ar[ld]\\
&{\cal E}xt^{p}({\cal O}_{X}, \Omega^{p+k}_{Z/S})&}$$
\smallskip\noindent Remarquons que pour $\alpha\in {\cal H}om(\Omega^{n-k}_{X/S}, \omega^{n}_{\pi})$ et $\beta\in \Gamma(X,\Omega^{j}_{X/S})$, on d\'efinit un \'el\'ement $\alpha\otimes\beta\in {\cal H}om(\Omega^{n-k-j}_{X/S}, \omega^{n}_{\pi})$, par la formule standard \par
\centerline{$(\alpha\otimes\beta)(\xi):=\alpha(\beta\wedge\xi),\,\,\,\,\forall\,\xi\in\Gamma(X, \Omega^{n-k-j}_{X/S})$}
\noindent La diff\'erentielle ext\'erieure relative usuelle du  complexe $\Omega^{\bullet}_{X/S}$ n'\'etant pas ${\cal
O}_{X}$-lin\'eaire, n'induira pas de fa\c con naive une diff\'erentielle sur  $\omega^{\bullet}_{\pi}$, comme on peut le constater en regardant [E] ou  la construction explicite de [Ke].
\smallskip\noindent
N'ayant pas, dans notre contexte, de complexe dualisant relatif ni de complexe r\'esiduel relatif, on aborde  la question en adoptant (avec pr\'ecaution) [B2] \`a notre situation relative.\par\noindent
Soit $x$ un point de $X$. Alors pour tout syst\`eme de fonctions
holomorphes ${\bf f}:=(f_{1},\dots,f_{p})$ sur $Z$, s'annulant sur
$X$ et donnant des \'equations r\'eduites en $x$ dans $Z$ au moins
g\'en\'eriquement sur $S$,  on a une suite exacte courte des germes
de faisceaux coh\'erents en $x$
$$\xymatrix{0\ar[r]&\omega^{k}_{\pi}\ar[r]&
{\cal E}xt^{p}({\cal O}_{X}, \Omega^{k+p}_{Z/S})\ar[r]^{\phi}&({\cal
E}xt^{p}({\cal O}_{X}, \Omega^{k+p+1}_{Z/S}))^{p}}$$ o\`u la
premi\`ere fl\`eche injective (pr\'ec\'edemment d\'ecrite)  n'est
rien d'autre que celle induite par le cup produit avec ${\cal
C}_{\pi, \goth{X}}$ et o\`u
 la seconde est donn\'ee par
$\phi(\alpha):=(\alpha\wedge df_{1},\cdots,\alpha\wedge df_{p})$.
Cette repr\'esentation (\`a priori d\'ependante des choix de la
suite r\'eguli\`ere ${\bf f}$) permet  d'identifier les sections de
$\omega^{k}_{\pi}$ aux sections de ${\cal
H}^{p}_{X}(\Omega^{p+k}_{Z/S})$ annul\'ees par $f$ et $df$. En
d'autre terme $\ds{\omega^{k}_{\pi}\simeq {\cal E}xt^{p} ({\cal
O}_{X}, \Omega^{k+p}_{Z/S})\cap {\rm Ker}(\wedge df})$.\par\noindent
Pour s'assurer de l'ind\'ependance de cette repr\'esentation
vis-\`a-vis de la suite r\'eguli\`ere choisie, on peut soit invoquer
le caract\`ere canonique  du morphisme classe fondamentale relative
${\cal C}_{\pi, \goth{X}}$ ou se r\'ef\'erer aux lemmes d'invariances (non
triviaux!) \'etablis par Kersken dans [Ke]. Ainsi, la diff\'erentielle
usuelle $d':{\cal H}^{p}_{X}(\Omega^{p+k}_{Z/S})\rightarrow {\cal H}^{p}_{X}(\Omega^{p+k+1}_{Z/S})$, rendant commutatif les diagrammes induits par les cup-produits par la classe fondamentales
$$\xymatrix{\Omega^{k}_{X/S}\ar[r]\ar[d]_{d_{X/S}}&{\cal H}^{p}_{X}(\Omega^{p+k}_{Z/S})\ar[d]^{d'}\\
\Omega^{k+1}_{X/S}\ar[r]&{\cal H}^{p}_{X}(\Omega^{p+k+1}_{Z/S})},$$
induit naturellement une diff\'erentielle
$\ds{D:\omega^{k}_{\pi}\rightarrow \omega^{k+1}_{\pi}}$ faisant de $
(\omega^{\bullet}_{\pi}, D)$ un complexe diff\'erentiel de
$(\Omega^{\bullet}_{X/S},d_{X/S})$-modules (gr\^ace aux morphismes traces).
De plus, elle satisfait la relation classique \par
\centerline{$D(\alpha\otimes\beta):=\alpha\otimes d'\beta \pm D\alpha\otimes \beta$}\noindent
pour tout $\alpha\in {\cal H}om(\Omega^{n-k}_{X/S}, \omega^{n}_{\pi})$ et $\beta\in \Gamma(X,\Omega^{j}_{X/S})$. Alors, la remarque {\bf(2.1.3)} (p.47)  suivant la {\it proposition 6} ou son {\it corollaire 6} (p.40) met en \'evidence, en vertu de ce qui pr\'ec\`ede, le  morphisme gradu\'e  de complexes diff\'erentiels
$${\cal C}^{\bullet}_{\pi, \goth{X}}: \Omega^{\bullet}_{X/S}\rightarrow
\omega^{\bullet}_{\pi}$$
\indent $\blacklozenge$ R\'eciproquement, soit $\pi:X\rightarrow S$
un morphisme universellement $n$-\'equidimensionnel muni d'une pond\'eration $\goth{X}$ et  d'un
morphisme canonique de faisceaux coh\'erents ${\cal
C}_{\pi, \goth{X}}:\Omega^{n}_{X/S}\rightarrow \omega^{n}_{\pi}$ compatible
aux changements de base r\'eduit et v\'erifiant la propri\'et\'e de
la trace relative. Alors, il existe un unique morphisme $$\int_{\pi, \goth{X}}:{\rm
I}\!{\rm R}^{n}\pi_{!}\Omega^{n}_{X/S}\rightarrow {\cal O}_{S}$$
v\'erifiant les  propri\'et\'es requises.\par\noindent En effet,
il r\'esulte des hypoth\`eses que le faisceau $\omega^{n}_{\pi}$ est n\'ecessairement stable par changement de base et v\'erifie la propri\'et\'e de la trace relative. De plus, le {\it th\'eor\`eme 1}  nous dit qu' il munit $\pi$ d'un morphisme  canonique $$\int_{\pi}:{\rm I}\!{\rm
R}^{n}\pi_{!}\omega^{n}_{\pi}\rightarrow {\cal O}_{S}$$
de formation compatible aux restrictions ouvertes sur $X$ et $S$. le morphisme  ${\cal C}_{\pi,\goth{X}}$ induira, alors, le  diagramme commutatif
$$\xymatrix{{\rm I}\!{\rm R}^{n}\pi_{!}\Omega^{n}_{X/S}
\ar[rd]_{\int_{\pi, \goth{X}}}\ar[rr]^{{\rm I}\!{\rm R}^{n}\pi_{!}({\cal C}_{\pi,
\goth{X}})} &&{\rm I}\!{\rm R}^{n}\pi_{!}\omega^{n}_{\pi}\ar[ld]^{\int_{\pi}}\\
&{\cal O}_{S}&}$$
C'est un exercice facile que de montrer que le  morphisme
$\int_{\pi, \goth{X}}$  satisfait toutes les propri\'et\'es voulues. Il suffit pour cela de les d\'eduire de celles de ${\cal C}_{\pi,\goth{X}}$. En effet, le proc\'ed\'e de d\'ecoupage et le lemme de Reiffen nous ram\`ene \`a la situation locale. Ainsi, pour des installations de type $(\spadesuit)_{\alpha}$ et gr\^ace \`a  ${\cal C}_{\pi, \goth{X}}$ (se localisant en ${\cal C}_{\pi_{\alpha}, \goth{X}_{\alpha}}$), on produit, pour chaque
$\alpha$,  une famille de
morphismes  gradu\'es $\ds{{\cal T}^{k}_{f_{\alpha}}:
{f_{\alpha}}_{*}\Omega^{k}_{X_{\alpha}/S}\rightarrow\Omega^{k}_{Y_{\alpha}/S}}$, donnant, en particulier pour $k=0$, les morphismes\par
\centerline{$\xymatrix{{f_{\alpha}}_{*}{\cal
O}_{X_{\alpha}}\ar@/_/[rr]_{{\cal
T}^{0}_{f_{\alpha}}}\ar[r]^{\!\!\!\!\!\!\!\!\!\wedge
 {\cal C}_{\pi_{\alpha}, \goth{X}_{\alpha}}}&{f_{\alpha}}_{*}{\cal H}^{p}_{|X_{\alpha}|}(\Omega^{p}_{Z_{\alpha}/S})\ar[r]^{\phi}&{\cal O}_{Y_{\alpha}}}$}\noindent
 qui sont des morphismes de type trace et assurent  de ce fait la platitude
  g\'eom\'etrique analytique de chaque $f_{\alpha}$. Il en r\'esulte alors, pour tout entier $k$ et tout indice $\alpha$, le morphisme trace usuelle
 ou naturel  $\ds{{\cal T}^{k}_{f_{\alpha}}:{f_{\alpha}}_{*}{f_{\alpha}}^{*}\Omega^{k}_{Y_{\alpha}/S}\rightarrow \Omega^{k}_{Y_{\alpha}/S}}$, rendant
commutatif le diagramme
$$\xymatrix{f_{*}\Omega^{k}_{X/S}\ar[rd]^{{\cal T}^{k}_{f}}\ar[r]^{\!\!\!\!\!\!\!\!\!\wedge {\cal C}_{\pi}}&f_{*}{\cal H}^{p}_{|X|}(\Omega^{p+k}_{Z/S})\ar[d]^{\phi}\\
f_{*}f^{*}\Omega^{k}_{Y/S}\ar[u]\ar[r]&\Omega^{k}_{Y/S}}$$ Comme ces
 op\'erations sont  stables par changement de base entre espaces
complexes r\'eduits, compatible aux localisations sur $X$ et \`a
l'additivit\'e des pond\'erations (qui, rappelons le, est une
propri\'et\'e inh\'erente aux morphismes trace), la {\it
proposition 4} permet de conclure. $\,\,\blacksquare$\smallskip\noindent
{\bf 2.1.5. Preuve du corollaire 2.1}\par\noindent
Comme $\pi_{1}$ et $\pi_{2}$ sont analytiquement g\'eom\'etriquement plats,
 ils sont, \`a fortiori, universellement \'equidimensionnels et donc leur
 compos\'ee $\pi$ aussi; la puret\'e de la dimension de ses  fibres est
  une cons\'equence de la propri\'et\'e d'ouverture que poss\`edent
  $\pi_{1}$ et $\pi_{2}$ (cf [KI], {\bf(1.2.1.2)}, p.34).  Il est facile de d\'efinir une pond\'eration contin\^ument
  g\'eom\'etriquement plate sur $\pi$.
  En effet, soit ${\goth{X}}_{1}$ (resp.  ${\goth{X}}_{2}$)
   la pond\'eration analytiquement g\'eom\'etriquement plate de $\pi_{1}$
   (resp. $\pi_{2}$) et consid\'erons le morphisme analytiquement
   g\'eom\'etriquement plat  $\Theta:=\pi_{2}\times id:X_{2}\times
   S\rightarrow X_{1}\times S$. Alors, utilisant la d\'efinition de
   l'image r\'eciproque d'un cycle par un morphisme analytiquement
   g\'eom\'etriquement plat donn\'ee dans [KI], {\bf(1.3.8.2)}, p.46, on voit que
   $\goth{X}:= \Theta^{*}({\goth{X}}_{1})$ est un cycle
    parfaitement d\'efini (la condition d'incidence $\Theta^{-1}(|\goth{X}_{1}|)\cap(\{s\}\times X_{2})=\pi^{-1}(s)$  est  satisfaite comme on peut
facilement s'en convaincre) et  fournit  une pond\'eration de $\pi$.
De plus, elle est contin\^ument g\'eom\'etriquement plate sans
 aucune condition sur $S$. En effet, la preuve de la {\it
 proposition 1}, p.47 de [KI] montre que, pour une certaine famille
 de cartes  sur $X_{2}$ ("adapt\'ees" \`a celles de $X_{1}$), les
 morphismes classifiants associ\'es \`a $\pi$ sont des applications
 holomorphes ou, ce qui revient au m\^eme, que les traces
 associ\'ees sont $S$-holomorphes. Il suffit, alors, d'invoquer le
 th\'eor\`eme du changement de projection (cf [B1], {\it
 th\'eor\`eme 2}, p.42) disant que le passage d'une carte \`a une autre est toujours
  m\'eromorphe continue pour en d\'eduire la continuit\'e de la
  famille de cycles $([\pi^{-1}(s)])_{s\in S}$. En effet, cette m\^eme proposition et le  {\it corollaire 2} ( p.34),  montre que, relativement
 \`a toute param\'etrisation locale $S$-relative
 $f:X_{2}\rightarrow S\times V$, le morphisme
  trace est \`a valeurs
 dans le faisceau ${\cal
O}^{c}_{S}{\widehat\otimes}\Omega^{n_{1}+n_{2}}_{V}$. On conatate, au passage,   que si $S$ est faiblement normal (i.e ${\cal O}^{c}_{S}={\cal
O}_{S}$), le morphisme est
 analytiquement g\'eom\'etriquement plat.
\par\noindent En vertu du  {\it th\'eoreme 2}, $\pi$ sera
analytiquement g\'eom\'etriquement plat si et seulement si  on
arrive \`a produire un morphisme d'int\'egration $\int_{\pi,
\goth{X}}:{\rm I}\!{\rm
R}^{n_{1}+n_{2}}{\pi}_{!}\Omega^{n_{1}+n_{2}}_{X_{2}/S}\rightarrow{\cal
O}_{S}$ ayant les propri\'et\'es fonctorielles cit\'ees  dans ce
th\'eor\`eme. Pour cela, on dispose d'un morphisme naturel
d'int\'egration (adapt\'e \`a cette pond\'eration) ${\rm I}\!{\rm
R}^{n_{1}+n_{2}}{\pi}_{!}( \Omega^{n_{2}}_{X_{2}/X_{1}} \otimes
\pi_{2}^{*}\Omega^{n_{1}}_{X_{1}/S})\rightarrow{\cal O}_{S}$ et
d'une fl\`eche naturelle $\Omega^{n_{2}}_{X_{2}/X_{1}}
\otimes\pi_{2}^{*}\Omega^{n_{1}}_{X_{1}/S}\rightarrow
\Omega^{n_{1}+n_{2}}_{X_{2}/S}$. Il parait, alors, clairement, que
la platitude g\'eom\'etrique analytique de $\pi$ est enti\`erement
conditionn\'ee par l'existence  d'un  unique prolongement du
morphisme naturel  d'int\'egration \`a  ${\rm I}\!{\rm
R}^{n_{1}+n_{2}}{\pi}_{!}\Omega^{n_{1}+n_{2}}_{X_{2}/S}$, via la
fl\`eche naturelle $\Omega^{n_{2}}_{X_{2}/X_{1}}
\otimes\pi_{2}^{*}\Omega^{n_{1}}_{X_{1}/S}\rightarrow
\Omega^{n_{1}+n_{2}}_{X_{2}/S}$, de sorte \`a avoir le diagramme
commutatif
$$\xymatrix{ {\rm I}\!{\rm R}^{n_{1}+n_{2}}{\pi}_{!}(\Omega^{n_{2}}_{X_{2}/X_{1}} \otimes
\pi_{2}^{*}\Omega^{n_{1}}_{X_{1}/S} )\ar[rd]\ar[rr]&&
 {\rm I}\!{\rm R}^{n_{1}+n_{2}}{\pi}_{!}\Omega^{n_{1}+n_{2}}_{X_{2}/S}\ar[ld]\\
&{\cal O}_{S}&}$$
 avec, \'evidemment, toutes les propri\'et\'es fonctorielles
 \'enonc\'ees dans le {\it th\'eoreme 2}. Ainsi, le corollaire se
 d\'eduit de facto de ce
 th\'eor\`eme.$\,\,\,\blacksquare$\bigskip\bigskip\noindent \centerline{{\tite III. La Preuve des
th\'eor\`emes  3 et 4.}}\bigskip\smallskip\noindent
{\bf 3.0.}  Le but de ce paragraphe
est de montrer  qu'un  morphisme universellement
$n$-\'equidimensionnel est analytiquement  g\'eom\'etriquement plat
si et seulement si le faisceau $\omega^{n}_{\pi}$ donn\'e par le
th\'eor\`eme 1 est enti\`erement caract\'eris\'e par la
propri\'et\'e de la trace relative et stable par changement de base.
 Apr\`es avoir donn\'e la preuve de ce fait, le d\'eduisant
 principalement du th\'eor\`eme {\bf 2}, nous proposons  une approche faisant
 abstraction de l'existence du faisceau $\omega^{n}_{\pi}$ et montrant
  que tout morphisme  analytiquement g\'eom\'etriquement plat
   porte naturellement un faisceau particulier contenant une double
   information.
\smallskip\noindent {\bf 3.1.
Preuve du th\'eor\`eme 3.}\par $\blacklozenge$ En supposons $\pi$
 analytiquement g\'eom\'etriquement plat, le {\it th\'eor\`eme 2} nous donne  un unique morphisme de faisceaux coh\'erents ${\cal
C}_{\pi,\goth{X}}:\Omega^{n}_{X/S}\rightarrow \omega^{n}_{\pi}$
compatible aux changements de bases entre espaces complexes
r\'eduits et v\'erifiant la propri\'et\'e de la trace relative.
Mais, comme nous l'avons d\'ej\`a  signal\'e , cela entraine
automatiquement que le faisceau $\omega^{n}_{\pi}$ (qui est de
profondeur au moins deux dans les fibres d'apr\`es le {\it lemme 3}) est stable par changements de bases et v\'erifie la
propri\'et\'e de la trace relative.\smallskip $\blacklozenge$
R\'eciproquement, soit $\pi\in {\cal E}(S,n)$
 muni d'une certaine pond\'eration $\goth X$ et  tel que le faisceau
 $\omega^{n}_{\pi}$ (donn\'e par le {\it th\'eor\`eme 1}) soit stable par changement
  de base et
caract\'eris\'e par la propri\'et\'e de la trace relative. Montrons,
alors, que $\pi$ est n\'ecessairement analytiquement
g\'eom\'etriquement plat.\par\noindent Le probl\`eme est de nature
locale sur $X$ et $S$. Soit
 $x$ un point de $X$ et $\ds{\xymatrix{X\ar@/_/[rr]_{\pi}\ar[r]^{f}&Y\ar[r]^{q}&S}}$
 une factorisation locale de $\pi$ en $x$ avec, comme de coutume,  $f$ fini et surjectif, $q$ lisse de
dimension relative $n$.  Comme on l'a vu dans le {\it th\'eor\`eme
1}, l'isomorphisme $f_{*}\omega^{n}_{\pi}\simeq {\cal
H}om(f_{*}{\cal O}_{X}, \Omega^{n}_{Y/S})$ nous donne une fl\`eche
$f_{*}\omega^{n}_{\pi}\rightarrow \Omega^{n}_{Y/S}$ qui est loin
d'\^etre un morphisme trace, en g\'en\'eral. Mais par hypoth\`ese,
 ce morphisme est une v\'eritable trace que l'on va noter ${\goth
T}^{n}_{f,\goth X}$, compatible \`a l'additivit\'e des
pond\'erations par essence. De plus, le morphisme naturel ${\cal
O}_{X}\rightarrow{\cal H}om(\Omega^{n}_{X/S}, \omega^{n}_{\pi})$,
donne, en appliquant le foncteur $f_{*}$ et de la trace ${\goth
T}^{n}_{f,\goth X}$, la fl\`eche $f_{*}{\cal O}_{X}\rightarrow{\cal
H}om(f_{*}\Omega^{n}_{X/S},\Omega^{n}_{Y/S})$ et, par suite, le
morphisme $f_{*}{\cal O}_{X}\rightarrow {\cal O}_{Y}$.\par\noindent
En vertu de la {\it proposition 4} (p.21)  ou du {\it corollaire 3} (p.35),  il nous
faut essentiellement v\'erifier que ${\goth T}^{n}_{f,\goth X}$
prolonge naturellement la trace usuelle
$f_{*}f^{*}\Omega^{n}_{Y/S}\rightarrow \Omega^{n}_{Y/S}$ ou, ce qui
revient au m\^eme, de montrer que le morphisme pr\'ec\'edemment
d\'ecrit par la fl\`eche ${f_{*}}{\cal O}_{X}\rightarrow {\cal
O}_{Y}$ est un morphisme trace d'ordre $0$ ( que l'on note ${\cal
T}^{0}_{f,\goth{X}}$).\par\noindent
 Pour ce faire, commen\c cons par compl\'eter la factorisation locale
 pr\'ec\'edente en l'ins\'erant dans le diagramme commutatif d\'esormais
 classique (cf $(\spadesuit)$ p.38)
$$\xymatrix{X\ar[rdd]_{\pi}\ar[rr]^{\sigma}\ar[rd]^{f}&&Z\ar[ld]_{q'}\ar[ldd]^{q}\\
&Y\ar[d]^{q''}&\\
&S&}$$ Alors la donn\'ee de ${\goth T}^{n}_{f,\goth{X}}$ est
strictement \'equivalente \`a la donn\'ee d'un morphisme
d'int\'egration ou r\'esidu
$${\rm Res}^{Z}_{Y}:{\rm I}\!{\rm R}q'_{!}{\cal E}xt^{p}(\sigma_{*}{\cal O}_{X}, \Omega^{n+p}_{Z/S})\rightarrow \Omega^{n}_{Y/S}$$
${\cal O}_{Y}$-lin\'eaire et compatible aux changements de bases
r\'eduits. Le terme ``int\'egration'' est proprement choisi (cf \S {\bf 2}, p.52  de [KI]). En effet, comme $X$ est fini sur $Y$ (donc $Y$- propre ), l'{\it
int\'egration} de classes de cohomologie de ${\rm H}^{p}_{X}(Z, \Omega^{n+p}_{Z/S})$ sur les
fibres $p$-\'equidimensionnelles du morphisme lisse $q':Z\rightarrow
Y$ donne le morphisme d'int\'egration  ${\rm H}^{p}_{X}(Z,\Omega^{n+p}_{Z/S}) \rightarrow \Gamma(Y, \Omega^{n}_{Y/S}$ que l'on peut ``faisceautiser'', gr\^ace \`a la variante relative du th\'eor\`eme de Siu-Trautmann (cf {\it lemme 1}, p.6),  pour obtenir  le  diagramme commutatif
$$\xymatrix{f_{*}{\cal E}xt^{p}(\sigma_{*}{\cal O}_{X}, \Omega^{n+p}_{Z/S})\ar[rd]\ar[rr]&&f_{*}{\cal H}^{p}_{X}(\Omega^{n+p}_{Z/S})\ar[ld]\\
&\Omega^{n}_{Y/S}&}$$
Il n'est pas difficile de voir que ${\goth
T}^{n}_{f,\goth{X}}$ contient d\'ej\`a l'information port\'ee par   ${\cal T}^{0}_{f,\goth{X}}$ et donc  ${\cal T}^{n}_{f,\goth{X}}$. En effet, soit  $x$ un point  g\'en\'erique sur $X$ et donnons nous $p$
fonctions holomorphes $(g_{1},\cdots,g_{p})$ sur $Z$ formant  une
suite r\'eguli\`ere sur $X$ en $x$. Consid\'erons  une  section arbitraire  $w$  de $\Omega^{n}_{Y/S}$. Alors,  pour toute fonction
holomorphe $h$ sur $Z$, on a l'\'egalit\'e
$${\goth
T}^{n}_{f,\goth{X}}(\left[\eqalign{hq'^{*}(w)dg_{1}\wedge\cdots\wedge
dg_{p}\cr g_{1}\cdots
g_{p}\,\,\,\,\,\,\,\,\,\,\,\,\,\,}\right]):={\rm
Res}^{Z}_{Y}\left[\eqalign{hq'^{*}(w)dg_{1}\wedge\cdots\wedge
dg_{p}\cr g_{1}\cdots g_{p}\,\,\,\,\,\,\,\,\,\,\,\,\,\,}\right] =
{\goth T}^{0}_{f,\goth{X}}(h|_{X}).w$$ de laquelle r\'esulte, gr\^ace \`a
l'abscence de torsion, le morphisme trace recherch\'e  ${f_{*}}{\cal
O}_{X}\rightarrow {\cal O}_{Y}$..
\par\noindent
 Par ailleurs,  $\omega^{n}_{\pi}$ \'etant caract\'eris\'e par  la propri\'et\'e de la trace relative, on en d\'eduit  un morphisme
$f^{*}\Omega^{n}_{Y/S}\rightarrow \omega^{n}_{\pi}$ rendant
commutatif le diagramme
$$\xymatrix{f_{*}f^{*}\Omega^{n}_{Y/S}\ar[rd]\ar[rr]&&f_{*}\omega^{n}_{\pi}\ar[ld]\\&\Omega^{n}_{Y/S}&}$$
\par\noindent Il en r\'esulte, en particulier, que le morphisme fini $f$ est
analytiquement g\'eom\'etriquement plat. Alors, utilisant les lemmes
{\bf (5.1)} ou {\bf (5.2)} de [KI] (p.86) permettant de  faire varier la
projection de sorte \`a engendrer le germe $\Omega^{n}_{X/S,x}$, on
en d\'eduit par des techniques analogues \`a celles d\'evelopp\'ees
dans la {\it proposition 4},  le diagramme commutatif
$$\xymatrix{f_{*}\Omega^{n}_{X/S}\ar[rdd]_{{\cal T}^{n}_{f, \goth X}}
\ar[rr]&&f_{*}\omega^{n}_{\pi}\ar[ldd]^{{\goth T}^{n}_{f, \goth X}}\\&f_{*}f^{*}\Omega^{n}_{Y/S}\ar[lu]\ar[ru]\ar[d]&\\
&\Omega^{n}_{Y/S}&}$$ prolongeant naturellement le pr\'ec\'edent.
Ceci \'etant valable en chaque point de $X$, la {\it proposition 4}
permet de conclure\footnote{$^{(14)}$}{On peut constater que la famille de morphismes
$f_{*}\Omega^{n}_{X/S}\rightarrow f_{*}\omega^{n}_{\pi}$, obtenue en
faisant varier le morphisme finie $f$,  d\'ecrit en fait un morphisme
global $\Omega^{n}_{X/S}\rightarrow \omega^{n}_{\pi}$ v\'erifiant
toutes les propri\'et\'es \'enonc\'ees dans le {\it th\'eor\`eme 2}.} puisque la notion de platitude g\'eom\'etrique analytique est une notion locale sur $X$. $\,\,\,\blacksquare$
\smallskip\par\noindent
{\bf 3.2. La preuve du corollaire 3.1.}\smallskip\noindent {\bf
3.2.0.} Rappelons (cf [KI], {\bf (3.3.1.2)},{\it corollaire 3}, p.73 ) que dans le cas d'un morphisme propre et plat, la theorie de la dualit\'e analytique relative de Ramis-Ruget-Verdier [R.R.V] et la construction de Kersken [Ke] permettent de voir, en particulier,  que le faisceau des formes r\'eguli\`eres  relatives est dualisant et que l'on a un isomorphisme canonique  $\omega^{n}_{X/S}\simeq{\cal H}^{-n}(\pi^{!}({\cal O}_{S}))$. De plus, si  $X$ est r\'eduit, le faisceau des formes m\'eromorphes r\'eguli\`eres $\tilde{\omega}^{n}_{X/S}$ coincide avec  le faisceau des formes r\'eguli\`eres.
 \smallskip\noindent {\bf 3.2.1. Preuve du corollaire 3.1.}\par\noindent Avant tout, disons que $X_{0}$ (resp.
$S_{0}$) est un ouvert dense de $X$ (resp. $S$) d'apr\`es le
th\'eor\`eme de Frisch ([Fr]). Apr\`es cela, il suffit d'appliquer
les th\'eor\`emes {\bf 1} et {\bf 3}. En effet, le premier  nous dit
que le faisceau $\omega^{n}_{\pi}$ est compatible aux restrictions
ouvertes sur $X$ et de ce fait
$\omega^{n}_{\pi}|_{X_{0}}=\omega^{n}_{\pi_{0}}$. Mais, $\pi_{0}$
\'etant plat, \`a fortiori analytiquement g\'eom\'etriquement plat,
le {\it th\'eor\`eme 3} dit, principalement, que
$\omega^{n}_{\pi_{0}}$ est enti\`erement caract\'eris\'e par la
propri\'et\'e de la trace relative. Alors, en tout $x$ de $X_{0}$
(qui est suppos\'e r\'eduit!), le germe en $x$ de $\pi$ est
naturellement associ\'e \`a un morphisme d'alg\`ebres analytiques
locales auquel on applique [Ke]. Il en r\'esulte, alors, par
unicit\'e d\^ue \`a la propri\'et\'e de la trace relative, une
identification canonique entre $\omega^{n}_{\pi_{0}}$ et le faisceau
de Kunz-Waldi-Kersken $\tilde{\omega}^{n}_{X_{0}/S_{0}}$.  On en
d\'eduit, alors,  que $\omega^{n}_{\pi}$ est l'unique prolongement
coh\'erent d\'esir\'e et qu'il m\'erite \`a juste titre le nom de
faisceau des formes m\'eromorphes r\'eguli\`eres.\par\noindent Les
autres propri\'et\'es d\'ecoulent \'evidemment du {\it th\'eor\`eme
3}.$\,\,\blacksquare$
\smallskip\par\noindent
 {\bf 3.3. Morphismes analytiquement g\'eom\'etriquement plats,
faisceaux dualisants et formes m\'eromorphes r\'eguli\`eres: le
th\'eor\`eme 4.}\smallskip\noindent {\bf 3.3.0.} En faisant
abstraction de l'existence de $\omega^{n}_{\pi}$ donn\'e par le {\it
th\'eor\`eme 1}, nous allons voir, dans ce qui suit,
 que la donn\'ee d'un morphisme analytiquement g\'eom\'etriquement
plat entraine l'existence d'un faisceau dualisant relatif naturel. Pour cela,
 on utilise une caract\'erisation obtenue dans le cas absolu dans [B2].
  On a, alors, le
\smallskip\noindent
\Th{4}{} Soit $\pi:X\rightarrow S$ un morphisme analytiquement
g\'eom\'etriquement plat. Alors, il existe un unique faisceau ${\cal
O}_{X}$ coh\'erent, $\Lambda^{n}_{X/S}$ v\'erifiant:\par\noindent
{\bf(i)} il est  de profondeur au moins deux fibre par fibre sur
$S$, caract\'eris\'e par la propri\'et\'e de la trace relative et
stable par changement de base entre espaces complexes r\'eduits de
dimension finie,\par\noindent {\bf(ii)} il munit $\pi$ d'un
morphisme d'int\'egration $\int_{\pi}:{\rm I}\!{\rm
R}\pi_{!}\Lambda^{n}_{X/S}\rightarrow {\cal O}_{S}$ compatible avec
l'additivit\'e des cycles et stable par changement de
base.\rm\smallskip\noindent dont la version locale est donn\'ee par
la  \Prop{7}{} Soit
$\ds{\xymatrix{X\ar@/_/[rr]_{\pi}\ar[r]^{\sigma}&Z\ar[r]^{q}&S}}$ un
diagramme commutatif d' espaces analytiques complexes r\'eduits o\`u
$\sigma$ est un plongement, $q$ un morphisme lisse de dimension
relative $n+p$ et $\pi$ analytiquement  g\'eom\'etriquement plat \`a
fibres de dimension pure $n$ muni d'une certaine pond\'eration
$\goth{X}$. Soit $\Sigma$, l'ensemble des points non r\'eguliers de
$\pi$, $j:X_{0}:= X-\Sigma \rightarrow X$ l'inclusion naturelle et
d\'esignons par
$$\tilde{\partial}^{n}_{\sigma}:
j_{*}j^{*}(\Omega^{n}_{X/S})\rightarrow {\cal H}^{1}_{\Sigma}({\cal
E}xt^{p}({\cal O}_{X}, \Omega^{n+p}_{Z/S}))$$
 le morphisme
induit par le morphisme classe fondamentale relative de la {\it
proposition 6}. Alors, $\ds{\Lambda^{n}_{\sigma}:= {\rm
Ker}\,\tilde{\partial}^{n}_{\sigma}}$  est un faisceau ${\cal
O}_{X}$-coh\'erent v\'erifiant les points {\bf(i)} et {\bf(ii)} du
{\it th\'eor\`eme 4}. \rm\smallskip\noindent \dem Soit $\ds{{\cal
C}^{\sigma}_{\pi, \goth{X}}:\Omega^{n}_{X/S}\rightarrow {\cal
E}xt^{p}({\cal O}_{X}, \Omega^{n+p}_{Z/S}),}$ le morphisme classe
fondamentale relative pond\'er\'ee de la {\it proposition.6}
relevant naturellement le morphisme cup-produit donn\'e par la
variante pond\'er\'ee (cf {\it th\'eor\`eme 0}, p.39) de la classe de
[B4]. Soit
$$\xymatrix{j_{*}j^{*}(\Omega^{n}_{X/S})
\ar[rr]^{{\partial}^{n}_{\sigma}}
\ar[rd]_{\tilde{\partial}^{n}_{\sigma}}&&{\cal H}^{1}_{\Sigma}
(\Omega^{n}_{X/S})\ar[ld]^{{\cal H}^{1}_{\Sigma}({\cal C}^{\sigma}_{\pi, \goth{X}})}\\
&{\cal H}^{1}_{\Sigma}({\cal E}xt^{p}({\cal O}_{X},
\Omega^{n+p}_{Z/S}))&}$$ le diagramme commutatif issu de suites
exactes courtes classiques  et posons
$\Lambda^{n}_{\sigma}:= {\rm Ker}\,\tilde{\partial}^{n}_{\sigma}$. Alors, en exploitant le  diagramme commutatif
$$\xymatrix{&&&0\ar[d]\\&&\Omega^{n}_{X/S}\ar[r]^{{\cal C}^{\sigma}_{\pi, \goth{X}}\,\,\,\,\,\,\,\,\,\,}\ar[d]&{\cal
E}xt^{p}({\cal O}_{X},
\Omega^{n+p}_{Z/S}))\ar[d]\\
0\ar[r]&\Lambda^{n}_{\sigma}\ar[r]&j_{*}j^{*}
\Omega^{n}_{X/S}\ar[rd]_{\!\!\!\!\!\!\tilde{\partial}^{n}_{\sigma}\,\,\,\,\,\,\,\,\,\,\,\,}\ar[r]&j_{*}j^{*}{\cal E}xt^{p}({\cal O}_{X},
\Omega^{n+p}_{Z/S})\ar[d]\\
&&&{\cal
H}^{1}_{\Sigma}({\cal E}xt^{p}({\cal O}_{X},
\Omega^{n+p}_{Z/S}))\ar[d]\\
&&&0}$$ et la  d\'efinition de $\Lambda^{n}_{\sigma}$,  il apparait clairement que
 $$\Lambda^{n}_{\sigma}\simeq {\cal E}xt^{p}({\cal
O}_{X}, \Omega^{n+p}_{Z/S})$$ Il r\'esulte, alors, du {\it lemme 3} (p.8)
 que ce faisceau  ${\cal O}_{X}$-coh\'erent est  de profondeur au
moins deux fibre par fibre sur $S$ et de la {\it proposition 6}
qu'il est stable par changement de base et v\'erifie la
propri\'et\'e de la trace relative. La fl\`eche donn\'ee dans
{\bf(ii)} est apparait d\'ej\`ale dans la preuve du point {\bf(v)}
du {\it th\'eor\`eme 1}. Mais il est plus judicieux de remarquer que
la {\it proposition 2} (p.17) peut \^etre appliqu\'ee en rempla\c cant
$\Omega^{n}_{/S}$ par $\Lambda^{n}_{\sigma}$.
$\,\,\,\blacksquare$\rm\smallskip\noindent {\bf 3.3.1. Preuve du
 th\'eor\`eme {\bf 4}.}\smallskip $\blacklozenge$ Comme $\pi$ est
universellement $n$-\'equidimensionnel, il admet des recouvrements
ouverts  $(X_{\alpha})_{\alpha\in A}$  munis  d'installations
locales (cf {\it th\'eor\`eme 2}, p.39) du type $(\spadesuit)_{\alpha\in A}$
$$\xymatrix{X_{\alpha}\ar[rdd]_{\pi_{\alpha}}\ar[rr]^{\sigma_{\alpha}}\ar[rd]^{f_{\alpha}}&&Z_{\alpha}\ar[ld]_{q'_{\alpha}}\ar[ldd]^{q_{\alpha}}\\
&Y_{\alpha}\ar[d]^{q''_{\alpha}}&\\&S&}$$ o\`u $\sigma_{\alpha}$ est
un plongement local, $\pi_{\alpha}$ la restriction de $\pi$ \`a
$X_{\alpha}$, $Y_{\alpha}$ et $Z_{\alpha}$ des espaces complexes
lisses sur $S$ de dimension relative $n+p$ et $n$ respectivement,
$f_{\alpha}$ fini, ouvert  et surjectif, $q_{\alpha}$,$q'_{\alpha}$
et $q''_{\alpha}$ lisses. \par\noindent $\pi$ \'etant analytiquement
g\'eom\'etriquement plat, il en est de m\^eme de  $\pi_{\alpha}$ et
$f_{\alpha}$.\par\noindent La {\it proposition 6} (p.40)  nous donne, pour
chaque $\alpha$, un morphisme classe  fondamentale relative
pond\'er\'e (on omet d\'elib\'eremment l'indexation en
$\sigma_{\alpha}$)
$${\cal C}_{\pi_{\alpha}, \goth{X}_{\alpha}}:\Omega^{n}_{X_{\alpha}/S}\rightarrow {\cal
E}xt^{p}({\cal O}_{X_{\alpha}}, \Omega^{n+p}_{Z_{\alpha}/S})$$
induisant  le diagramme
commutatif
$$\xymatrix{{j_{\alpha}}_{*}{j_{\alpha}}^{*}(\Omega^{n}_{X_{\alpha}/S})
\ar[rr]^{\partial^{n}_{\alpha}}
\ar[rd]_{\tilde{\partial}^{n}_{\alpha}}&&{\cal H}^{1}_{\Sigma_{\alpha}}
(\Omega^{n}_{X_{\alpha}/S})\ar[ld]^{{\cal H}^{1}_{\Sigma}({\cal C}_{\pi_{\alpha}, \goth{X}_{\alpha}})}\\
&{\cal H}^{1}_{\Sigma_{\alpha}}({\cal
E}xt^{p}({\sigma_{\alpha}}_{*}{\cal O}_{X_{\alpha}},
\Omega^{n+p}_{Z_{\alpha}/S}))&}$$ La  {\it proposition 7} (p.58) nous dit que
$\ds{\Lambda^{n}_{\alpha}:= {\rm
Ker}\,{\tilde{\partial}^{n}}_{\alpha}}$ est un faisceau ${\cal
O}_{X_{\alpha}}$-coh\'erent dot\'e des propri\'et\'es qui y sont
\'enum\'er\'ees. On obtient, ainsi, une donn\'ee
 simplicielle $(\Lambda^{n}_{\alpha})_{\alpha\in A}$ dont il nous faut montrer qu'elle
 constitue une donn\'ee de recollement d\'efinissant, ainsi,  un faisceau intrins\`eque sur $X$. Pour cela, on v\'erifie son ind\'ependance, vis-\`a-vis du plongement (ou de l'installation
locale) choisi. Comme la justification utilise, presqu'en tout point, des arguments semblables \`a ceux d\'ej\`a rencontr\'es dans
la preuve du {\it th\'eor\`eme 1}, nous omettrons certains d\'etails.\par\noindent Soient donn\'ees deux installations locales de $\pi$  munies de deux plongements $\sigma$ et $\rho$.
L\`a encore, quelques r\'eductions
standards, nous conduisent  \`a la situation classique donn\'ee par
le diagramme commutatif de $S$-espaces complexes
$$\xymatrix{X\ar[r]_{\sigma}\ar@/^/[rr]^{\rho}\ar[rd]_{\pi}&Z\ar[r]^{\!\!\!\!\!\!\!\!h}\ar[d]_{q}&Z\times V\ar[ld]^{q'}\ar@/^/[l]^{g}\\
&S&}$$
dans lequel  $V$ est  un ouvert d'un certain espace num\'erique ${\Bbb C}^{N}$ , $h$ un plongement muni d'une section $g$ (i.e $hog=Id$)
faisant de $Z$ un r\'etract de $Z_{1}:= Z\times V$, $q$ et $q'$ des morphismes lisses de dimension relatives respectives  $n+p$ et $n+p+p'$.\par\noindent
Remarquons que si $\pi$ est muni de la pond\'eration  $\goth{X}$ vue comme cycle de $S\times Z$, il peut aussi \^etre consid\'er\'e comme pond\'er\'e par le cycle ${\goth{X}}_{1}:=(h\times Id_{S})_{*}\goth{X}$ de $S\times Z_{1}$. Cela traduit simplement le fait que le morphisme classe fondamentale relative pond\'er\'ee est intrins\`eque sur $X$ comme l'annonce le {\it th\'eor\`eme 2}.
\par\noindent
Alors les morphismes classes fondamentales relatives de $X$ dans $Z$ et $Z$ dans $Z_{1}$ donnent le diagramme commutatif
$$\xymatrix{\Omega^{n}_{X/S}\ar[r]^{{\cal
C}^{\sigma}_{\pi,\goth{X}}}
 \ar[d]_{\phi}&{\cal E}xt^{p}({\cal O}_{X}, \Omega^{n+p}_{Z/S})
 \ar[d]\ar[ld]_{\beta}
 \\{\cal E}xt^{p+p'}
 ({\cal O}_{X},
  \Omega^{n+p+p'}_{Z_{1}/S})&{\cal E}xt^{p}({\cal O}_{X}, {\cal E}xt^{p'}
  ({\cal O}_{Z},
  \Omega^{n+p+p'}_{Z_{1}/S}))\ar[l]&},$$
dans lequel la fl\`eche horizontale inf\'erieure est un isomorphisme
 puisque
$${\cal E}xt^{j}({\cal O}_{Z}, \Omega^{n+p+p'}_{Z_{1}/S})=0,\,\,\,\forall\,j\not= q$$
D\'esignons par  $\Lambda^{n}_{\sigma}$ et $\Lambda^{n}_{\rho}$ les faisceaux coh\'erents  d\'efinis  par la
{\it proposition 7}. Le diagramme ci-dessus induit, en utilisant les suites exactes courtes classiques,   le diagramme commutatif
$$\xymatrix{&{\Lambda^{n}_{\rho}}\ar[d]&\\
{\Lambda^{n}_{\sigma}}\ar[r]&j_{*}j^{*}(\Omega^{n}_{X/S})
\ar[d]_{\tilde{\partial}^{n}_{\rho}}\ar[r]^{\tilde{\partial}^{n}_{\sigma}}&{\cal
H}^{1}_{\Sigma} ({\cal E}xt^{p}({\cal O}_{X},
\Omega^{n+p}_{Z/S}))
\ar[ld]^{{\cal H}^{1}_{\Sigma}(\beta)}\\
&{\cal H}^{1}_{\Sigma}({\cal E}xt^{p+p'}({\cal O}_{X},
\Omega^{n+p+p'}_{Z_{1}/S}))\ar@<6pt>[ru]^{{\cal
H}^{1}_{\Sigma}(Res^{Z_{1}}_{Z})}&}$$
Dont il  est facile de  d\'eduire, pour des raisons d\'ej\`a invoqu\'ees dans la preuve du {\it th\'eor\`eme 1},  l'isomorphisme $\ds{{\Lambda^{n}_{\sigma}}\simeq \Lambda^{n}_{\rho}}$ traduisant  le recollement de ces donn\'ees locales en un faisceau coh\'erent
  global $\Lambda^{n}_{X/S}$. Le fait qu'il soit  de formation compatible aux changements de bases entre espaces complexes r\'eduits et qu'il v\'erifie la propri\'et\'e de la trace relative d\'ecoule naturellement de la  {\it  proposition 7}.\smallskip\noindent
Il nous reste \`a voir que cette derni\`ere propri\'et\'e  caract\'erise
compl\`etement ce faisceau. Pour cela, consid\'erons une section
$\xi\in \Gamma(V,j_{*}j^{*}\Omega^{n}_{X/S})$ satisfaisant la
propri\'et\'e de la trace relative et montrons qu'elle d\'efinit
n\'ecessairement une (unique) section du faisceau
$\Lambda^{n}_{X/S}$.\par\noindent Comme le  probl\`eme \'etant  de nature
locale sur $X$, on peut supposer donn\'ee  une factorisation locale
$\ds{\xymatrix{X\ar@/_/[rr]_{\pi}\ar[r]^{f}&Y\ar[r]^{q}&S}}$ de
$\pi$  avec $f$ fini et g\'eom\'etriquement plat  sur $Y$ lisse sur
$S$ et de dimension relative $n$. Par ailleurs,  la notion de platitude
g\'eom\'etrique, et de ce fait, la construction du faisceau
$\Lambda^{n}_{X/S}$, \'etant compatible  aux changements de bases
``r\'eduits'', on peut supposer $S$ irr\'eductible et raisonner sur
la composante g\'en\'erique du morphisme classifiant (cf [KI], {\bf
(1.0.3.2))}, p.24) de sorte que la trace relative ait un sens. Alors,
il suffit simplement d' utiliser la classe fondamentale relative et
la propri\'et\'e de la trace qui la caract\'erise (cf {\it th\'eor\`eme 2}, {\bf (ii)}).\par\noindent En effet, le morphisme classe
fondamentale relative pond\'er\'e (locale)\par
\centerline{${\cal C}^{\sigma}_{\pi, \goth{X}}:\Omega^{n}_{X/S}\rightarrow {\cal E}xt^{p}({\cal O}_{X}, \Omega^{n+p}_{Z/S})$}
\noindent se factorise au travers de $\Lambda^{n}_{\sigma}$ (et ce par
d\'efinition m\^eme de ce faisceau!) et permet d'\'ecrire
 \par\centerline{${\cal T}^{n}_{f}(\xi)\in
\Gamma(Y, \Omega^{n}_{Y/S})\Longrightarrow {\cal C}^{\sigma}_{\pi,
\goth{X}}(\xi)\in \Gamma(X,{\cal E}xt^{p}({\cal O}_{X},
\Omega^{n+p}_{Z/S}))$} \noindent Pour s'en convaincre, il suffit
simplement de reprendre les arguments utilisant la dualit\'e
analytique locale d\'evelopp\'es dans la preuve de la {\it proposition 6}. Alors, l'hypoth\`ese faite sur $\xi$ et la
propri\'et\'e de la trace relative  de la classe fondamentale
relative pond\'er\'ee (cf {\it th\'eor\`eme 2}) montre clairement que ``
la forme '' $\xi\wedge {\cal C}^{\sigma}_{\pi,\goth{X}}$ se prolonge
en une unique section globale du faisceau ${\cal E}xt^{p}({\cal
O}_{X}, \Omega^{n+p}_{Z/S})$. Ceci \'etant valable pour toute
installation locale, on en d\'eduit, gr\^ace \` a {\it la proposition 7}, que
 $\xi$ d\'efinit une section du faisceau $\Lambda^{n}_{X/S}$. De
 plus, ce dernier \'etant de profondeur au moins deux fibre par fibre et
donc sans torsion fibre par fibre, la section ainsi d\'efinie est
unique.\par\noindent
Pour \'etablir {\bf(ii)}, il suffit de se r\'ef\'erer \`a la construction (locale) donn\'ee  dans le {\it th\'eor\`eme 1}. On peut aussi s'inspirer de la {\it proposition 1} en y rempla\c cant le faisceau $\Omega^{n}_{X/S}$ par $\Lambda^{n}_{X/S}$. En effet, la nature locale du faisceau $\Lambda^{n}_{X/S}$ et  sa caract\'erisation par la propri\'et\'e de la trace relative entraine automatiquement l'existence d'une fl\`eche
$${\rm I}\!{\rm R}^{n}{\pi}_{!}\Lambda^{n}_{X/S}\rightarrow {\cal O}_{S}$$
On s'en convainc ais\'ement en invoquant  le th\'eor\`eme de Reiffen
 permettant  de localiser sur $X$ gr\^ace \`a la suite exacte \`a droite (cf [KI], {\bf (2.0)}, p.53 ) $$ \bigoplus_{\alpha\in A}{\rm I}\!{\rm R}^{n}{\pi_{\alpha}}_{!}
\Lambda^{n}_{{X_{\alpha}}/S} \rightarrow {\rm I}\!{\rm
R}^{n}{\pi_{!}}{\Lambda^{n}_{X/S}}\rightarrow 0$$
 produisant, par recollement, la  fl\`eche d\'esir\'ee.\smallskip\noindent
Soit $\ds{\xymatrix{X_{2}\ar@/_/[rr]_{\pi_{2}}\ar[r]^{\Psi}&X_{1}\ar[r]^{\pi_{1}}&S}}$ un diagramme commutatif d'espaces analytiques complexes r\'eduits avec $\pi_{1}$ (resp. $\pi_{2}$)  analytiquement g\'eom\'etriquement plat  de dimension relative $n_{1}$ (resp. $n_{2}$) et $\Psi$ propre universellement \'equidimensionnel de dimension relative $d:=n_{2}-n_{1}$.\par\noindent
Pour montrer l'existence du morphisme ${\rm I}\!{\rm
 R}^{d}{\Psi}_{*}\Lambda^{n_{1}+d}_{{X_{2}}/S}\rightarrow
 \Lambda^{n_{1}}_{{X_{1}}/S}$, on proc\`ede exactement
  de la m\^eme fa\c con que pour le point {\bf (vii)} du {\it th\'eor\`eme 1}. Le probl\`eme \'etant de nature locale sur
    $X_{1}$ et $X_{2}$, il nous suffit de v\'erifier que
 pour tout germe de param\'etrisation locale $f:X_{1}\rightarrow S\times
 U$, on a un morphisme (bien d\'efini) $f_{*}{\rm I}\!{\rm
 R}^{d}{\Psi}_{*}\Lambda^{n_{1}+d}_{{X_{2}}/S}\rightarrow
 \Omega^{n_{1}}_{{S\times U}/S}$ et qui
 est un v\'eritable morphisme de type trace. On peut donc supposer
 $Y:=S\times U$ et $\pi_{1}$ la projection canonique $q:S\times
 U\rightarrow S$. Comme $\Lambda^{n_{j}}_{{X_{j}}/S}\simeq
 \omega^{n_{j}}_{\pi_{j}}$ (cf {\it corollaire 7}), pour $j=1,2$, l'existence de la fl\`eche
 r\'esulte de {\bf(vii)} du {\it th\'eor\`eme 1}. Pour s'assurer que
 le faisceau coh\'erent ${\rm I}\!{\rm
 R}^{d}{\Psi}_{*}\Lambda^{n_{1}+d}_{{X_{2}}/S}$ v\'erifie la
 propri\'et\'e de la trace relative, il suffit, gr\^ace \`a la
 description du morphisme d'int\'egration faite dans [KI], \S{\bf 2}, p.52, d'expliciter en terme de cochaines de \v Cech le morphisme ${\rm I}\!{\rm
 R}^{d}{\Psi}_{*}\Lambda^{n_{1}+1}_{{X_{2}}/S}\rightarrow
 \Omega^{n_{1}}_{{S\times U}/S}$. En effet, soit $t\in S\times U$ et $({\goth U})_{\alpha\in A}$ un recouvrement ouvert de $X_{2}$ adapt\'e \`a la fibre $\Psi^{-1}(t)$. Soit $\xi\in {\rm H}^{d}({\goth U},\Lambda^{n_{1}+d}_{{X_{2}}/S})$ repr\'esent\'e par le cocyle que l'on note encore $\xi:=(\xi_{\alpha_{0}\cdots\alpha_{q}})\in {\cal Z}^{q}({\goth U},
{\cal F})$ dont un repr\'esentant $\bar\partial$-ferm\'e est donn\'e, au signe pr\`es,  par
 $$\psi^{n_{1}+d,n_{1}}:=\ds{\sum_{\alpha_{0}\cdots\alpha_{n_{1}+d}}}
\rho_{\alpha_{0}}{\bar\partial}\rho_{\alpha_{1}}\wedge\cdots
\wedge{\bar\partial}\rho_{\alpha_{n_{1}+d}}\wedge
\xi_{\alpha_{0}\cdots\alpha_{n_{1}}\lambda_{0}\cdots\lambda_{d-1}}$$
et que l'on int\`egrera judicieusement relativement \`a $S$. Nous \'evitons les d\'etails techniques dont la lourdeur, inh\'erente  au calcul ``\v Cechiste'', n'apporte que tr\`es peu de substance \`a  la compr\'ehension du texte.\par\noindent
Comme $\pi_{2}$ est analytiquement g\'eom\'etriquement plat, pour tout $s\in S$,  on peut choisir un recouvrement ouvert  $(\tilde{\goth U})_{\beta\in B}$ adapt\'e \`a $\pi_{2}^{-1}(s)$ et tel que ${\goth U}=\tilde{\goth U}\cap \Psi^{-1}(t)$. Alors, utilisant la caract\'erisation du faisceau
$\Lambda^{n_{1}+d}_{{X_{2}}/S}$ relativement au recouvrement $\tilde{\goth U}$, on se ram\`ene essentiellement \`a la situation simple  donn\'ee par $f:Z\rightarrow Y$  un morphisme lisse entre vari\'et\'es complexes $S$-lisses pour laquelle on v\'erifie l'existence d'une fl\`eche d'int\'egration ${\rm I}\!{\rm
 R}^{d}{f}_{!}\Omega^{n_{1}+d}_{Z/S}\rightarrow\Omega^{n_{1}}_{Y/S} $; ce dont il est facile de se convaincre en la construisant explicitement sans aucune difficult\'e. En fait, quitte \`a r\'etr\'ecir les donn\'ees, on installe la factorisation pr\'ec\'edente dans le diagramme commutatif
$$\xymatrix{&X_{2}\ar[dd]^{\pi_{2}}\ar[ld]_{F}\ar[rd]^{\Psi}&\\
Z:=S\times U\times V\ar[rd]_{g}\ar[rr]_{h}\ar[rd]&&Y:=S\times U\ar[ld]^{\pi_{1}}\\
&S&}$$
avec $U$ (resp. $V$) un ouvert relativement compact de ${\Bbb C}^{n_{1}}$ (resp.  ${\Bbb C}^{d}$), $F$ analytiquement $0$-g\'eom\'etriquement plat, $h$ et $g$ les projections canoniques. Alors, le morphisme trace $F_{*}\Lambda^{n_{1}+d}_{{X_{2}}/S}\rightarrow\Omega^{n_{1}+d}_{Z/S}$ et l'int\'egration ${\rm I}\!{\rm
 R}^{d}{h}_{!}\Omega^{n_{1}+d}_{Z/S}\rightarrow\Omega^{n_{1}}_{Y/S}$  induiront, gr\^ace \`a l'isomorphisme ${\rm I}\!{\rm  R}^{d}{h}_{!}F_{*}\Lambda^{n_{1}+d}_{{X_{2}}/S}\simeq {\rm I}\!{\rm  R}^{d}\Psi_{*}\Lambda^{n_{1}+d}_{{X_{2}}/S}$, le morphisme attendu  ${\rm I}\!{\rm  R}^{d}\Psi_{*}\Lambda^{n_{1}+d}_{{X_{2}}/S}\rightarrow \Omega^{n_{1}}_{Y/S}$ $\,\,\,\blacksquare$\bigskip\noindent
 \Cor{7}{} Avec les
hypoth\`eses du {\it th\'eor\`eme 4}, on a des isomorphismes
canoniques
$$\Lambda^{n}_{X/S}\simeq \tilde{\omega}^{n}_{X/S}\simeq \omega^{n}_{\pi}$$
\rm\dem
L'isomorphisme canonique
$\omega^{n}_{\pi}\simeq \tilde{\omega}^{n}_{X/S}$ faisait l'objet du
{\it corollaire (3.1)}. Il nous reste simplement \`a v\'erifier
l'existence d'un isomorphisme canonique entre les faisceaux
$\omega^{n}_{\pi}$ et $\Lambda^{n}_{X/S}$.\par\noindent Or  (avec
les notations du {\it th\'eor\`eme 4} et de la {\it proposition 7}), on avait constat\'e que $\ds{\Lambda^{n}_{\alpha}:= {\rm
Ker}\,{\tilde{\partial}^{n}}_{\alpha}\simeq\omega^{n}_{\pi_{\alpha}}}$.
Mais les propri\'et\'es du faisceau $\omega^{n}_{\pi}$ donn\'e par
le {\it th\'eor\`eme 1}, et, en particulier, sa nature locale, permettent
d'\'ecrire
$\ds{\omega^{n}_{\pi_{\alpha}}=\omega^{n}_{\pi}|_{X_{\alpha}}}$.
D'o\`u, pour chaque ouvert d'un recouvrement ouvert localement fini
adapt\'e \`a $\pi$, les isomorphismes
$\ds{\Lambda^{n}|_{X_{\alpha}}\simeq
\omega^{n}_{\pi}|_{X_{\alpha}}}$. Par ailleurs, le  morphisme
canonique
 $\ds{{\cal C}_{\pi,\goth{X}}: \Omega^{n}_{X/S}\rightarrow
 \omega^{n}_{\pi}}$
 donn\'e dans le {\it th\'eor\`eme 2}  induit le diagramme commutatif
$$\xymatrix{&&\Omega^{n}_{X/S}\ar[r]^{{\cal C}_{\pi, \goth X}}\ar[d]&\omega^{n}_{\pi}\ar[d]\\
0\ar[r]&\Lambda^{n}_{X/S}\ar[r]&j_{*}j^{*}\Omega^{n}_{X/S}\ar[rd]_{{\tilde\partial}^{n}} \ar[r]&j_{*}j^{*}\omega^{n}_{\pi}\ar[d]\\
&&&{\cal H}^{1}_{\Sigma}(\omega^{n}_{\pi})\ar[d]\\
&&&0}$$
 mettant clairement  en \'evidence une fl\`eche naturelle $\ds{\Lambda^{n}_{X/S}\rightarrow \omega^{n}_{\pi}}$ injective puisque ces faisceaux coincident g\'en\'eriquement et sont sans ${\cal O}_{X}$-torsion.\par\noindent
 En effet, le {\it th\'eor\`eme 4} montre  que ${\cal C}_{\pi, \goth{X}}$
est en fait \`a valeurs dans le faisceau $\Lambda^{n}_{X/S}$ qui est
caract\'eris\'e par la propri\'et\'e de la trace relative. On en
d\'eduit, donc, par universalit\'e, un morphisme canonique
$\Lambda^{n}_{X/S}\rightarrow \omega^{n}_{\pi}$. De l'isomorphisme local
et de cette fl\`eche injective d\'ecoule l'isomorphisme global et
canonique (d\'efini \`a isomorphisme canonique pr\`es!)
$$\Lambda^{n}_{X/S}\simeq \omega^{n}_{\pi}\,\,\,\,\blacksquare$$\smallskip\noindent
\noindent{\bf 3.3.2. $k$-formes m\'eromorphes r\'eguli\`eres.}\smallskip\noindent
Pour obtenir des formes m\'eromorphes r\'eguli\`eres en tout degr\'e $k$, on utilise essentiellement la {\it proposition 6}.  On a, alors,
 \Prop{7bis}{} Avec les hypoth\`eses et notations de la {\it proposition 7}, soit
$$\tilde{\partial}^{k}_{\sigma}:
j_{*}j^{*}(\Omega^{k}_{X/S})\rightarrow {\cal H}^{1}_{\Sigma}({\cal
E}xt^{p}({\cal O}_{X}, \Omega^{k+p}_{Z/S}))$$
 le morphisme
induit par le morphisme classe fondamentale relative de la {\it
proposition 6}. Alors, $\ds{\Lambda^{k}_{\sigma}:= {\rm
Ker}\,\tilde{\partial}^{k}_{\sigma}}$  est un faisceau ${\cal
O}_{X}$-coh\'erent  de profondeur au moins deux fibre par fibre sur
$S$, caract\'eris\'e par la propri\'et\'e de la trace relative et
stable par changement de base entre espaces complexes r\'eduits de
dimension finie. \rm
\dem Soit ${\cal C}^{k}_{\pi, \goth
X}:\Omega^{k}_{X/S}\rightarrow {\cal E}xt^{p}({\cal O}_{X},
\Omega^{p+k}_{Z/S})$ le morphisme induit par le morphisme classe
fondamentale relative pond\'er\'e de la {\it proposition 6}. On
obtient, alors, un diagramme analogue \`a celui de la page 57
$$\xymatrix{&&&0\ar[d]\\&&\Omega^{k}_{X/S}\ar[r]^{{\cal C}^{k}_{\pi, \goth{X}}\,\,\,\,\,\,\,\,\,\,}\ar[d]&{\cal
E}xt^{p}({\cal O}_{X},
\Omega^{n+k}_{Z/S}))\ar[d]\\
0\ar[r]&\Lambda^{k}_{\sigma}\ar[r]&j_{*}j^{*}
\Omega^{k}_{X/S}\ar[rd]_{\!\!\!\!\!\!\tilde{\partial}^{k}_{\sigma}\,\,\,\,\,\,\,\,\,\,\,\,}\ar[r]&j_{*}j^{*}{\cal E}xt^{p}({\cal O}_{X},
\Omega^{n+k}_{Z/S})\ar[d]\\
&&&{\cal
H}^{1}_{\Sigma}({\cal E}xt^{p}({\cal O}_{X},
\Omega^{n+k}_{Z/S}))\ar[d]\\
&&&0}$$ dans lequel on  pose  $ \Lambda^{k}_{\sigma}:= {\rm
Ker}\,\tilde{\partial}^{k}_{\sigma}$. et dont on va montrer les
propri\'et\'es annonc\'ees.\smallskip Pour cela, commen\c cons par
remarquer que le diagramme ci-dessus nous donne une fl\`eche
$\phi:\Lambda^{k}_{\sigma}\rightarrow{\cal E}xt^{p}({\cal O}_{X},
\Omega^{k+p}_{Z/S})$ qui est injective puisque les faisceaux
$\Lambda^{k}_{\sigma}$ et ${\cal E}xt^{p}({\cal O}_{X},
\Omega^{k+p}_{Z/S})$ sont g\'en\'eriquement isomorphes et  sans
${\cal O}_{X}$-torsion.\par\noindent
 Consid\'erons un syst\`eme $(f_{1},\cdots,f_{p})$ de fonctions
analytiques sur $Z$ s'annulant sur $X$ et donnant des \'equations
r\'eduites de $X$ dans $Z$ et le morphisme
$$\Phi: {\cal E}xt^{p}({\cal O}_{X}, \Omega^{k+p}_{Z/S})\rightarrow({\cal
E}xt^{p}({\cal O}_{X}, \Omega^{k+p+1}_{Z/S}))^{p}$$ donn\'e par
$\Phi(\alpha):=(\alpha\wedge df_{1},\cdots,\alpha\wedge df_{p})$.
Alors, $\Lambda^{k}_{\sigma}\simeq {\rm Ker}\Phi$.\par\noindent Par
d\'efinition m\^eme, on a clairement une
 fl\`eche ${\rm Ker}\Phi\rightarrow \Lambda^{k}_{\sigma}$ d'ailleurs
  injective pour les m\^emes raisons que  $\phi$. D'autre part, une section $\sigma$ du faisceau
$j_{*}j^{*}\Omega^{k}_{X/S}$ qui a une image nulle  par
$\tilde{\partial}^{k}_{\sigma}$  provient n\'ecessairement  d'une
section
   $\tilde\sigma$ du faisceau ${\cal E}xt^{p}({\cal O}_{X},
\Omega^{p+k}_{Z/S})$. Or, par construction, elle  est
g\'en\'eriquement annul\'ee par le morphisme $\Phi$ et donc
 globalement annul\'ee par ce morphisme en raison de l'abscence
 de ${\cal O}_{X}$- torsion.  On d\'efinit ainsi
une section du faisceau  ${\rm Ker}\Phi$ et donc une fl\`eche
$\Lambda^{k}_{\sigma}\rightarrow {\rm Ker}\Phi$ qui est
n\'ecessairement injective. D'o\`u l'isomorphisme
$\Lambda^{k}_{\sigma}\simeq{\rm Ker}\Phi$ et par suite la ${\cal
O}_{X}$- coh\'erence du faisceau
$\Lambda^{k}_{\sigma}$.\par\noindent Pour montrer qu'il est de
profondeur au moins deux fibre par fibre, il suffit de montrer
l'isomorphisme
$$\Lambda^{k}_{\sigma}\simeq {\cal H}om_{{\cal
O}_{X}}(\Omega^{n-k}_{X/S}, \Lambda^{n}_{\sigma})$$ Pour ce faire,
on commence par remarquer que la stabilit\'e par cup produit des
sections de $\Lambda^{k}_{\sigma}$ par les sections du faisceau
$\Omega^{k}_{X/S}$, donne imm\'ediatement une fl\`eche
$$\Lambda^{k}_{\sigma}\rightarrow
{\cal H}om_{{\cal O}_{X}}(\Omega^{n-k}_{X/S},
\Lambda^{n}_{\sigma})$$ qui est  automatiquement  injective puisque
les faisceaux coincident g\'en\'eriquement et sont sans ${\cal
O}_{X}$-torsion. La fl\`eche inverse se construit exactement comme
dans le cas absolu trait\'e dans [B2].\par\noindent
Alors, comme  $\Lambda^{n}_{\sigma}$ est de profondeur au moins deux fibre par fibre et qu'il est v\'erifie la propri\'et\'e de la trace relative, en vertu de la  {\it proposition 7}, on en d\'eduit que $\Lambda^{k}_{\sigma}$ satisfait aussi  ces propri\'et\'es. En effet, il suffit simplement, par exemple,  d'exploiter  (comme on l'a expliqu\'e plusieurs fois) l'isomorphisme de dualit\'e pour $f$ donn\'e par  $f_{*}{\cal
H}om_{{\cal O}_{X}}(\Omega^{n-k}_{X/S},\Lambda^{n}_{\sigma})= {\cal
H}om_{{\cal O}_{Y}}({f}_{*}\Omega^{n-k}_{X/S}, \Omega^{n}_{Y/S})$ et le morphisme ${f}_{*}{\cal H}om_{{\cal O}_{X}}(\Omega^{n-k}_{X/S},\Lambda^{n}_{\sigma}
)\rightarrow \Omega^{k}_{Y/S}$ (obtenu gr\^ace \`a l'image r\'eciproque naturelle $\Omega^{n-k}_{Y/S}\rightarrow
{f}_{*}\Omega^{n-k}_{X/S}$). \par\noindent
Il est facile de v\'erifier que la propri\'et\'e de la trace relative caract\'erise compl\`etement les sections du faisceaux $\Lambda^{k}_{\sigma}$. En effet, une section $\xi$ de
$j_{*}j^{*}\Omega^{k}_{X/S}$, v\'erifiant cette propri\'et\'e d\'efinit naturellement une section de $\Lambda^{k}_{\sigma}$ puisque, au vu de l'interpr\'etation pr\'ec\'edente,  l'image de $\xi$ par ${\cal C}^{k}_{\pi,\goth{X}}$ se prolonge en une section
    globale de ${\cal E}xt^{p}({\cal O}_{X}, \Omega^{k+p}_{Z/S})$. D'ailleurs, en utilisant l'isomorphisme $\Lambda^{k}_{\sigma}\simeq{\cal H}om_{{\cal
O}_{X}}(\Omega^{n-k}_{X/S}, \Lambda^{n}_{\sigma})$, on se ram\`ene \`a v\'erifier le fait, d\'ej\`a \'etabli pour  $\Lambda^{n}_{X/S}$,
$${\cal T}^{n}_{f,\goth{X}}
    (\xi)\in
\Gamma(Y, \Omega^{n}_{Y/S})\Longrightarrow {\cal C}^{\sigma}_{\pi,
\goth{X}}(\xi)\in \Gamma(X,{\cal E}xt^{p}({\cal O}_{X},
\Omega^{n+p}_{Z/S}))\,\,\,\,\,\,\,\blacksquare$$\smallskip\noindent
\Cor{8}{}
 Soit $\pi:X\rightarrow S$  un
morphisme $n$- analytiquement g\'eom\'etriquement plat muni d'une
pond\'eration $\goth{X}$. Alors, pour tout entier $k\geq 0$, il
exite un unique faisceau ${\cal O}_{X}$-coh\'erent
$\Lambda^{k}_{X/S}$ v\'erifiant :\par\noindent {\bf(i)}
$\Lambda^{k}_{X/S}$ est de profondeur au moins deux fibre par fibre
sur $S$ et coincide avec le faisceau des formes holomorphes
$S$-relatives sur la partie r\'eguli\`ere de
$\pi$.\smallskip\noindent {\bf(ii)} il est muni d'un isomorphisme
canonique $\Lambda^{k}_{X/S}\simeq {\cal H}om(\Omega^{n-k}_{X/S},
\Lambda^{n}_{X/S})$
\smallskip\noindent {\bf(iii)} Les faisceaux $\Lambda^{k}_{X/S}$
sont dot\'es d'une diff\'erentielle ${\rm D}$ et d'un morphisme
canonique de complexes diff\'erentielles
 $\Omega^{\bullet}_{X/S}\rightarrow \Lambda^{\bullet}_{X/S}$
 compatible avec la propri\'et\'e de la trace relative.
\rm\dem
 Tout comme dans le {\it th\'eor\`eme 4}, on proc\'ed\'e par recollement des donn\'ees locales collect\'ees, cette fois, gr\^ace \`a  la {\it proposition 7bis}. Ainsi, tout recouvrement ouvert  arbitraire $(X_{\alpha})_{\alpha}$ localement fini et adapt\'e \`a
   $\pi$ sur $X$, fournit  une donn\'ee simplicielle
   $(\Lambda^{k}_{\alpha})_{\alpha}$ dont on montre la propri\'et\'e
   de recollement. On aboutira ainsi \`a un faisceau intrins\`eque
   ${\cal O}_{X}$-coh\'erent
   $\Lambda^{k}_{X/S}$ satisfaisant \'evidemment
   {\bf(i)} et {\bf(ii)}.\par\noindent
De plus, on a un isomorphismes canonique
$\Lambda^{k}_{X/S}\simeq \omega^{k}_{\pi}$. En effet, le {\it corollaire (1.1)} nous donne  $\omega^{k}_{\pi}\simeq
{\cal H}om(\Omega^{k}_{X/S},\omega^{n}_{\pi} )$ et le  {\it corollaire 7}
assure  l'identification $\Lambda^{n}_{X/S}\simeq \omega^{n}_{\pi}$. D'o\`u, en vertu de {\bf(ii)}, $\Lambda^{k}_{X/S}\simeq\omega^{k}_{\pi}$.\par\noindent
{\bf(iii)} d\'ecoule, alors, de la {\it proposition 5}.$\,\,\,\blacksquare$\smallskip\smallskip\smallskip\bigskip\noindent
\centerline{\bf R\'ef\'erences bibliographiques}

\medskip
{\baselineskip=9pt\eightrm
\item{[A.E]}Ang\'eniol~B, Elzein~F., La classe fondamentale relative d'un cycle,Bull.Soc.Math.France Mem~58 (1978) 63---93.
\smallskip

\item{[A.L.J]}Ang\'eniol~B, Lejeune-Jalabert M., Calcul diff\'erentiel et classes caract\'eristiques en g\'eom\'etrie alg\'ebrique. Travaux en Cours~38, Hermann, Paris (1989).
\smallskip

\item {[A.N]}Andreotti~A., Norguet~F., La
convexit\'e holomorphe dans l'espace analytique des cycles d'une
vari\'et\'e  alg\'ebrique, Ann. Sc. Norm. Pisa, tome  21, (1965),
811---842.
\smallskip

\item {[A.K]}Andreotti~A., Kaas~A. Duality on complex spaces, Ann. Sc. Norm. Pisa, tome  27, (1973) 187---263.
\smallskip

\item{[Ba]}Banica~C., Sur les ext globaux dans une d\'eformation. Compositio.Math, ~44, Fasc.1-3,(1981), p.17---27
\smallskip

\item{[B1]}Barlet~D.,
Espace analytique r\'eduit des cycles analytiques complexes compacts
d'un espace analytique
 r\'eduit,  S\'em. F. Norguet-Lecture Notes in Mathematics.~482,
 Springer Verlag, (1975)  1---158.
\smallskip

\item{[B2]}Barlet~D.Faisceau $\omega^{.}_{X}$ sur un espace analytique de dimension pure,  S\'em. F. Norguet-Lecture Notes in Mathematics ~670, Springer Verlag, (1978) 187---204.
\smallskip

\item{[B3]}Barlet~D. Convexit\'e au voisinage d'un cycle. Sem.Fran\c cois Norguet-Lecture Notes in Mathematics ~807, Springer-Verlag, (1977-1979),p.102---121.
\smallskip

\item{[B4]}Barlet~D., Famille analytique de cycles et classe
fondamentale relative, S\'em. F. Norguet-Lecture Notes in
Mathematics ~807, Springer Verlag,  (1980)  1---24.
\smallskip

\item{[B5]}Barlet~D., Majoration du volume des fibres g\'en\'eriques et forme g\'eom\'etrique du th\'eor\`eme d'applatissement. Seminaire P. Lelong-H.Skoda. Lecture Notes in Mathematics~822, Springer verlag, (1980), p.1---17
\smallskip

\item{ [B.M]}Barlet~D.,  Magnusson~J.,
Integration de classes de cohomologie m\'eromorphes et diviseur
d'incidence, Ann. sc. de l'E.N.S: S\'erie 4, tome 31 fasc. 6, (1998)
811---842.
\smallskip

\item{[B.V]}Barlet~D., Varouchas~J.  Fonctions holomorphes sur l'espace des cycles, Bulletin de la Soci\'et\'e Math\'ematique de France ~117, (1989) 329---341\smallskip

\item{[E]}Elzein~F., Complexe dualisant et applications \`a la classe fondamentale d'un cycle, Bull.Soc.Math.France Mem~58 (1978).
\smallskip

\item{[C]}Cassa~A.,The Topology of the Space of Positive Analytic
Cycles. Annali di Math.Pura.Appl, tome 112,~1 (1977) p1---12.
\smallskip

\item{[Fr]}Frsich~J., Points de platitude d'un morphisme d'espace  analytique complexe. IInv.Math~4, (1967), p.118---138
\smallskip

\item{[Fu]}Fujiki~A., Closedeness of the Douady Spaces of compact
K\"ahler Spaces. Publ.Rims, Kyoto-Univ,tome 14, (1978) p1---52.
\smallskip

\item{[Go]}Golovin~V.D.,On the homology theory of analytic sheaves, Math.USSR.Izvestija~16,(2)  (1981),p239---260.
\smallskip

\item{[H]}Harvey~F.Reese.,Integral formulae connected by Dolbeault's isomorphism. Rice.Univ.Studies~56,(n2), (1970) p77---97.
\smallskip

\item{[H.S]} H\"ubl~R, Sastry~P., Regular differential forms and relative duality. American Journal of Mathematics~115, (1993), p.749---787.
\smallskip

\item{[Ju]}Jurchescu~M.,On the canonical topology of an analytic algebra of an analytic module, Bull.Soc.Math>France,~93, (1965), p 129---153.
\smallskip

\item{[K1]}Kaddar~M.Classe fondamentale relative en cohomologie de Deligne et application. Math.Ann,~306,(1986),p285---322.
\smallskip

\item{[K2]}Kaddar~M., Int\'egration d'ordre sup\'erieure sur les cycles en g\'eom\'etrie analytique complexe, Ann. Sc. Norm. Pisa Cl.Sci(4), tome  29, (2000)
187---263.
\smallskip

\item{[Ke]}Kersken~M.,Der Residuencomplex in der lokalen algebraischen und analytischen Geometrie, Math.Ann~265,(1983) 423---455.
\smallskip

\item{[K.W]}Kunz~E, Waldi~R.,  Regular differential forms, Contemporary Mathematics~79, (1988), Amer.Math.Soc., Providence.
\smallskip

\item{[Kl]}Kleiman, Steven~L., Relative duality for quasicoherent sheaves, Compositio.Math.~41 (1980),(1) 39---60.
\smallskip

\item{[L]}Lipman~J.,  Dualizing  sheaves, differentials and residues on algebraic varieties,  Asterisque~117, (1984).
\smallskip

\item{[Ma]}Mathieu~D., Universal reparametrizationof family of
analytic cycles: a new approach to meromorphicc equivalence
relations. Analles de l'Inst.Fourier~50 (44),(2000) p1115---1189.
\smallskip

\item{[Mo]}Mochizuki~N., Quasi-normal analytic spaces, Proc.Japan Acad.,~48, (1972), p.181---185.
\smallskip

\item{[M.F.K]}Mumford~D, Fogarty~J, Kirwan~F., Geometric invariant theory. Third edition. Ergebnisse der Mathematik und ihrer Grenzgebiete(2),~34, Sprinmger-Verlag, Berlin,(1994).
\smallskip

\item{[R]}Reiffen~H.J., Riemmansche Hebbartkeitss\"atze f\"r Kohomologieklassen mit kompactem Tr\"ager, Math.Ann~164, (1966) 272---279.
\smallskip

\item{[RRV]}Ramis~J.P, Ruget~G, Verdier~J.L., Dualit\'e relative en g\'eom\'etrie analytique complexe, Invent.Math ~13 (1971) 261---283.
\smallskip

\item{[Ry]}Rydh~D., Families of zero cycles and divided power I. Arxiv 0803.0618v1[math AG], 5mai 2008.
\smallskip

\item{[Si]}Siebert~B.,Fiber cycles of holomorphic maps I-II, Math.Ann~296 et 300,(1993,1994),p269----283, p243---271.
\smallskip

\item{[V]}Verdier~J.L. Base change for twisted inverse image of coherent sheaves, Internat.Colloq, Tata Inst.Fund.Res, Bombay,(1968) 393---408.
\smallskip

\item{[Va]}Varouchas~J,.K\"ahler spaces and proper open morphisms. Math.Ann ~283 (1), (1989), p.13---52. 

\bigskip\smallskip\noindent
M. Kaddar,\par\noindent
Institut Elie Cartan, UMR 7502\par\noindent
Universit\'e Nancy 1, BP 239\par\noindent
Vandoeuvre-l\`es-Nancy Cedex. France.\smallskip\noindent
e-mail: kaddar@iecn.u-nancy.fr\end
 La
construction, de nature purement alg\'ebrique,   du faisceau des
formes  m\'eromorphes r\'eguli\`eres relatives, entreprise par  Kunz
et Waldi([K.W]), se transporte \`a la cat\'egorie analytique
complexe pour un morphisme plat gr\^ace \`a l'important et
cons\'equent travail de Kersken ([Ke]). Ainsi, \`a tout morphisme
plat, \`a fibres de dimension pure $n$, d'espaces analytiques
complexes ( r\'eduits) $ g:Y\rightarrow T$, est canoniquement
associ\'e un unique faisceau (\`a isomorphisme canonique pr\`es)
${\cal O}_{X}$-coh\'erent  $\tilde{\omega}^{n}_{Y/T}$ enti\`erement
caract\'eris\'e par la propri\'et\'e de la trace
relative.\par\noindent Une fa\c con directe pour nous a \'et\'e de
v\'erifier point par point les d\'etails de la construction de
Kersken pour nous appercevoir qu'elle n'est pas le privil\`ege de la
platitude mais de la platitude g\'eom\'etrique analytique. Mais il
est moins douloureux de construire un tel faisceau pour un morphisme
analytiquement  g\'eom\'etriquement plat gr\^ace aux diff\'erentes
caract\'erisations que l'on donne de ces morphismes.
 En d'autre termes, on peut \'enoncer le
Plus pr\'ecisemment, pour tout $s_{0}$ dans $S$, il
existe voisinage ouvert $S_{0}$ de $s_{0}$ dans $S$ et une famille
d'\'ecailles ou de cartes  $(E_{\alpha}:=(U_{\alpha}, B_{\alpha},
\sigma_{\alpha}))_{\alpha\in A}$ $S_{0}$- adapt\'ees  \`a
$X_{s_{0}}$ de sorte que, pour tout $\alpha\in A$,  on ait une
installation (locale) ($\clubsuit$)
$$\xymatrix{X_{\alpha}\ar[d]^{\pi_{\alpha}} \ar [dr]_{f_{\alpha} }\ar [r]^{\sigma_{\alpha}}
 & S_{0}\times U_{\alpha}\times B_{\alpha} \ar[d]^{p_{\alpha}}\\
 S_{0}&\ar[l]^{q_{\alpha}} S_{0}\times U_{\alpha}}$$
dans laquelle $f_{\alpha}$ est un morphisme fini et surjectif,
$\sigma_{\alpha}$ est un plongement local, $p_{\alpha}$ et
$q_{\alpha}$ \'etant les projections canoniques. De telles cartes, que l'on peut choisir centr\'ees par un rapport \`a un point privil\'egi\'e sur la fibre, existent toujours (cf [A.S]).\par\noindent
 Par
construction m\^eme de `` l'int\'egration locale `` sur les cycles,
les morphismes ${f_{\alpha}}_{*}\Omega^{n}_{X_{\alpha}/S}\rightarrow
\Omega^{n}_{Y_{\alpha}/S}$, obtenus pr\'ec\'edemment, sont  de
v\'eritables morphismes traces et v\'erifient, de fa\c con
\'evidente, les propri\'et\'es \'enonc\'ees dans le lemme {\bf
(2.4)}.\par\noindent
  Mais dans chacune de ces situations
locales, [B4] donne une classe fondamentale relative ${\cal
C}_{X_{\alpha}/S}$ induisant un morphisme cup-produit
$\ds{\wedge{\cal C}_{X_{\alpha}/S}:
{\sigma_{\alpha}}_{*}\Omega^{n}_{X_{\alpha}/S}\rightarrow {\cal
H}^{p}_{|X_{\alpha}|}(\Omega^{n+p}_{Z_{\alpha}/S})}$. Utilisant,
alors, le morphisme canonique $\ds{{\cal
E}xt^{p}({\sigma_{\alpha}}_{*}{\cal O}_{X_{\alpha}},
\Omega^{n+p}_{Z_{\alpha}/S})\rightarrow{\cal
H}^{p}_{|X_{\alpha}|}(\Omega^{n+p}_{Z_{\alpha}/S})}$, on en d\'eduit
le diagramme\par
\centerline{$\xymatrix{{\sigma_{\alpha}}_{*}\Omega^{n}_{X_{\alpha}/S}\ar[rr]^{{\cal
C}_{\pi_{\alpha}}}\ar[rd]_{{\cal C}_{{X_{\alpha}}/S}}&&{\cal
E}xt^{p}({\sigma_{\alpha}}_{*}{\cal O}_{X_{\alpha}},
\Omega^{n+p}_{Z_{\alpha}/S})\ar[ld]\\
&{\cal
H}^{p}_{|X_{\alpha}|}(\Omega^{n+p}_{Z_{\alpha}/S})&}$}\noindent qui
est commutatif puisque ${\cal C}_{X_{\alpha}/S}$ \'etant annul\'ee,
au moins g\'en\'eriquement sur $S$, par tout id\'eal  de
d\'efinition de $X_{\alpha}$ dans $Z_{\alpha}$, ${\cal
C}_{\pi_{\alpha}}$ en est un rel\`evement naturel.\par\noindent
D'apr\`es le th\'eor\`eme 1, le faisceau $\omega^{n}_{\pi}$ muni
$\pi$ d'une fl\`eche globale ${\rm I}\!{\rm
R}^{n}{\pi}_{!}\omega^{n}_{\pi}\rightarrow {\cal O}_{S}$ qui, par
localisation, donne ${\rm I}\!{\rm
R}^{n}{\pi_{\alpha}}_{!}\omega^{n}_{\pi_{\alpha}}\rightarrow {\cal
O}_{S}$. Mais la donn\'ee de cette derni\`ere est strictement
\'equivalente \`a la donn\'ee  d'un  morphisme
${f_{\alpha}}_{*}{\cal E}xt^{p}({\sigma_{\alpha}}_{*}{\cal
O}_{X_{\alpha}}, \Omega^{n+p}_{Z_{\alpha}/S})\rightarrow
\Omega^{n+p}_{Y_{\alpha}/S}$ qui n'a, \`a priori, aucune raison
d'\^etre un morphisme trace!\par\noindent Il se trouve que dans la
situation particuli\`ere envisag\'ee, c'est effectivement le cas.
Pour s'en assurer, il suffit simplement de se rappeler que la classe
fondamentale relative
 ${\cal C}_{X_{\alpha}/S}$ induit d\'ej\`a, pour tout entier $k\leq n$,  un morphisme trace
$${f_{\alpha}}_{*}{\cal H}^{p}_{|X_{\alpha}|}(\Omega^{k+p}_{Z_{\alpha}/S})
\rightarrow f_{*}{\cal
H}^{p}_{{Y_{\alpha}}-propre}(\Omega^{k+p}_{Z/S})\rightarrow
\Omega^{k}_{Y_{\alpha}/S}$$ pouvant \^etre vu comme un morphisme
{\it r\'esidu} ou  d'{\it int\'egration} de classes de cohomologie
\`a support $|X_{\alpha}|$ (qui est  $Y_{\alpha}$- propre puisque
fini sur $Y_{\alpha}$)  sur les fibres $p$-\'equidimensionnelles du
morphisme lisse $q'_{\alpha}:Z_{\alpha}\rightarrow Y_{\alpha}$. Les
fl\`eches ainsi d\'efinies prolonge naturellement les traces
usuelles
$\ds{{f_{\alpha}}_{*}{f_{\alpha}}^{*}\Omega^{k}_{Y_{\alpha}/S}\rightarrow
\Omega^{k}_{Y_{\alpha}/S}}$ et rendent commutatif le diagramme
$$\xymatrix{{f_{\alpha}}_{*}\Omega^{k}_{X/S}\ar[rd]^{{\cal T}^{k}_{{f_{\alpha}}}}
\ar[r]^{\!\!\!\!\!\!\!\!\!\wedge {\cal C}_{\pi_{\alpha}}}&{f_{\alpha}}_{*}{\cal H}^{p}_{|X_{\alpha}|}(\Omega^{p+k}_{Z_{\alpha}/S})\ar[d]^{\phi}\\
{f_{\alpha}}_{*}{f_{\alpha}}^{*}\Omega^{k}_{Y_{\alpha}/S}\ar[u]\ar[r]&\Omega^{k}_{Y_{\alpha}/S}}$$
avec stabilit\'e par changement de base entre espaces complexes
r\'eduits, compatibilit\'e  aux localisations sur $X$ et \`a
l'additivit\'e des pond\'erations A la vue de tout ceci  et le lien
intime entre la trace d'un morphisme fini et le morphisme r\'esidu
(cf p.49), il apparait que la  construction  s'adapte   \`a
n'importe quel morphisme d'alg\`ebres analytiques locales provenant
de morphisme universellement \'equidimensionnel et  ayant la vertu
de produire de v\'eritables morphismes traces relatives dans toute
installation locale  (i.e factorisation locale en morphisme fini $f$
suivi d'un lisse) $f_{*}\Omega^{\bullet}_{X/S}\rightarrow
\Omega^{\bullet}_{Y/S}$. Il parait alors \'evident, \`a la lumi\`ere
de la proposition {\bf(2.1)} donnant une des caract\'eristaion  des
morphismes analytiquement g\'eom\'etriquement plats, que [Ke]
s'applique mutatis-mutandis \`a ce cadre en d\'ecr\'etant m\^eme
qu'{\it{un morphisme d'alg\`ebres analytiques locales  $\phi:
\phi:{\cal P}\rightarrow {\cal A}$  est analytiquement
g\'eom\'etriquement plat si et seulement si pour toute ${\cal
P}$-param\'etrisation locale $f:{\cal R}\rightarrow {\cal A}$, il
existe un morphisme trace de ${\cal R}$-modules ${\cal T}r:{\cal
A}\rightarrow {\cal R}$}}.\par\noindent Comme notre probl\`eme est
de  nature purement locale sur $X$ et $S$, on se ram\`ene \`a l'
\'etude d'une situation  ``germifi\'ee'' et donc \`a  des
consid\'erations d'alg\`ebres analytiques locales nous permettant
d'utiliser [Ke]  que l'on globalise gr\^ace au th\'eor\`eme de
noeth\'eriannit\'e de Frisch-Grothendieck  en se rappelant que la
notion de  sch\'ema noeth\'erien est l'analogue de celle de compact
de Stein dans le cadre analytique complexe.\par\noindent Signalons,
enfin, que toutes les alg\`ebres analytiques peuvent \^etre munies
d'une topologie canonique (cf [Ju]) La diff\'erentielle dont on peut
munir $\omega^{\bullet}_{\pi}= {\cal H}om(\Omega^{n-\bullet}_{X/S},
\omega^{n}_{\pi})$ ne peut pas \^etre naivement d\'eduite de la
diff\'erentielle ext\'erieure relative usuelle puisque, relativement
\`a un plongement local de $X$ dans $Z$ lisse sur $S$, la
diff\'erentielle n'est pas ${\cal O}_{Z}$-lin\'eaire. D'ailleurs la
construction de Kersken [Ke] dans le cas d'un morphisme d'alg\`ebres
analytiques plates montre bien que sur un objet de ce type la
diff\'erentielle n'est jamais aussi naive ! Soient  $n$ un entier
naturel  et  $\pi:X\rightarrow S$ un morphisme universellement
$n$-\'equidimensionnel d'espaces analytiques complexes r\'eduits.
Soient  $X_{0}$ (resp. $S_{0}$) l' ouvert dense de $X$ (resp. $S$)
sur lequel la restriction  $\pi_{0}$ de $\pi$ est plate sur $S_{0}$
et $\omega^{n}_{X_{0}/S_{0}}$ le faisceau coh\'erent caract\'eris\'e
par la propri\'et\'e de la trace relative de [Ke] . Alors
\par\noindent $\omega^{n}_{\pi}$ est un prolongement
 canonique et coh\'erent  du faisceau $\omega^{n}_{X_{0}/S_{0}}$  si
et seulement si $\pi$ est analytiquement g\'eom\'etriquement plat.
et caract\'eris\'e par la propri\'et\'e de la trace relative disant
qu'une section $\xi$ de $j_{*}j^{*}\Omega^{k}_{X/S}$ sur un ouvert
de Stein de $X$ d\'efinit une section du faisceau
$\Lambda^{k}_{X/S}$ si et seulement si, pour toute $n-k$-forme
holomorphe relative $\alpha$, $\xi\wedge \alpha$ v\'erifie la
propri\'et\'e de la trace relative pr\'ec\'edemment d\'efinie.
l'isomorphisme
$$\Lambda^{k}_{\sigma}\simeq{\cal H}om_{{\cal
O}_{X}}(\Omega^{n-k}_{X/S}, \Lambda^{n}_{\sigma})$$ puisque l'on
dispose d\'ej\`a d'une fl\`eche injective
$\Lambda^{k}_{\sigma}\rightarrow{\cal H}om_{{\cal
O}_{X}}(\Omega^{n-k}_{X/S}, \Lambda^{n}_{\sigma})$ et que le
faisceau ${\cal H}om_{{\cal O}_{X}}(\Omega^{n-k}_{X/S},
\Lambda^{n}_{\sigma})$ v\'erifie manifestement la propri\'et\'e de
la trace relative

\par\noindent On peut,
bien s\^ur, exploiter le diagramme
$$\xymatrix{0\ar[r]&\Lambda^{k}\ar[r]&j_{*}j^{*}\Omega^{k}_{X/S}\ar[d]\ar[rd]\\
0\ar[r]&\omega^{k}_{\pi}\ar[r]&j_{*}j^{*}\omega^{k}_{\pi}\ar[r]&{\cal
H}^{1}_{\Sigma}(\omega^{k}_{\pi})}$$ et utiliser la m\^eme ligne de
d\'emonstration que dans la {\it proposition 7}.
$\,\,\,\blacksquare$\smallskip\smallskip\noindent
{\bf 3.3.3.
Remarques.}\par\noindent
{\bf(i)} Evidemment ces faisceaux sont canoniquement
isomorphes aux faisceaux $\omega^{k}_{\pi}$ donn\'es par le {\it
corollaire (1.1)} du {\it th\'eor\`eme 1}.\par\noindent
{\bf(ii)}  Rappelons que pour tout morphisme
$f:X\rightarrow Y$  fini, surjectif et analytiquement
g\'eom\'etriquement plat sur $Y$ lisse sur $S$ et dimension relative
$n$, on dispose d' un  morphisme r\'esidu (cf \S 2 de [KI])  (ou
int\'egration sur les fibres de la projection $Z\rightarrow Y$ de
classes de cohomologie de ${\cal H}^{p}_{X}(\Omega^{n+p}_{Z/S})$)
$${\goth R}^{X}_{Y}:f_{*}{\cal H}^{p}_{X}(\Omega^{n+p}_{Z/S})\rightarrow \Omega^{n}_{Y/S}$$
rendant commutatif le diagramme
$$\xymatrix{f_{*}{\cal E}xt^{p}({\cal O}_{X}, \Omega^{n+p}_{Z/S})\ar[rr]\ar[rd]&&f_{*}{\cal H}^{p}_{X}(\Omega^{n+p}_{Z/S})\ar[ld]\\
& \Omega^{n}_{Y/S}&}$$ dans lequel la fl\`eche ${\goth
T}^{n}_{f}:f_{*}\Lambda^{n}\rightarrow \Omega^{n}_{Y/S}$
 est un v\'eritable morphisme trace (cf le {\it th\'eor\`eme 2}).
 En notant $\tilde{\Sigma}:=f(\Sigma)$ et
$\tilde{j}:Y-{\tilde{\Sigma}}\rightarrow Y$ l'inclusion naturelle,
on a un diagramme commutatif  de suites exactes courtes
$$\xymatrix{0\ar[r]&f_{*}\Lambda^{n}\ar[r]\ar[d]_{\goth{T}^{n}_{\pi}}&f_{*}j_{*}j^{*}(\Omega^{n}_{X/S})\ar[r]\ar[d]&
f_{*}{\cal H}^{1}_{\Sigma}(\Lambda^{n}) \ar[d]\\
0\ar[r]&\Omega^{n}_{Y/S}\ar[r]&{\tilde{j}}_{*}{\tilde{j}}^{*}(\Omega^{n}_{Y/S})\ar[r]&{\cal
H}^{1}_{\tilde\Sigma}(\Omega^{n}_{Y/S})}$$ dans lequel les fl\`eches
verticales sont des morphismes traces.
 On peut,
bien s\^ur, exploiter le diagramme
$$\xymatrix{0\ar[r]&\Lambda^{k}\ar[r]&j_{*}j^{*}\Omega^{k}_{X/S}\ar[d]\ar[rd]\\
0\ar[r]&\omega^{k}_{\pi}\ar[r]&j_{*}j^{*}\omega^{k}_{\pi}\ar[r]&{\cal
H}^{1}_{\Sigma}(\omega^{k}_{\pi})}$$ et utiliser la m\^eme ligne de
d\'emonstration que dans la {\it proposition 7}.
 Pour \^etre plus pr\'ecis, dans le diagramme
$$\xymatrix{&&X_{2}\times X_{1}\times S\ar[ld]_{p_{2}}\ar[d]^{p}\ar[rd]^{p_{1}}&&\\
&X_{1}\times S\ar[d]&X_{2}\times S\ar[l]_{\pi_{2}\times id}\ar[d]&X_{2}\times X_{1}\ar[l]_{id\times \pi_{1}}\ar[d]&\\
&X_{1}\ar[r]_{\pi_{1}}&S&Y\ar[l]^{\pi_{1}}}$$ Alors la pond\'eration
de $\pi$ est donn\'ee par le cycle $\goth{ X}$ de $X_{2}\times S$
v\'erifiant
$$p^{*}({\goth{X}}) = p^{*}_{1}({\goth X}_{1})\bullet p^{*}_{2}({\goth X}_{2})$$
{\bf 2.1.6. Remarque.}\par\noindent
Si l'on applique [K2] \`a $\pi_{1}$ et $\pi_{2}$, on obtient les
fl\`eches ${\rm I}\!{\rm R}^{n_{2}}{\pi_{2}}_{!}(
\Omega^{n_{2}+n_{1}}_{X_{2}/S})\rightarrow {\cal L}^{n_{1}}_{X_{1}}$
puis ${\rm I}\!{\rm R}^{n_{1}}{\pi_{1}}_{!}{\cal
L}^{n_{1}}_{X_{1}}\rightarrow {\cal L}^{0}_{S}$.\par\noindent Mais
un examen d\'etaill\'e de la premi\`ere int\'egration montre,
gr\^ace au {\it corollaire 4} (p.35), que son image est en fait un
sous faisceau de ${\cal L}^{n_{1}}_{X_{1}}$ dont les sections sont
caract\'eris\'ees par le fait que leurs traces, dans toute
param\'etrisation locale $S$-relative $f:X_{1}\rightarrow S\times
U$, sont des sections du faisceau ${\cal
O}^{c}_{S}{\widehat\otimes}\Omega^{n_{1}}_{U}$. Il est, alors facile
de v\'erifier que la compos\'ee  est \`a valeurs dans le faisceau
${\cal O}^{c}_{S}$. Cela montre que la compos\'ee de deux morphismes
analytiquement g\'eom\'etriquement plats est toujours contin\^ument
g\'eom\'etriquement plate.

%% file: zineb.tex

	%
	%
	%
	\newfam\symboles
	

%
 %

 %

	\font\ninerm=cmr9
\font\eightrm=cmr8
 \font\ninei=cmti9

%


	%
	%

\def\flead{\leaders\hbox to 5pt {\hss.}\hfill}
\def\somr#1|#2|{#1\flead \rlap{\hbox to 25pt{\hfill #2}}\par
\goodbreak}
\def\ssomr#1|#2|{\qquad #1\flead \rlap{\hbox to 25pt {\hfill #2}}\par
\goodbreak}

 \magnification=1200 
\pageno=1
	\baselineskip7mm
	\hsize=130mm \vsize=200mm 
	\hoffset=35mm \voffset=10mm
	\hfuzz=50pt \vfuzz=10pt
	
	\headline={\hfill-- \tenrm\folio\hbox{ }--\hfill}
	\footline={\hfill}

\newtoks\auteurcourant      \auteurcourant={\hfil}
\newtoks\titrecourant       \titrecourant={\hfil}

\newtoks\hautpagetitre      \hautpagetitre={\hfil}
\newtoks\baspagetitre       \baspagetitre={\hfil}

\newtoks\hautpagegauche   \newtoks\hautpagedroite 
  
\hautpagegauche={\tenrm\rlap{\folio}\tenit\hfil\the\auteurcourant\hfil}
\hautpagedroite={\tenit\hfil\the\titrecourant\hfil\tenrm\llap{\folio}}

\newtoks\baspagegauche      \baspagegauche={\hfil} 
\newtoks\baspagedroite      \baspagedroite={\hfil}

\newif\ifpagetitre          \pagetitretrue
\nopagenumbers  

\def\nopagenumbers{\def\folio{\hfil}}  

\headline={\ifpagetitre\the\hautpagetitre
            \else\ifodd\pageno\the\hautpagedroite
             \else\the\hautpagegauche
              \fi\fi}

\footline={\ifpagetitre\the\baspagetitre\else
            \ifodd\pageno\the\baspagedroite
             \else\the\baspagegauche
              \fi\fi
               \global\pagetitrefalse}


\def\raggedbottom{\topskip 10pt plus 36pt\r@ggedbottomtrue}

	%
	%
	%
	%
	%
	%
	
	%
        \def\page#1{\leaders\hbox to 5mm{\hfil\hbox{\ninerm.}\hfil}
        \hfill\rlap{\hbox to 5mm{\hfill#1}}\par}

%
  \def\Th#1{\medskip\goodbreak \nod{\bf  \underbar{Th\'eor\`eme} #1.}
	\quad\nobreak \sl }
%
  

%
\def\Prop#1{\medskip\goodbreak \nod{\bf
	 Proposition #1.} \quad\nobreak \sl }
	
%
	\def\Lemme#1{\medskip\goodbreak \nod{\bf  Lemme #1.}
	\quad\nobreak \sl }
%

%

%

%

%

%
	
%

%
        
%

%
\def\Cor#1{\medskip\goodbreak \nod{\bf
	Corollaire#1.} \quad\nobreak \sl }
%
	\def\dem{\medskip\goodbreak \nod \ninei D\'emonstration. \rm}

  \def\ds{\displaystyle}
	\def\ss{\scriptstyle}
	
	\def\nod{\noindent}
	
	\def\hb{\hbox}

	\def\dbar{d\!\!\hbox to 4.5pt{\hfill\vrule height 5.5pt
	depth -5.3pt width 3.5pt}}

	\def\tas#1_#2^#3{\mathrel{\mathop{\kern 0pt#1}
					\limits_{#2}^{#3}  }}

	\def\up#1{\raise \hbox{\ss \rm #1}}


	\def\hfl#1#2{\smash{\mathop{\hbox to 10mm{\rightarrowfill}}
	\limits^{\ss#1}_{\ss#2}}}
        \def\hfla#1#2{\smash{\mathop{\hbox to 6mm{\rightarrowfill}}
	\limits^{\ss#1}_{\ss#2}}}
	\def\hflo#1#2{\smash{\mathop{\hbox to 10mm{\leftarrowfill}}
	\limits^{\ss#1}_{\ss#2}   }}

	\def\pile#1#2{\smash{\mathop{\hbox to 20mm{}}
	\limits^{\hb{\rm#1}}_{\hb{\rm#2}} }}
	\def\pila#1#2#3{\smash{\mathop{\hbox to 20mm{#2}}
	\limits^{\hb{\rm#1}}_{\hb{\rm#3}}   }}